\documentclass{article}
\usepackage[utf8]{inputenc}
\usepackage[T1]{fontenc}
\usepackage[english]{babel}
\usepackage{amsmath}
\usepackage{amsthm,amssymb}
\usepackage{authblk}
\usepackage{amssymb}
\usepackage{xspace}
\usepackage{hyperref}
\usepackage{verbatim}
\usepackage{tikz}
\usepackage[procnumbered,linesnumbered,ruled,vlined]{algorithm2e}
\usepackage{tabularray}
\usepackage[textsize=footnotesize,color=green!40]{todonotes}
\usepackage{stmaryrd}
\usepackage{comment}
\usepackage{subcaption}
\usepackage{thmtools}
\usepackage{cleveref}
\usepackage[left= 3 cm,right=3 cm,top=3 cm,bottom=3 cm]{geometry}

\usetikzlibrary{shapes.misc}
\usetikzlibrary{calc}
\tikzstyle{noeud}=[circle,inner sep=2, minimum size =3 pt, line width = 1pt, draw=black, fill=white]
\newcommand{\convexpath}[2]{
	[   
	create hullnodes/.code={
		\global\edef\namelist{#1}
		\foreach [count=\counter] \nodename in \namelist {
			\global\edef\numberofnodes{\counter}
			\node at (\nodename) [draw=none,name=hullnode\counter] {};
		}
		\node at (hullnode\numberofnodes) [name=hullnode0,draw=none] {};
		\pgfmathtruncatemacro\lastnumber{\numberofnodes+1}
		\node at (hullnode1) [name=hullnode\lastnumber,draw=none] {};
	},
	create hullnodes
	]
	($(hullnode1)!#2!-90:(hullnode0)$)
	\foreach [
	evaluate=\currentnode as \previousnode using \currentnode-1,
	evaluate=\currentnode as \nextnode using \currentnode+1
	] \currentnode in {1,...,\numberofnodes} {
		-- ($(hullnode\currentnode)!#2!-90:(hullnode\previousnode)$)
		let \p1 = ($(hullnode\currentnode)!#2!-90:(hullnode\previousnode) - (hullnode\currentnode)$),
		\n1 = {atan2(\y1,\x1)},
		\p2 = ($(hullnode\currentnode)!#2!90:(hullnode\nextnode) - (hullnode\currentnode)$),
		\n2 = {atan2(\y2,\x2)},
		\n{delta} = {-Mod(\n1-\n2,360)}
		in 
		{arc [start angle=\n1, delta angle=\n{delta}, radius=#2]}
	}
	-- cycle
}

\newcommand{\outcomeP}{$\mathcal{P}$\xspace}
\newcommand{\outcomeN}{$\mathcal{N}$\xspace}
\newcommand{\gr}{\mathcal{G}}

\DeclareMathOperator{\mex}{mex}
\DeclareMathOperator{\opt}{opt}
\DeclareMathOperator{\lt}{lt}
\newcommand{\AGAG}{\textsc{Closed Geodetic Game}\xspace}
\newcommand{\GAG}{\textsc{Geodetic Game}\xspace}
\newcommand{\lovesme}{\textsc{She Loves Me-She Loves Me Not}\xspace}
\newcommand{\grFBG}{GrundyBlock}
\newcommand{\grFBV}{GrundyBlockWithVertex}
\newcommand{\grFCG}{GrundyCactus}
\newcommand{\grFCVA}{GrundyCactusOneVertex}
\newcommand{\grFCVB}{GrundyCactusTwoVertices}

\newtheorem{theorem}{Theorem}

\newtheorem{proposition}[theorem]{Proposition}

\newtheorem{lemma}[theorem]{Lemma}

\newtheorem{claim}{Claim}[theorem]
\newtheorem{notation}[theorem]{Notation}

\title{The Closed Geodetic Game: algorithms and strategies\thanks{A short version of this paper was presented at IWOCA 2025~\cite{dailly2024closed2}.} \thanks{The research of the first author was supported by the International Research Center "Innovation Transportation and Production Systems" of the I-SITE CAP 20-25. The research of the second author was partially supported by ANR project GRALMECO (ANR-21-CE48-0004-01). The research of the third author was supported by the "Doctoral Fellowship in India for ASEAN (DIA:2020-25)".}}
\author[1,2]{Antoine Dailly}
\author[1,3,4]{Harmender Gahlawat}
\author[5,6]{Zin Mar Myint}
\affil[1]{Université Clermont Auvergne, CNRS, Mines de Saint-Étienne, Clermont-Auvergne-INP, LIMOS, 63000 Clermont-Ferrand, France}
\affil[2]{Université Clermont Auvergne, INRAE, UR TSCF, 63000, Clermont-Ferrand, France}
\affil[3]{Laboratoire G-SCOP (CNRS, Univ. Grenoble Alpes), Grenoble, France}
\affil[4]{Ben-Gurion University of the Negev, Beersheba, Israel}
\affil[5]{Indian Institute of Technology Dharwad, India}
\affil[6]{Polytechnic University (Kyaing Tong), Myanmar}

\date{}

\begin{document}

\maketitle

\begin{abstract}
	The \emph{geodetic closure} of a set $S$ of vertices of a graph $G$ is the set of all vertices in shortest paths between pairs of vertices of $S$.
	A set $S$ of vertices in a graph is \emph{geodetic} if its geodetic closure contains all the vertices of the graph.
	
	Several authors have studied variants of games around constructing geodetic sets. The most studied of those, \GAG, was introduced by Harary in 1984 and developed by Buckley and Harary in 1985. It is an \emph{achievement} game: both players construct together a geodetic set by alternately adding vertices to the set, the winner being the one who plays last. However, this version of the game allows the players to select vertices that already are in the geodetic closure of the current set.
	
	We study the more natural variant, called \AGAG, where the players alternate adding to $S$ vertices that are not in the geodetic closure of $S$. This variant was also introduced in another, less-noticed, paper by Buckley and Harary in 1985, and only studied since then in the context of trees by Araujo \emph{et al.} in 2024.
	We provide a full characterization of the Sprague-Grundy values (equivalence values for games) of graph classes such as paths and cycles, of the outcomes of the Cartesian product of several graphs in function of their individual outcomes, and give polynomial-time algorithms to determine the Sprague-Grundy values in cactus and block graphs.
	
\end{abstract}

\section{Introduction}

\subsection{A history of geodetic games}

Let $G(V,E)$ be a simple, undirected graph with vertex set $V$ and edge set $E$. Given two vertices $u,v \in V$, let $\mathcal{I}(u,v)$ denote the set of vertices that lie on some shortest path with $u$ and $v$ as endpoints.  Given a subset $S \subseteq V$ of vertices, its \emph{geodetic closure}, denoted by $(S)$, is defined as the set of vertices in all shortest paths between every pair of vertices of $S$. More formally, $(S) = \bigcup_{u,v\in S} \mathcal{I}(u,v)$.

\begin{figure}[h]
	\centering
	\begin{subfigure}[b]{0.2\textwidth}
		\centering
		\scalebox{1}{
			\begin{tikzpicture}[scale=0.85, transform shape]
				\node[noeud,fill=black] (0) at (0,0) {};
				\node[noeud,cross out] (1) at (1,0) {};
				\node[noeud, cross out] (2) at (2,0) {};
				\node[noeud,fill=black] (3) at (3,0) {};
				\draw (0)to(1)to(2)to(3);
				\node (empty) at (0,-1) {};
			\end{tikzpicture}
		}
	\end{subfigure}\hfill
	\begin{subfigure}[b]{0.4\textwidth}
		\centering
		\scalebox{1}{
			\begin{tikzpicture} [scale=0.85, transform shape]
				\node[noeud, fill=black] (0) at (-1,0) {};
				\node[noeud, cross out] (1) at (0,1) {};
				\node[noeud, cross out] (2) at (0,-1) {};
				\node[noeud] (3) at (1,0) {};
				\node[noeud, cross out] (5) at (2,0) {};
				\node[noeud, cross out] (6) at (3,1) {};
				\node[noeud, cross out] (7) at (3,-1) {};
				\node[noeud, cross out] (8) at (4,0) {};
				\node[noeud, fill=black] (9) at (3,0) {};
				\draw[fill=none, dotted, very thick](1,0) circle (0.2) node {};
				\draw (0)to(1)to(3)to(2)to(0);
				\draw (5)to(6)to(8)to(7)to(5);
				\draw (5)to(9)to(8);
				\draw (1)to(6);
				\draw (2)to(7);
				\draw (0)to(3);
			\end{tikzpicture}
		}
	\end{subfigure}
	\begin{subfigure}[b]{0.3\textwidth}
		\centering
		\scalebox{1}{
			\begin{tikzpicture} [scale=0.85, transform shape]
				\node[noeud] (0) at (0.2,0.2) {};
				\node[noeud] (1) at (1,0.7) {};
				\node[noeud, cross out] (2) at (1.8,0.2) {};
				\node[noeud] (3) at (1.5,-0.8) {};
				\node[noeud,fill=black] (4) at (0.5,-0.8) {};
				\draw (0)to(3)to(1)to(4)to(2)to(0);
				\node[noeud,fill=black] (5) at (-0.5,0.3) {};
				\node[noeud, cross out] (6) at (1,1.5) {};
				\node[noeud,fill=black] (7) at (2.5,0.3) {};
				\node[noeud] (8) at (1.8,-1.3) {};
				\node[noeud, cross out] (9) at (0.2,-1.3) {};
				\draw (5)to(6)to(7)to(8)to(9)to(5);
				\draw[fill=none, dotted, very thick](1.5,-0.8) circle (0.2) node {};
				\draw (1)to(6);
				\draw (2)to(7);
				\draw (0)to(5);
				\draw (8)to(3);
				\draw (9)to(4);
			\end{tikzpicture}
		}
	\end{subfigure}
	\caption{Black vertices belong to the set $S$, crossed vertices belong to the geodetic closure $(S)$, white vertices are uncovered vertices yet. However, if we add the dotted vertices to $S$, then $V=(S)$, \emph{i.e.}, $S$ becomes a geodetic set.} 
	\label{fig-geodeticSets}
\end{figure}

A set $S$ of vertices of $G(V,E)$ is called a \emph{geodetic set of $G$} if $(S)=V$ (see \Cref{fig-geodeticSets} for examples). The problem of finding a geodetic set of minimum size was defined in the 1980s~\cite{harary1981convexity,buckley1988extremal,buckley1990distance,harary1993geodetic} and has been well-studied and is known to be NP-hard for general graphs~\cite{harary1993geodetic,atici2002computational}. The problem has been also well-studied on graph classes and displays quite an interesting behaviour. In particular, the problem stays NP-hard even for subcubic partial grids~\cite{chakraborty2023algorithms}, for interval graphs~\cite{chakraborty2020algorithms}, and admits a polynomial time algorithm for solid grids~\cite{chakraborty2023algorithms}, outerplanar graphs~\cite{mezzini2018polynomial}, and proper interval graphs~\cite{ekim2012computing}.

Interestingly, the game version of the geodetic number was introduced in parallel of its combinatorial counterpart. The first mention of a game based around geodetic sets was in 1984, by Harary~\cite{harary1984convexity}: two players alternate adding new vertices of a graph $G(V,E)$ to a set $S$, the game ends when $(S)=V$. Two winning conditions are proposed: in the \emph{achievement} setting, the player who plays the last move wins; while they lose in the \emph{avoidance} setting. Achievement and avoidance correspond to the \emph{normal} and \emph{misère} settings of combinatorial games, respectively. We will focus exclusively on achievement games.

With respect to combinatorial game theory, the first question to ask is, \textbf{given an input graph, decide which player wins}. This is called computing the \emph{outcome} of the game. Note that, since the game is \emph{impartial} (both players are indistinguishable), there are only two possible outcomes: either the first player wins (denoted by \outcomeN) or the second player wins (denoted by \outcomeP). Furthermore, when playing on disjoint connected components, each component is independent (selecting a vertex in a component has no influence whatsoever on the other components). If two components are \outcomeP, then, playing on their disjoint union is also \outcomeP (the second player always applies their winning strategy on the component on which the first player played a move); if one component is \outcomeP and the other is \outcomeN, then the first player wins by applying their winning strategy on the \outcomeN component, which becomes a \outcomeP component, and thus leaving the disjoint union of two \outcomeP components on which they are now the second player (so they win by the previous case). However, outcomes are not enough to solve the disjoint union of two \outcomeN components: we need to use \emph{Sprague-Grundy values}, which are equivalence classes for games (and can be mapped to nonnegative integers). More details will be given in \Cref{sec-games}. Hence, a stronger question is, \textbf{for a given input graph, compute its Sprague-Grundy value}.

\subsection{\GAG, \AGAG and other related variants}

In his seminal paper~\cite{harary1984convexity}, Harary discussed several possible rulesets for games based around the notion of geodesics. Two of those rulesets were then expanded in two 1985 papers by Buckley and Harary, and constitute the main variants. All games are played on a graph $G(V,E)$.

In the first paper~\cite{buckley1985geodetic}, the authors define \GAG. In this game, two players alternately add vertices to a set $S$, until $(S)=V$. They give the outcome of the game for complete graphs, cycles, (generalized) wheels, complete bipartite graphs and $n$-cubes. This line of research was expanded later in~\cite{necaskova1993note} which improved on the wheels' result, and in~\cite{haynes2003geodetic}, which studied complete multipartite graphs, hypercubes, and graphs with a unique minimum-size geodetic set (such as trees, split graphs or coronas).

In the second paper~\cite{buckley1985closed}, they introduce \AGAG. In this game, the players alternatively add to $S$ vertices that are not in $(S)$ as it is at the moment they are selecting it. They study the same graph classes as in~\cite{buckley1985geodetic}. Surprisingly, this more natural variant was left unstudied for decades, probably due to the difficulty of access of~\cite{buckley1985closed} compared to~\cite{buckley1985geodetic}. Recently, Araujo \emph{et al.}~\cite{araujo2024graph} studied several geodetic game variants, including \AGAG (which they called \emph{geodesic closed interval game}), for which they provided a linear-time algorithm computing the Sprague-Grundy value of a tree.

Note that \GAG and \AGAG can have different outcomes on even very simple graphs. One such example is illustrated in \Cref{fig-differencesActiveNonActive}. Since the players are allowed to pick covered vertices in \GAG, the parity of the graph (and of covered vertices) plays a bigger role in determining the outcome compared to \AGAG. Indeed, such is the case for Proposition~6 in~\cite{haynes2003geodetic}, where graphs with a unique minimum-size geodetic set are solved due to parity, since the simili-passing moves (adding to $S$ vertices that are in the current geodetic closure) are available.

\begin{figure}[h]
	\centering
	\begin{subfigure}[b]{0.30\textwidth}
		\centering
		\scalebox{1}{
			\begin{tikzpicture}
				\node[noeud] (0) at (0,1) {};
				\node[noeud,fill=black] (1) at (1,1) {};
				\node[noeud] (2) at (2,1) {};
				\node[noeud] (3) at (3,1) {};
				\node[noeud] (4) at (1,0) {};
				\draw (0)to(1)to(2)to(3);
				\draw (1)to(4);
				%\node (empty) at (0,-0.5) {};
			\end{tikzpicture}
		}
		\caption{For \GAG, the first player wins by picking the black vertex, ensuring that every vertex will be in $S$ (if the second player picks a leaf, they pick the degree~2 vertex, and conversely).}
		\label{fig-differencesGag}
	\end{subfigure}\hfill
	\begin{subfigure}[b]{0.69\textwidth}
		\centering
		\scalebox{0.65}{
			\begin{tikzpicture}
				\node (g1a1) at (0,0) {
					\begin{tikzpicture}
						\node[noeud,fill=black] (0) at (0,1) {};
						\node[noeud] (1) at (1,1) {};
						\node[noeud] (2) at (2,1) {};
						\node[noeud] (3) at (3,1) {};
						\node[noeud] (4) at (1,0) {};
						\draw (0)to(1)to(2)to(3);
						\draw (1)to(4);
					\end{tikzpicture}
				};
				\node (g1a2) at (4,0) {
					\begin{tikzpicture}
						\node[noeud] (0) at (0,1) {};
						\node[noeud] (1) at (1,1) {};
						\node[noeud,fill=black] (2) at (2,1) {};
						\node[noeud] (3) at (3,1) {};
						\node[noeud] (4) at (1,0) {};
						\draw (0)to(1)to(2)to(3);
						\draw (1)to(4);
					\end{tikzpicture}
				};
				\node (g1b) at (2,-2.5) {
					\begin{tikzpicture}
						\node[noeud,fill=black] (0) at (0,1) {};
						\node[noeud,cross out] (1) at (1,1) {};
						\node[noeud,fill=black] (2) at (2,1) {};
						\node[noeud] (3) at (3,1) {};
						\node[noeud] (4) at (1,0) {};
						\draw (0)to(1)to(2)to(3);
						\draw (1)to(4);
					\end{tikzpicture}
				};
				\node (g2a1) at (8,0) {
					\begin{tikzpicture}
						\node[noeud] (0) at (0,1) {};
						\node[noeud,fill=black] (1) at (1,1) {};
						\node[noeud] (2) at (2,1) {};
						\node[noeud] (3) at (3,1) {};
						\node[noeud] (4) at (1,0) {};
						\draw (0)to(1)to(2)to(3);
						\draw (1)to(4);
					\end{tikzpicture}
				};
				\node (g2a2) at (12,0) {
					\begin{tikzpicture}
						\node[noeud] (0) at (0,1) {};
						\node[noeud] (1) at (1,1) {};
						\node[noeud] (2) at (2,1) {};
						\node[noeud,fill=black] (3) at (3,1) {};
						\node[noeud] (4) at (1,0) {};
						\draw (0)to(1)to(2)to(3);
						\draw (1)to(4);
					\end{tikzpicture}
				};
				\node (g2b) at (10,-2.5) {
					\begin{tikzpicture}
						\node[noeud] (0) at (0,1) {};
						\node[noeud,fill=black] (1) at (1,1) {};
						\node[noeud,cross out] (2) at (2,1) {};
						\node[noeud,fill=black] (3) at (3,1) {};
						\node[noeud] (4) at (1,0) {};
						\draw (0)to(1)to(2)to(3);
						\draw (1)to(4);
					\end{tikzpicture}
				};
				\draw[->,double,line width=0.25mm] (g1a1)to(g1b);
				\draw[->,double,line width=0.25mm] (g1a2)to(g1b);
				\draw[->,double,line width=0.25mm] (g2a1)to(g2b);
				\draw[->,double,line width=0.25mm] (g2a2)to(g2b);
			\end{tikzpicture}
		}
		\caption{For \AGAG, the second player always has an answer to the first player's first move, ensuring that exactly two moves will remain afterwards (the last two leaves), and thus that the second player wins.}
		\label{fig-differencesAgag}
	\end{subfigure}
	\caption{An example of a graph with different outcomes for \GAG and \AGAG.}
	\label{fig-differencesActiveNonActive}
\end{figure}

Other variants were later introduced, such as the \emph{hull game}, where the \emph{convex hull} (instead of the geodetic closure) of $S$ is considered~\cite{benesh2025impartial,fraenkel1989geodetic}. The outcomes for trees, cycles and complete graphs with rays were determined in~\cite{fraenkel1989geodetic}, while~\cite{benesh2025impartial} focused on Sprague-Grundy values for (among others) block graphs, lattice graphs, (generalized) wheels and complete multipartite graphs, as well as algebraic properties of the game.

Note that those games can be generalized as hypergraph games~\cite{sieben2023impartial}, where two players alternate selecting vertices of a hypergraph (potentially modifying it or under conditions) until those satisfy some condition related to the hyperedges. In particular, in \GAG, the players add hyperedges corresponding to the shortest paths between the selected vertices, until every vertex is in at least one hyperedge.

\subsection{Game definitions}
\label{sec-games}

\GAG and \AGAG are impartial games, a subset of combinatorial games which have been extensively studied.
Since both players have the same possible moves, impartial games can have two possible \emph{outcomes}: either the first player has a winning strategy (\outcomeN) or the second player does (\outcomeP). A position $P'$ that can be reached by a move from a position $P$ is called an \emph{option} of $P$, the set of options of $P$ is denoted by $\opt(P)$. If every option of $P$ is \outcomeN, then $P$ is \outcomeP; conversely, if any option is \outcomeP, then $P$ is \outcomeN (the winning strategy for the first player is to play to a \outcomeP option). A refinement of outcomes are the \emph{Sprague-Grundy values}~\cite{grundy1939mathematics,sprague1935mathematische} (also called \emph{nim-values} or \emph{nimbers}), denoted by $\gr(P)$ for a given position $P$ and which can be mapped to nonnegative integers. The value of a position can be computed inductively, with $\gr(P)=\mex(\{\gr(P')~:~P' \in \opt(P)\})$, where $\mex(X)$ (for a set $X$ of nonnegative integers) is the smallest nonnegative integer not in $X$. A position $P$ is \outcomeP if and only if $\gr(P) = 0$. Two positions $P$ and $P'$ are said to be \emph{equivalent}, denoted by $P \equiv P'$, if $\gr(P)=\gr(P')$.

The \emph{disjoint sum} of two positions $P$ and $Q$, denoted by $P+Q$, is a position where each player chooses to play on one of the two positions at their turn. When one of them is finished, the game continues on the other one; the last player to play wins. For \AGAG, a disjoint sum can, for example, represent the disjoint union of two graphs. Computing the Sprague-Grundy value of a disjoint sum $P+Q$ can be done by applying the so-called \emph{nim-sum}, which is the bitwise XOR: $\gr(P+Q)=\gr(P) \oplus \gr(Q)$. An example is given below for $7 \oplus 10$.
\begin{center}
	\begin{tikzpicture}
		\draw (0,0) node {7};
		\draw (0.5,0) node {$\rightarrow$};
		\draw (1.25,0) node {1};
		\draw (1.5,0) node {1};
		\draw (1.75,0) node {1};
		\draw (0,-0.5) node {10};
		\draw (0.5,-0.5) node {$\rightarrow$};
		\draw (1,-0.5) node {1};
		\draw (1.25,-0.5) node {0};
		\draw (1.5,-0.5) node {1};
		\draw (1.75,-0.5) node {0};
		\draw (0.875,-0.75) to (1.875,-0.75);
		\draw (1,-1) node {1};
		\draw (1.25,-1) node {1};
		\draw (1.5,-1) node {0};
		\draw (1.75,-1) node {1};
		\draw (2.25,-1) node {$\rightarrow$};
		\draw (2.75,-1) node {13};
	\end{tikzpicture}
\end{center}
For more information and history about combinatorial games, we refer the reader to~\cite{albert2019lessons,berlekamp2004winning,conway2000numbers,siegel2013combinatorial}.

Due to the definition of the sum of positions, we can limit our study to connected graphs: studying a graph with several connected components is equivalent to studying each of the components independently and applying the nim-sum.
Furthermore, we will always consider that both players play optimally. Optimal play is understood as, applying a winning strategy using as few rounds as possible if possible, and trying to make the game last as long as possible if no winning strategy exists.
We will also make use of the following notation:

\begin{notation}
	Let $G$ be a graph. We also denote the game position by $G$. Furthermore, let $S$ be a subset of vertices of $G$, then the position denoted by $G,S$ is the graph $G$ where the vertices in $S$ have already been selected.
\end{notation}

\subsection{Structure of the paper}

We will begin in \Cref{sec-generalAndProducts} by proving general results on \AGAG as well as studying how to compute the outcome for the Cartesian product of several graphs in function of the outcomes of the individual graphs. We will also show that the tensor product cannot be analyzed as straightforwardly. Then, in \Cref{sec-lovesMeLovesMeNot}, we will study some graph classes where we can predict exactly which vertices will be picked. Afterwards, \Cref{sec-symmetry} will be focused on highly symmetric graphs such as paths, cycles and grids, and we compute the Sprague-Grundy values for the former two. Finally, in \Cref{sec-blockAndCacti}, we present polynomial-time algorithms based on dynamic programming computing the Sprague-Grundy values of block graphs and cacti.

\section{General results and graph products}
\label{sec-generalAndProducts}

In this section, we present a few general results that help the study of \AGAG, be it by stating that some vertices have to be selected, by helping decompose a graph, or by showing how some graph products affect the outcome of graphs.

First, we note that, in \AGAG, due to the fact that vertices can be in the geodetic closure $(S)$ without being in the geodetic set $S$ only if they are on some shortest path between two vertices, some vertices always have to be selected. For example, a leaf (\emph{i.e.}, a degree~1 vertex), can never be in a shortest path between two other vertices. More generally, consider \emph{simplicial} vertices, that is, vertices whose neighbourhood induces a complete subgraph.

\begin{lemma}
	\label{lem-simplicial}
	In \AGAG, all simplicial vertices of the graph will have to be selected.
\end{lemma}

\begin{proof}
	A simplicial vertex can never be on a shortest path between two other vertices, since all its neighbours are fully connected to each other. Hence, it can never be in a geodetic closure $(S)$ if it is not in $S$ itself.
\end{proof}

\emph{Articulation points}, that is, vertices whose removal disconnects the graph (also called \emph{cut vertices} or \emph{separating vertices}), are important for the study of \AGAG. The following lemma will often be used to allow for decompositions:

\begin{lemma}
	\label{cut_vertex_in_closure}
	Let $G$ be a graph with an articulation point $u$, such that components $G_1,\ldots, G_k$ are attached at $u$. In \AGAG, $\gr(G,\{u\})=\gr(G_1,\{u\})\oplus \ldots \oplus \gr(G_k,\{u\})$.
\end{lemma}

\begin{proof}
	When $u$ is selected, any selected vertex $x$ in $G_i$ will never help cover any vertex in $G_j$ for $j \neq i$. Indeed, all the shortest paths from $x$ to any vertex in $G_j$ will go through $u$, so, from the point of view of $G_j$ selecting $x$ after $u$ has no impact. Hence, each $G_i,\{u\}$ is an independent position, and we are playing on their disjoint union. See \Cref{fig-ArticulationPoint} for an illustration.
\end{proof}

\begin{figure}[h]
	\centering
	\scalebox{1}{
		\begin{tikzpicture}
			\node (graph) at (0,0) {
				\begin{tikzpicture}
					\draw (0,0) ellipse (1 and 0.5);
					\draw (0,0) node {$G_1$};
					\draw (2.25,0) ellipse (1.25 and 0.625);
					\draw (2.25,0) node {$G_2$};
					\node[noeud,fill=black] (u) at (1,0) {};
					\draw (u) node[above,yshift=1.5mm] {$u$};
				\end{tikzpicture}
			};
			
			\draw (2.875,0) node {$\equiv$};
			
			\node (graphs) at (6,0) {
				\begin{tikzpicture}
					\draw (0,0) ellipse (1 and 0.5);
					\draw (0,0) node {$G_1$};
					\draw (3,0) ellipse (1.25 and 0.625);
					\draw (3,0) node {$G_2$};
					\node[noeud,fill=black] (u) at (1,0) {};
					\draw (u) node[above,yshift=1.5mm] {$u$};
					\node[noeud,fill=black] (v) at (1.75,0) {};
					\draw (v) node[above,yshift=1.5mm] {$u$};
					\draw (1.375,0) node {+};
				\end{tikzpicture}
			};
		\end{tikzpicture}
	}
	\caption{Selecting an articulation point is equivalent to splitting the graph and playing independently on the components.}
	\label{fig-ArticulationPoint}
\end{figure}

We also focus on a classical question whenever studying games: if we know the outcomes of two graphs, can we obtain the outcome of some composition of those graphs? This question stems from the disjoint sum, which is managed through the nim-sum of Sprague-Grundy values. We study how the Cartesian product of graphs interacts with the outcomes for \AGAG.

Given two graphs $G$ and $H$, their \emph{Cartesian product} $G \Box H$ has vertex set $V(G \Box H)= \{(u, v)~|~u \in V(G)~\&~v \in V(H)\}$, and two vertices $(u_1, v_1)$ and $(u_2, v_2)$ are adjacent in $G \Box H$ if and only if
\begin{itemize}
	\item $u_1 = u_2$ and $v_1$ is adjacent to $v_2$ in $H$, or
	\item $v_1 = v_2$ and $u_1$ is adjacent to $u_2$ in $G$.
\end{itemize}

\begin{theorem}
	\label{thm-cartesianProductOfNandN}
	Let $G = G_1 \Box G_2 \Box \cdots \Box G_n$ be the Cartesian product of $n$ graphs. The graph $G$ is \outcomeN for \AGAG  if and only if  $G_i$ is \outcomeN for \AGAG  for all  $i \in \{1, 2, \cdots, n\}$. 
\end{theorem}

\begin{proof}
	We prove the statement by induction on \(n\).
	For \(n=1\), the statement is trivial.
	
	Assume the result holds for \(n-1\), and let
	\(
	G = G_1 \Box G_2 \Box \cdots \Box G_n.
	\)
	Let \(G' = G_1 \Box G_2 \Box \cdots \Box G_{n-1}\). By the induction hypothesis, \(G'\) is \outcomeN\ if and only if each \(G_i\) is \outcomeN\ for \(i=1,\dots,n-1\). Hence, it suffices to prove the statement for $n=2$. That is, it is enough to show that for any graphs \(H_1\) and \(H_2\),
	\[
	H_1 \Box H_2 \text{ is \outcomeN} \quad \text{if and only if} \quad H_1 \text{ and } H_2 \text{ are \outcomeN}.
	\]
	
	We now prove this statement.
	Let \(G = G_1 \Box G_2\). We use the following well-known property of the Cartesian product of graphs (\cite{klavvzar2000product}),
	for any \(x_1,x_2 \in V(G_1)\) and \(y_1,y_2 \in V(G_2)\),
	\[
	d_G\big((x_1,y_1),(x_2,y_2)\big)
	= d_{G_1}(x_1,x_2) + d_{G_2}(y_1,y_2).
	\]
	
	Next, we have the following claim.

	\begin{claim}
		\label{clm-cartesianProductClaim2}
		Let $G = G_1 \Box G_2$, $S_1 \subseteq V(G_1)$, $S_2 \subseteq V(G_2)$, and $S = S_1 \times S_2 \subseteq V(G)$. We have $(u,v) \in (S)$ if and only if $(u,v) \in (S_1) \times (S_2)$. 
	\end{claim}
	
	\begin{proof}
		In one direction, assume that $(u,v) \in (S)$. Thus, 
		there exists $(x_1,y_1), (x_2,y_2) \in S$ such that $(u,v)$ lies on some shortest path $P$ between  $(x_1,y_1)$ and $(x_2,y_2)$. Let $P'$ be the subpath of $P$ that connects $(x_1, y_1)$ to $(u,v)$. Similarly, let   $P''$ be the subpath of $P$ that connects $(u,v)$ to $(x_2, y_2)$. Note that both $P'$ and $P''$ are shortest paths.

		Denote the length of a path $X$ (that is, its number of edges) by $\lt(X)$.
		By ~\cite{klavvzar2000product}, we know that there are shortest paths $P'_1$ in $G_1$ and $P'_2$ in $G_2$ connecting $x_1$ to $u$ and $y_1$ to $v$, respectively, satisfying 
		$\lt(P'_1)+\lt(P'_2) = \lt(P')$. Similarly, we know 
		that there are shortest paths $P''_1$ and $P''_2$ in $G_2$ connecting $u$ to $x_2$ and $v$ to $y_2$, respectively, satisfying $\lt(P''_1)+\lt(P''_2) = \lt(P'')$. 
		Let $P_1$ (resp., $P_2$) 
		be the path obtained through concatenation of $P'_1$ and $P''_1$
		(resp., $P'_2$ and $P''_2$). 
		Thus, we must have 
		$\lt(P_1)+\lt(P_2) = \lt(P)$. This implies that $P_1$ and $P_2$ are shortest paths due to ~\cite{klavvzar2000product}. Thus, $(u,v) \in (S)$ implies $(u,v) \in (S_1) \times (S_2)$.

		In the other direction, assume that $(u,v) \in (S_1) \times (S_2)$. Thus, there are $x_1, x_2 \in S_1$ such that $u$ lies on the path $P_1$ connecting $x_1$ and $x_2$. Moreover, without loss of generality, we may assume that the path is the concatenation of the subpaths $P'_1$ connecting $x_1$ and $u$ and $P''_1$ connecting $u$ and $x_2$. 
		Similarly, 
		there are $y_1, y_2 \in S_2$ such that $v$ lies on the path $P_2$ connecting $x_1$ and $x_2$ where $P_2$ is the concatenation of the subpaths $P'_2$ connecting $y_1$ and $v$ and $P''_2$ connecting $v$ and $y_2$. By ~\cite{klavvzar2000product}, we know that there are two paths $P'$ connecting $(x_1, y_1)$ and $(u,v)$ and 
		$P''$ connecting $(u,v)$ and $(x_2, y_2)$ in $G$, and
		satisfying 
		$\lt(P'_1)+\lt(P'_2) = \lt(P')$ and $\lt(P''_1)+\lt(P''_2) = \lt(P'')$. In particular, this implies that $P'$ and $P''$ are shortest paths. 
		
		Let $P$ be the path obtained by concatenation of $P'$ and $P''$. Observe that $P$ connects $(x_1,y_1)$ and $(x_2, y_2)$ and passes through the vertex $(u,v)$.  Note that  $$\lt(P) = \lt(P') + \lt(P'') = \lt(P'_1)+\lt(P'_2)+\lt(P''_1)+\lt(P''_2) = \lt(P_1) + \lt(P_2).$$
		Thus, by ~\cite{klavvzar2000product}, $P$ must be a shortest path.
	\end{proof}
	
	\medskip
	\noindent
	In one direction, assume that $G_1$ and $G_2$ are \outcomeN. We show that
	$G = G_1 \Box G_2$ is also \outcomeN.

	We will play the game simultaneously on $G$, $G_1$ and $G_2$. 
	The vertex sets that will be constructed during the games are $S$, $S_1$ and $S_2$, respectively. We will maintain $S = S_1 \times S_2$. 
	The strategy of the first player on $G$ will be derived from the strategies of the first player on $G_1$ and $G_2$. The second player will play on $G$, and the corresponding moves will be replicated on $G_1$ and $G_2$. 
	
	Let us explain first how a move is replicated or translated. 
	At any point in the game, if a player selects the vertex $(u,v)$ in $G$, we know that $(u,v) \not\in (S)$ at the time it is selected. However, by \Cref{clm-cartesianProductClaim2}, it is possible to have either $u \in (S_1)$ or $v \in (S_2)$, but not both. If $u \in (S_1)$, then the selection of $(u,v)$ corresponds to a ``dummy move'' in $G_1$. 
	If $u \not\in (S_1)$, then the selection of $(u,v)$ corresponds to selecting $u$ in $G_1$. Similarly, one can translate the moves on $G$ to $G_2$.

	The selection of $u$ and $v$ in $G_1$ and $G_2$, respectively, corresponds to the selection of $(u,v)$ in $G$. Note that either $u \not\in (S_1)$ or $v \not\in (S_2)$ at the time of their respective selection, implying $(u,v) \not\in (S) = (S_1) \times (S_2)$ due to \Cref{clm-cartesianProductClaim2}, and thus selecting $(u,v)$ is a valid move.

	The first player will begin the game by selecting vertices $u_1, v_1$ in $G_1, G_2$, respectively, according to the corresponding winning strategies. In $G$, the vertex $(u_1,v_1)$ is selected by the first player. After that the players take turns to make their moves. 
	Let us assume that the second player has made a move by selecting the vertex $(u,v)$ in $G$. This moves translates to moves in $G_1, G_2$. However due to \Cref{clm-cartesianProductClaim2}, we know that at most one of these translated moves can be a dummy move. Without loss of generality assume that the second player's move translates to a dummy move in $G_1$, and a valid move in $G_2$. In that case, the first player makes a dummy move in $G_1$ and provides a response by selecting $v'$ in $G_2$ according to winning strategy. These moves are translated to $G$ as selection of the vertex $(u,v')$, which is a valid move due to \Cref{clm-cartesianProductClaim2}. If the second player's move in $G$ translates to valid moves in both $G_1$ and $G_2$, then the first player responds in both the graphs according to corresponding winning strategies by selecting $u',v'$ in $G_1, G_2$, respectively. This translates to the valid move of selecting $(u',v')$ in $G$ by the first player due to \Cref{clm-cartesianProductClaim2}.  
	
	Hence, the first player always has a valid response to the second player's move. Thus, $G$ is \outcomeN.

	\medskip

	In the other direction, assume that at least one of $G_1$ or $G_2$ is \outcomeP. 
	Without loss of generality, assume that $G_2$ is \outcomeP. We show that 
	$G = G_1 \Box G_2$ is also \outcomeP.
	
	The second player simulates the game on both $G_1$ and $G_2$, following a winning strategy on $G_2$. 
	Whenever the first player selects a vertex $(u,v)$ in $G$, this induces corresponding moves in $G_1$ and $G_2$. 
	The second player responds in $G_2$ according to the winning strategy, selecting a vertex $v'$. 
	In $G_1$, the second player either mirrors the first player's move or plays arbitrarily if the move is a dummy move. 
	
	These choices correspond to selecting a vertex $(u',v')$ in $G$, which is valid by Claim~\ref{clm-cartesianProductClaim2}. 
	Thus, the second player can always respond to any move of the first player while maintaining consistency with the winning strategy in $G_2$. 
	Since the game on $G_2$ ends with a win for the second player, the same holds for the game on $G$. 
	Therefore, $G$ is \outcomeP.
\end{proof}

\Cref{thm-cartesianProductOfNandN} shows that for the Cartesian product $G\Box H$, the outcome of the Closed Geodetic Game is determined exactly by the outcomes on the individual graphs of $G\Box H$: $G\Box H$ is \outcomeN if and only if both $G$ and $H$ are \outcomeN.  In what follows, we specialize to the case where $G=P_n$ and $H=P_m$ are paths.
This result naturally raises the question: does a similar outcome preservation hold for other products? This does not seem to be the case in general, as experienced by the case of the tensor product $G\times H$, which we are going to develop.

Even when both $G$ and $H$ are \outcomeN individually, the tensor product $G\times H$ can be \outcomeP. Thus, unlike the Cartesian case, outcomes are not preserved under the tensor product. This fundamental distinction illustrates that strategies in the tensor product setting cannot be straightforwardly composed from the component strategies, thus denying a straightforward solution for graph products in general.

We begin by recalling structural properties of the tensor product. Given two graphs $G$ and 
$H$, their \emph{tensor product} $G \times H$  has vertices of the form $(u,v)$ where $u \in V(G)$ and $v \in V(H)$, such that two vertices $(u_1, v_1)$ and $(u_2, v_2)$ are adjacent if and only if $u_1$ is adjacent to $u_2$ in $G$ and $v_1$ is adjacent to $v_2$ in $H$.

\begin{lemma}\cite{weichsel1962kronecker}\label{tensorlemma1}
	The tensor product $G \times H$ of the connected graphs $G$ and $H$ is connected if and only if at least one of the individual graphs $G$ or $H$ is non-bipartite.
\end{lemma} 

\begin{lemma}\cite{weichsel1962kronecker}\label{tensorlemma2}
	If both $G$ and $H$ are connected bipartite graphs, then the tensor product $G \times H$ consists of exactly two components.
\end{lemma}

\begin{theorem}\cite{jha1997isomorphic}\label{isomorphic}
	If $G$ is a bipartite graph and $H$ is a path with even number of vertices, then the two components of $G \times H$ are isomorphic.
\end{theorem}

The following serves as a concrete witness, illustrating that even when two graphs \( G \) and \( H \) have the same outcome individually, their tensor product may yield different outcomes. We begin by showing a graph family where the tensor products preserve the outcome: paths of a given parity.

First, consider that if $G=P_n$ and $H=P_m$ are odd paths with the vertex sequences
$u_1,\dots,u_n$ and $v_1,\dots,v_m$, respectively, then we denote
\[
c_n := \Big\lceil\frac{n}{2}\Big\rceil,\qquad
c_m := \Big\lceil\frac{m}{2}\Big\rceil.
\]
Define the \emph{central reflection} mapping
$\sigma:V(G\times H)\to V(G\times H)$ such that
\(    \sigma\big(u_{c_n+k},v_{c_m+\ell}\big)\;=\;\big(u_{c_n-k},v_{c_m-\ell}\big),
\)
for $k \in \{-\lfloor\frac{n}{2} \rfloor,\cdots, \lfloor \frac{n}{2} \rfloor\}$ and $\ell \in \{-\lfloor\frac{m}{2} \rfloor,\cdots, \lfloor \frac{m}{2} \rfloor\}$. Note that
$\sigma^2=\mathrm{id}$ and that $(u_{c_n},v_{c_m})$ is the unique fixed
point of $\sigma$. 

\medskip

We recall some standard definitions. Two graphs \(G\) and \(H\) are \emph{isomorphic}, denoted by \(G \cong H\), if there exists a bijection
\(\phi : V(G) \to V(H)\)
such that for all \(u,v \in V(G)\),
\[
uv \in E(G) \quad \iff \quad \phi(u)\phi(v) \in E(H).
\]
The function \(\phi\) is called a \emph{graph isomorphism} between \(G\) and \(H\). Given two isomorphic graphs \(G\) and \(H\), there may be many different isomorphisms \(\phi: G \to H\). A \emph{canonical isomorphism} (if it exists) is a distinguished isomorphism determined by the natural structure of \(G\) and \(H\), without needing arbitrary choices. 

An \emph{automorphism} of a graph \(G\) is a graph isomorphism from \(G\) to itself, that is, a bijection \(\sigma : V(G) \to V(G)\) such that for all \(u,v \in V(G)\),
\[
uv \in E(G) \quad \iff \quad \sigma(u)\sigma(v) \in E(G).
\]

Here, $\sigma$ denotes the central reflection mapping defined earlier for $G\times H$.

\begin{lemma}\label{lem:sigma-automorphism}
	The map $\sigma$ is a graph automorphism of $G\times H$. Moreover,  
	for every $x,y\in V(G\times H)$,
	\[
	d_{G\times H}(\sigma(x),\sigma(y))=d_{G\times H}(x,y).
	\]
	Consequently, $\sigma$ maps shortest paths to shortest paths.
\end{lemma}

\begin{proof}
	The map $\sigma$ is a bijection and $\sigma^2=\mathrm{id}$ (\emph{i.e.}, reflection twice is the identity), so it is invertible.
	Now take any adjacent vertices \((u_i,v_j), (u_{i'},v_{j'}) \in V(G\times H)\). By the tensor product definition, we have
	\[
	u_i u_{i'} \in E(G) \quad \text{and} \quad v_j v_{j'} \in E(H).
	\]
	
	Apply \(\sigma\) to these vertices such that
	\(
	\sigma(u_i,v_j) = (u_{i^*}, v_{j^*}), \quad
	\sigma(u_{i'},v_{j'}) = (u_{i'^*}, v_{j'^*}),
	\)
	where \(i^* = 2c_n - i\) and \(j^* = 2c_m - j\), and similarly for \(i'^*, j'^*\). 
	Since \(u_iu_{i'}\in E(G)\), their indices differ by one. Reflection about \(c_n\) preserves this difference, so \(u_{i^*}u_{i'^*} \in E(G)\) (i.e., $d_G(u_{i^*},u_{i'^*} )=1$). Similarly, \(v_jv_{j'} \in E(H)\) implies \(v_{j^*}v_{j'^*} \in E(H)\).  
	
	Thus, the images under \(\sigma\) are also adjacent in \(G \times H\):
	\[
	\sigma(u_i,v_j) \sigma(u_{i'},v_{j'}) \in E(G \times H).
	\]
	Since \(\sigma\) is a bijection that preserves adjacency, it is a graph automorphism. As an automorphism, it preserves lengths of paths. Therefore, if \(P\) is a shortest path between two vertices in \(G \times H\), then \(\sigma(P)\) is a shortest path between the images of those vertices. In particular, this ensures that \(\sigma\) maps shortest paths to shortest paths.
\end{proof}

\begin{proposition}\label{P_ntimes P_m}
	Let $G = P_n$ and $H = P_m$ be paths of length $n$ and $m$, respectively. If both $G$ and $H$ have the same outcome in \AGAG, then their tensor product $G \times H$ will also yield the same outcome.
\end{proposition}

\begin{proof}
	We assume that $G$ and $H$ have the same outcome in \AGAG. In particular, throughout this proof we specialize to $G=P_n$ and $H=P_m$. %We aim to show that the tensor products $G\times H$ have same outcome. 
	The tensor product $G \times H$ has vertex set $\{(u,v)\in V(G\times H):u\in V(G),v\in V(H)\}$, with adjacency defined as $(u_1,v_1)(u_2,v_2)\in E(G\times H)$ if and only if $u_1u_2 \in E(G)$ and $v_1,v_2\in E(H)$.
	Due to \Cref{tensorlemma1,tensorlemma2}, since path graphs are bipartite, $P_n \times P_m$ is disconnected for all $m,n\geq 2$. Consider $P_n \times P_m$, and let $u_1, \dots, u_n$ and $v_1, \dots, v_m$ denote the vertices of $P_n$ and $P_m$, respectively. A vertex in the tensor product is denoted $(u_i, v_j)$.
	By \Cref{thm-paths}, the outcome of a path depends only on the parity of its order: odd paths are \outcomeN and even paths are \outcomeP. Hence, if $G$ and $H$ have the same outcome, then both are odd paths or both are even paths. We treat these two cases separately below. Before starting the game, we construct a matrix $M$ with the rows representing the copies of vertices of $G$, and the columns representing the copies of vertices of $H$; hence each cell of $M$ is a vertex $(u_{k},v_{\ell})$ of $G \times H$ where $u_{k}\in V(G)$ and $v_{\ell}\in V(H)$.
	
	Suppose first that $n$ and $m$ are odd. Then both $P_n$ and $P_m$ are \outcomeN for \AGAG on both $G$ and $H$.  
	Define the map \(
	\sigma: V(G)\times V(H)\to V(G)\times V(H)
	\)
	by
	\[
	\sigma\big(u_{c_n+k},v_{c_m+\ell}\big) \;=\; \big(u_{c_n-k},v_{c_m-\ell}\big)
	\]
	for $k \in \{0,1,2,\cdots, \lfloor \frac{n}{2} \rfloor\}$ and $\ell \in \{0,1,2,\cdots, \lfloor \frac{m}{2} \rfloor\}$. Denote \(\sigma^{-1}=\sigma\), and note that \(\sigma(u_{c_n},v_{c_m})=(u_{c_n},v_{c_m})\) is a fixed point and the central vertex with both \(k=\ell=0\).
	
	We will use \(\sigma\) to define a mirror strategy. To do so, we establish the following two claims:

	\begin{claim}
		\label{lem:sigma-closure}
		For any set \(S\subseteq V(G\times H)\), if  \(\sigma(S)=S\), then \(\sigma\big((S)\big)=(S)\).
	\end{claim}
	
	\begin{proof}
		Let \(x\in(S)\). By definition there exist \(a,b\in S\) such that \(x\) lies on some shortest path \(P\) between \(a\) and \(b\). Apply \(\sigma\) to every vertex of \(P\): because \(\sigma\) reflects indices coordinatewise about the centers, edges map to edges and shortest paths map to shortest paths (length is preserved). Since \(a,b\in S\) and \(\sigma(S)=S\), we have \(\sigma(a),\sigma(b)\in S\). By Lemma~\ref{lem:sigma-automorphism},
		$\sigma(P)$ is a shortest path between $\sigma(a)$ and $\sigma(b)$ and has
		the same length, and then $\sigma(x)\in\sigma(P)$. 
		This shows \(\sigma\big((S)\big)\subseteq (S)\). The reverse inclusion follows by applying the same argument to \(\sigma\big((S)\big)\) since \(\sigma^{-1}=\sigma\). Thus \(\sigma\big((S)\big)=(S)\).
	\end{proof}
	
	\begin{claim}
		\label{lem:legal-reply}  
		Let the current position of the game is \((G\times H,S)\) with $S$ satisfying $\sigma(S)=S$. If the opponent plays
		a legal move \(x\notin(S)\), then the symmetric mate $\sigma(x)$ is also a legal move: $\sigma(x)\notin S$
		and $\sigma(x)\notin (S)$.
	\end{claim}
	
	\begin{proof}
		Suppose, for contradiction, that \(\sigma(x)\in(S)\). By Claim~\ref{lem:sigma-closure} we have 
		$\sigma\big((S)\big)=(S)$, hence applying $\sigma$ we
		obtain $x=\sigma(\sigma(x))\in (S)$, contradicting the hypothesis \(x\notin (S)\). Therefore, \(\sigma(x)\notin (S)\).
		Also $x\notin S$ and $\sigma(S)=S$ imply $\sigma(x)\notin S$. This proves
		the claim.  %\(\square\) %\hfill $\circ$
	\end{proof}

	Now, we will show how to apply the winning strategies on $G \times H$. The first player starts by playing the unique fixed central vertex \((u_{c_n},v_{c_m})\). Let $S=\{(u_{c_n},v_{c_m})\}$. Clearly $\sigma(S)=S$.
	By Claim~\ref{lem:legal-reply}, whenever the second player later plays some vertex $x\notin (S)$, the first player can immediately reply with the symmetric mate $\sigma(x)$,
	and this reply is legal. Thus, the first player maintains the invariant that
	the current selected set is $\sigma$-invariant (after each pair of moves, the
	set of selected vertices is closed under $\sigma$).
	By doing so, due to Claim~\ref{lem:sigma-closure}, the first player ensures that they will always have an answer to the second player's move in every step. Therefore, they will always win in $G \times H$. Hence, $G \times H$ also has outcome \outcomeN.

	Suppose now that $n$ and $m$ are even. Then both $G$ and $H$ are outcome \outcomeP for \AGAG. Since both $P_n$ and $P_m$ are bipartite graphs and $n, m$ are even, 
	it follows from \Cref{isomorphic} that 
	$G \times H$ has exactly two isomorphic connected components, denoted by $C_1$ and $C_2$.
	
	Since $C_1$ and $C_2$ are isomorphic by Theorem~\ref{isomorphic}, there exists a graph isomorphism 
	\[
	\phi: V(C_1) \longrightarrow V(C_2)
	\]
	such that for all $(x_i,v_j),(x_{i'},v_{j'})\in V(C_1)$, 
	\((x_i,v_j)(x_{i'},v_{j'})\in E(C_1)
	\quad\Longleftrightarrow\quad
	\phi(x_i,v_j)\phi(x_{i'},v_{j'})\in E(C_2).\) Equivalently, we can describe $\phi$ as follows:
	\[
	\phi(x_i,v_j) = (x',v') 
	\quad \text{where } (x_i,v_j)\in V(C_1),~(x',v')\in V(C_2),
	\]
	and $(x_i,v_j)$ and $(x',v')$ correspond under the canonical isomorphism between the two components.
	In particular, $\phi$ is a bijection satisfying $\phi^{-1}:V(C_2)\to V(C_1)$.
	
	We now use this isomorphism to apply a symmetry strategy for \AGAG:
	In the initial move, the first player chooses any vertex $(x_i,v_j)$ in one component, say $C_1$, then
	the second player responds by selecting its image $\phi(x_i,v_j)=(x',v')\in C_2$. Now, Let $S_1=\{(x_i,v_j)\}$ in $C_1$ and $S_2=\{(x',v')\}$ in $C_2$. Thereafter, the first player can select any vertex not in $(S)$, either $C_1$ or $C_2$, and the second player responds by the image of the bijection mapping $\phi$.

	\begin{claim}\label{lem:phi-preserves-geodetic}
		Let $\phi:V(C_1)\to V(C_2)$ be a graph isomorphism between the two
		components of $G\times H$.  For any subset $S_1\subseteq V(C_1)$ and set
		$S_2:=\phi(S_1)\subseteq V(C_2)$, 
		\[
		\phi\big( (S_1)_{C_1} \big) \;=\; (S_2)_{C_2},
		\]
		where \((\cdot)_{C_r}\) denotes geodetic closure taken inside component
		\(C_r\).
	\end{claim}
	
	\begin{proof}
		Let $x\in (S_1)_{C_1}$. By definition, there exist $a,b\in S_1$ such that
		$x$ lies on a shortest path $P$ between $a$ and $b$ in $C_1$. Since
		$\phi$ is an isomorphism and it
		preserves adjacency and path-lengths, $\phi(P)$ is a shortest path between
		$\phi(a),\phi(b)\in S_2$ in $C_2$ and $d_{C_1}(P)=d_{C_2}(\phi(P))$. Thus, $\phi(x)$ lies on $\phi(P)$ . Hence
		$\phi(x)\in (S_2)_{C_2}$, proving $\phi\big((S_1)_{C_1}\big)\subseteq (S_2)_{C_2}$.
		Since $\phi$ is a bijection, there is the same argument for $\phi^{-1}$. Thus equality holds.
	\end{proof}

	Since $\phi$ is a bijection, by doing the above process, the second player can always choose a vertex after the first player's move. After a finite number of moves, \Cref{lem:phi-preserves-geodetic} ensures that the game can only end after a second player's move. Hence, the second player will always win in $G \times H$.
\end{proof}

It is worth observing that the above argument extends naturally to all cases in which at least one of the factors has an even order. Indeed, whenever either 
$P_n$ or $P_m$ is an even path, Theorem~\ref{isomorphic} ensures that the tensor product $P_n\times P_m$ decomposes into two isomorphic connected components. Hence, there is a canonical isomorphism between these components. Therefore, this case will follow the even-even case.

This phenomenon, however, relies crucially on the parity structure of paths and does not extend to general bipartite graphs.
In fact, even when both factors share the same individual outcome, the tensor product may exhibit a different behaviour, as shown below.

\begin{figure}
	\centering
	\begin{subfigure}[b]{\textwidth}
		\centering
		\scalebox{0.85}{
			\begin{tikzpicture} [scale=0.85, transform shape]
				\node (A) at (-1,0) {
					\begin{tikzpicture}
						\node[noeud] (au) at (0,4) {};
						\node[noeud] (av) at (1,4) {};
						\node[noeud] (aw) at (2,4) {};
						\node[noeud] (bu) at (0,3) {};
						\node[noeud] (bv) at (1,3) {};
						\node[noeud] (bw) at (2,3) {};
						\node[noeud] (cu) at (0,2) {};
						\node[noeud,fill=black] (cv) at (1,2) {};
						\node[noeud] (cw) at (2,2) {};
						\node[noeud] (du) at (0,1) {};
						\node[noeud] (dv) at (1,1) {};
						\node[noeud] (dw) at (2,1) {};
						\node[noeud] (eu) at (0,0) {};
						\node[noeud] (ev) at (1,0) {};
						\node[noeud] (ew) at (2,0) {};
						\draw[blue, thick] (au)to(bv)to(cu)to(dv)to(eu);
						\draw[red, thick] (av)to(bu)to(cv)to(du)to(ev);
						\draw [red, thick] (av)to(bw)to(cv)to(dw)to(ev);
						\draw[blue, thick] (aw)to(bv)to(cw)to(dv)to(ew);
						
						\foreach \I in {0,...,4} {
							\node[noeud] (u\I) at (-1,\I) {};
						}
						\draw (u0)to(u1)to(u2)to(u3)to(u4);
						\foreach \I in {0,1,2} {
							\node[noeud] (v\I) at (\I,5) {};
						}
						\draw (v0)to(v1)to(v2);
						
						\draw[line width=0.25mm,dashed] (u2) circle (0.25);
						\draw[line width=0.25mm,dashed] (cv) circle (0.25);
						\draw[line width=0.25mm,dashed] (v1) circle (0.25);
				\end{tikzpicture}};
				
				\node (B) at (6,0) {
					\begin{tikzpicture}[scale=0.85, transform shape]
						\node[noeud] (au) at (-1,4) {};
						\node[noeud] (bv) at (0,3) {};
						\node[noeud] (cu) at (-1,2) {};   
						\node[noeud] (dv) at (0,1) {};
						\node[noeud] (eu) at (-1,0) {};
						
						\node[noeud] (aw) at (1,4) {};
						\node[noeud] (cw) at (1,2) {};
						\node[noeud] (ew) at (1,0) {};
						
						\node[noeud] (bu) at (2,3) {};
						\node[noeud] (av) at (3,4) {};
						\node[noeud] (bw) at (4,3) {};
						\node[noeud,fill=black] (cv) at (3,2) {};
						\node[noeud] (du) at (4,1) {};
						\node[noeud] (ev) at (3,0) {};
						\node[noeud] (dw) at (2,1) {};
						
						\draw[blue, thick] (au)to(bv)to(cu)to(dv)to(eu);
						\draw[red, thick] (av)to(bu)to(cv)to(du)to(ev);
						\draw [red, thick] (av)to(bw)to(cv)to(dw)to(ev);
						\draw[blue, thick] (aw)to(bv)to(cw)to(dv)to(ew);

				\end{tikzpicture}};
				\draw[->] (1.2,0) to (2.8,0);
		\end{tikzpicture}}
		\caption{The first move: applying the winning strategy on both graphs (the vertices circled with dashes) and this tensor product $P_5 \times P_3$ has $2$ components.}
		\label{fig-tensorProductNandN-1}
	\end{subfigure}
	
	\begin{subfigure}[b]{0.45\textwidth}
		\centering
		\scalebox{0.85}{
			\begin{tikzpicture}[scale=0.85, transform shape]
				\node[noeud] (au) at (0,4) {};
				\node[noeud] (av) at (1,4) {};
				\node[noeud] (aw) at (2,4) {};
				\node[noeud] (bu) at (0,3) {};
				\node[noeud,red,fill=red] (bv) at (1,3) {};
				\node[noeud] (bw) at (2,3) {};
				\node[noeud] (cu) at (0,2) {};
				\node[noeud,fill=black] (cv) at (1,2) {};
				\node[noeud] (cw) at (2,2) {};
				\node[noeud] (du) at (0,1) {};
				\node[noeud,blue,fill=blue] (dv) at (1,1) {};
				\node[noeud] (dw) at (2,1) {};
				\node[noeud] (eu) at (0,0) {};
				\node[noeud] (ev) at (1,0) {};
				\node[noeud] (ew) at (2,0) {};
				\draw[blue, thick] (au)to(bv)to(cu)to(dv)to(eu);
				\draw[red, thick] (av)to(bu)to(cv)to(du)to(ev);
				\draw [red, thick] (av)to(bw)to(cv)to(dw)to(ev);
				\draw[blue, thick] (aw)to(bv)to(cw)to(dv)to(ew);
				
				\foreach \I in {0,...,4} {
					\node[noeud] (u\I) at (-1,\I) {};
				}
				\draw (u0)to(u1)to(u2)to(u3)to(u4);
				\foreach \I in {0,1,2} {
					\node[noeud] (v\I) at (\I,5) {};
				}
				\draw (v0)to(v1)to(v2);
				\draw[thick,line width=0.25mm,dotted,red] (u3) circle (0.25);
				\draw[line width=0.25mm,dotted,red] (bv) circle (0.25);
				\draw[line width=0.25mm,dashed,blue] (u1) circle (0.25);
				\draw[line width=0.25mm,dashed,blue] (v1) circle (0.35);
				\draw[line width=0.25mm,dashed,blue] (dv) circle (0.25);
				\draw[line width=0.25mm,dotted,red] (v1) circle (0.25);
				\draw[<->] (2.5,3)to(3,3)to(3,1)to(2.5,1);
			\end{tikzpicture}
		}
		\caption{If the second player plays in the same column (dotted vertices), answer on the same (dashed vertices), and associate the rows.}
		\label{fig-tensorProductNandN-2}
	\end{subfigure}\hfil
	\begin{subfigure}[b]{0.45\textwidth}
		\centering
		\scalebox{0.85}{
			\begin{tikzpicture}[scale=0.85, transform shape]
				\node[noeud] (au) at (0,4) {};
				\node[noeud] (av) at (1,4) {};
				\node[noeud] (aw) at (2,4) {};
				\node[noeud] (bu) at (0,3) {};
				\node[noeud] (bv) at (1,3) {};
				\node[noeud] (bw) at (2,3) {};
				\node[noeud, red, fill=red] (cu) at (0,2) {};
				\node[noeud,fill=black] (cv) at (1,2) {};
				\node[noeud, blue,fill=blue] (cw) at (2,2) {};
				\node[noeud] (du) at (0,1) {};
				\node[noeud] (dv) at (1,1) {};
				\node[noeud] (dw) at (2,1) {};
				\node[noeud] (eu) at (0,0) {};
				\node[noeud] (ev) at (1,0) {};
				\node[noeud] (ew) at (2,0) {};
				\draw[blue, thick] (au)to(bv)to(cu)to(dv)to(eu);
				\draw[red, thick] (av)to(bu)to(cv)to(du)to(ev);
				\draw [red, thick] (av)to(bw)to(cv)to(dw)to(ev);
				\draw[blue, thick] (aw)to(bv)to(cw)to(dv)to(ew);
				
				\foreach \I in {0,...,4} {
					\node[noeud] (u\I) at (-1,\I) {};
				}
				\draw (u0)to(u1)to(u2)to(u3)to(u4);
				\foreach \I in {0,1,2} {
					\node[noeud] (v\I) at (\I,5) {};
				}
				\draw (v0)to(v1)to(v2);
				
				\draw[line width=0.25mm,dotted,red] (u2) circle (0.25);
				\draw[line width=0.25mm,dotted,red] (cu) circle (0.25);
				\draw[line width=0.25mm,dashed,blue] (cw) circle (0.25);
				\draw[line width=0.25mm,dashed,blue] (u2) circle (0.35);
				\draw[line width=0.25mm,dotted,red] (v0) circle (0.25);
				\draw[line width=0.25mm,dashed,blue] (v2) circle (0.25);
				\draw[<->] (0,5.5)to(0,6)to(2,6)to(2,5.5);
			\end{tikzpicture}
		}
		\caption{If the second player plays in the same row (dotted vertices), answer on the same (dashed vertices), and associate the columns.}
		\label{fig-tensorProductNandN-3}
	\end{subfigure}
	
	\begin{subfigure}[b]{0.45\textwidth}
		\centering
		\scalebox{0.85}{\begin{tikzpicture}[scale=0.85, transform shape]
				\node[noeud] (au) at (0,4) {};
				\node[noeud] (av) at (1,4) {};
				\node[noeud] (aw) at (2,4) {};
				\node[noeud] (bu) at (0,3) {};
				\node[noeud] (bv) at (1,3) {};
				\node[noeud,red,fill=red] (bw) at (2,3) {};
				\node[noeud] (cu) at (0,2) {};
				\node[noeud,fill=black] (cv) at (1,2) {};
				\node[noeud] (cw) at (2,2) {};
				\node[noeud,blue,fill=blue] (du) at (0,1) {};
				\node[noeud] (dv) at (1,1) {};
				\node[noeud] (dw) at (2,1) {};
				\node[noeud] (eu) at (0,0) {};
				\node[noeud] (ev) at (1,0) {};
				\node[noeud] (ew) at (2,0) {};
				\draw[blue, thick] (au)to(bv)to(cu)to(dv)to(eu);
				\draw[red, thick] (av)to(bu)to(cv)to(du)to(ev);
				\draw [red, thick] (av)to(bw)to(cv)to(dw)to(ev);
				\draw[blue, thick] (aw)to(bv)to(cw)to(dv)to(ew);

				\foreach \I in {0,...,4} {
					\node[noeud] (u\I) at (-1,\I) {};
				}
				\draw (u0)to(u1)to(u2)to(u3)to(u4);
				\foreach \I in {0,1,2} {
					\node[noeud] (v\I) at (\I,5) {};
				}
				\draw (v0)to(v1)to(v2);
				
				\draw[line width=0.25mm,dotted,red] (v2) circle (0.25);
				\draw[line width=0.25mm,dotted,red] (u3) circle (0.25);
				\draw[line width=0.25mm,dotted,red] (bw) circle (0.25);
				\draw[line width=0.25mm,dashed,blue] (v0) circle (0.25);
				\draw[line width=0.25mm,dashed,blue] (du) circle (0.25);
				\draw[line width=0.25mm,dashed,blue] (u1) circle (0.25);
				\draw[<->] (0,5.5)to(0,6)to(2,6)to(2,5.5);
				\draw[<->] (2.5,3)to(3,3)to(3,1)to(2.5,1);
			\end{tikzpicture}
			
		}
		\caption{If the second player plays distinctly (dotted vertices), answer on the same component (dashed vertices), and associate the rows and columns.}
		\label{fig-tensorProductNandN-4}
	\end{subfigure}\hfil
	\begin{subfigure}[b]{0.45\textwidth}
		\centering
		\scalebox{0.85}{\begin{tikzpicture}[scale=0.85, transform shape]
				\node[noeud,red, fill=red] (au) at (0,4) {};
				\node[noeud] (av) at (1,4) {};
				\node[noeud] (aw) at (2,4) {};
				\node[noeud] (bu) at (0,3) {};
				\node[noeud] (bv) at (1,3) {};
				\node[noeud] (bw) at (2,3) {};
				\node[noeud] (cu) at (0,2) {};
				\node[noeud,fill=black] (cv) at (1,2) {};
				\node[noeud] (cw) at (2,2) {};
				\node[noeud] (du) at (0,1) {};
				\node[noeud] (dv) at (1,1) {};
				\node[noeud] (dw) at (2,1) {};
				\node[noeud] (eu) at (0,0) {};
				\node[noeud] (ev) at (1,0) {};
				\node[noeud,blue, fill=blue] (ew) at (2,0) {};
				\draw[blue, thick] (au)to(bv)to(cu)to(dv)to(eu);
				\draw[red, thick] (av)to(bu)to(cv)to(du)to(ev);
				\draw [red, thick] (av)to(bw)to(cv)to(dw)to(ev);
				\draw[blue, thick] (aw)to(bv)to(cw)to(dv)to(ew);

				\foreach \I in {0,...,4} {
					\node[noeud] (u\I) at (-1,\I) {};
				}
				\draw (u0)to(u1)to(u2)to(u3)to(u4);
				\foreach \I in {0,1,2} {
					\node[noeud] (v\I) at (\I,5) {};
				}
				\draw (v0)to(v1)to(v2);
				
				\draw[line width=0.25mm,dashed,blue] (v2) circle (0.25);
				\draw[line width=0.25mm,dotted,red] (v0) circle (0.25);
				\draw[line width=0.25mm,dotted,red] (au) circle (0.25);
				\draw[line width=0.25mm,dotted,red] (u4) circle (0.25);
				\draw[line width=0.25mm,dashed,blue] (u0) circle (0.25);
				\draw[line width=0.25mm,dashed,blue] (ew) circle (0.25);
				\draw[<->] (0,5.5)to(0,6)to(2,6)to(2,5.5);
				\draw[<->] (2.5,4)to(3,4)to(3,0)to(2.5,0);
			\end{tikzpicture}
		}
		\caption{If the second player plays distinctly (dotted vertices), answer on the same component (dashed vertices), and associate the rows and columns.}
		\label{fig-tensorProductNandN-5}
	\end{subfigure}
	\caption{A depiction of the first player's strategy on $P_5 \times P_3$, the Tensor product of two \outcomeN graphs.}
	\label{fig-tensorProductNandN}
\end{figure}

\begin{proposition}
	Let $F = P_{5}$, $G = P_{3}$, and $H = S_{1,1,3}$, where $S_{1,1,3}$ denotes the spider tree with a degree-$3$ root having three branches of lengths $1,1$ and $3$. Then $F$, $G$, and $H$ are all outcome \outcomeN\ in the \AGAG, yet their tensor products satisfy
	\(
	F \times G \;\text{ is \outcomeN }
	\quad\text{while}\quad 
	G \times H \;\text{ is \outcomeP }.
	\)
	In particular, this shows that even if two bipartite graphs have the same outcome individually, their tensor product may produce a different outcome.
\end{proposition}

\begin{proof}
	\begin{figure}[ht]
		\centering
		\begin{tikzpicture}[scale=0.85, transform shape]
			\node at (0,0) {
				\begin{tikzpicture}
					\node (u0) at (0,5.5) {$G$};
					\node (u'0) at (0.5,6) {$H$};
					
					\node[noeud] (u1) at (1,6) {};
					\node[noeud] (u2) at (2,6) {};
					\node[noeud] (u3) at (3,6) {};
					
					\node[noeud] (v0) at (0,0) {};
					\node[noeud] (v1) at (0,1) {};
					\node[noeud] (v2) at (0,2) {};
					\node[noeud] (v3) at (0,3) {};
					\node[noeud] (v4) at (0,4) {};
					\node[noeud] (v5) at (0,5) {};
					
					\draw (u1)to(u2)to(u3);
					\draw (v0)to(v1)to(v2);
					\draw[bend left] (v2)to(v5);
					\draw[bend left] (v3)to(v5);
					\draw (v4)to(v5);
					
					\node[noeud] (w01) at (1,0) {};
					\node[noeud] (w02) at (2,0) {};
					\node[noeud] (w03) at (3,0) {};
					\node[noeud] (w11) at (1,1) {};
					\node[noeud] (w12) at (2,1) {};
					\node[noeud] (w13) at (3,1) {};
					\node[noeud] (w21) at (1,2) {};
					\node[noeud] (w22) at (2,2) {};
					\node[noeud] (w23) at (3,2) {};
					\node[noeud] (w31) at (1,3) {};
					\node[noeud] (w32) at (2,3) {};
					\node[noeud] (w33) at (3,3) {};
					\node[noeud] (w41) at (1,4) {};
					\node[noeud] (w42) at (2,4) {};
					\node[noeud] (w43) at (3,4) {};
					\node[noeud] (w51) at (1,5) {};
					\node[noeud] (w52) at (2,5) {};
					\node[noeud] (w53) at (3,5) {};
					
					\draw [red,thick] (w01)to(w12)to(w03);
					\draw [blue,thick](w11)to(w02)to(w13);
					\draw [blue,thick](w11)to(w22)to(w13);
					\draw [red,thick](w21)to(w12)to(w23);
					\draw [red,thick](w21)to(w52)to(w23);
					\draw [blue,thick](w51)to(w22)to(w53);
					\draw [red,thick](w31)to(w52)to(w33);
					\draw [blue,thick](w51)to(w32)to(w53);
					\draw[red,thick] (w41)to(w52)to(w43);
					\draw [blue,thick](w51)to(w42)to(w53);
				\end{tikzpicture}
			};
			
			\node at (8,0) {
				\begin{tikzpicture}[scale=0.85, transform shape]
					
					\node (u'_1) at (-0.3,5) {$u_1$};
					\node (u'_2) at (1,6.3) {$u_2$};
					\node (u'_3) at (2,5.5) {$u_3$};
					\node (u'_4) at (1,3.7) {$u_4$};
					\node (u'_5) at (0.6,4.8) {$u_5$};
					\node (u'_6) at (4.3,5) {$u_6$};
					\node (u'_7) at (3,3.7) {$u_7$};
					\node (u'_8) at (3,6.3) {$u_8$};
					\node[noeud] (c0) at (1,2) {};
					\node (c') at (-1,3) {$R$};
					\node[noeud] (c1) at (2,3) {};
					
					\node (c"1) at (2,3.3) {$v_6$};
					\node (c"2) at (1,2.3) {$v_5$};
					\node (c"3) at (2,0.7) {$v_8$};
					\node (c"4) at (3,2.3) {$v_7$};
					\node (c"5) at (4,3.3) {$v_9$};
					\node (c"6) at (4,0.7) {$v_10$};
					\node (c"7) at (-0.5,2.5) {$v_2$};
					\node (c"8) at (-0.5,3.5) {$v_1$};
					\node (c"9) at (-0.5,1.5) {$v_3$};
					\node (c"10) at (-0.5,0.5) {$v_4$};
					
					\node[noeud] (c2) at (3,2) {};
					\node[noeud] (c3) at (2,1) {};
					
					\draw [red,thick] (c0)to(c1)to(c2)to(c3)to(c0);
					\foreach \I in {0,1,2,3} {
						\node[noeud] (l\I) at (0,\I+0.5) {};
						\draw [red,thick] (l\I) to (c0);
					}
					\node[noeud] (m1) at (4,1) {};
					\node[noeud] (m3) at (4,3) {};
					\draw [red,thick](m1)to(c2)to(m3);
					
					\node[noeud] (d0) at (0,5) {};
					\node (d') at (-1,6) {B};
					\node[noeud] (d1) at (1,6) {};
					\node[noeud] (d2) at (2,5) {};
					\node[noeud] (d3) at (1,4) {};
					\node[noeud] (d4) at (1,5) {};
					\draw [blue,thick](d0)to(d1)to(d2)to(d3)to(d0);
					\draw [blue,thick](d1)to(d4)to(d3);
					\node[noeud] (e0) at (3,6) {};
					\node[noeud] (e1) at (4,5) {};
					\node[noeud] (e2) at (3,4) {};
					\draw [blue,thick](d2)to(e0)to(e1)to(e2)to(d2);
				\end{tikzpicture}
			};
			\draw[->,line width=0.25mm] (3,0) to (4.5,0);
		\end{tikzpicture} 
		\caption{The tensor product $G \times H$ which has $2$ components.}
		\label{fig-tensorProductNandN-diff}
	\end{figure}
	%%%%%%%%%%%%%%%%%%%%%%%%%%%%%%%%%%%%%%%%%%%%%%%%%%%%%%%%%%%%%%%%%%%%%%%%%%%%%%%%%%%%%%%%%%%%%
	
	\begin{figure}[ht]
		\centering
		\begin{tikzpicture}[scale=0.85, transform shape]
			\node (initial) at (0,1) {
				\begin{tikzpicture}[scale=0.85, transform shape]
					\node[noeud] (d0) at (0,5) {};
					\node[noeud] (d1) at (1,6) {};
					\node[noeud] (d2) at (2,5) {};
					\node[noeud] (d3) at (1,4) {};
					\node[noeud] (d4) at (1,5) {};
					\draw [blue, thick](d0)to(d1)to(d2)to(d3)to(d0);
					\draw [blue, thick](d1)to(d4)to(d3);
					\node[noeud] (e0) at (3,6) {};
					\node[noeud] (e1) at (4,5) {};
					\node[noeud] (e2) at (3,4) {};
					\draw[blue, thick] (d2)to(e0)to(e1)to(e2)to(d2);
					
					\node (u_1) at (-0.5,5) {$u_1$};
					\node (u_2) at (1,6.5) {$u_2$};
					\node (u_3) at (2,5.5) {$u_3$};
					\node (u_4) at (1,3.5) {$u_4$};
					\node (u_5) at (0.6,4.8) {$u_5$};
					\node (u_6) at (4.5,5) {$u_6$};
					\node (u_7) at (3,3.5) {$u_7$};
					\node (u_8) at (3,6.5) {$u_8$};
					
				\end{tikzpicture}
			};
			
			\node (opt1) at (-5.5,-2.5) {
				\begin{tikzpicture}[scale=0.85, transform shape]
					\node[noeud] (d0) at (0,5) {};
					\node[noeud] (d1) at (1,6) {};
					\node[noeud,fill=red] (d2) at (2,5) {};
					\node[noeud] (d3) at (1,4) {};
					\node[noeud] (d4) at (1,5) {};
					\draw [blue, thick](d0)to(d1)to(d2)to(d3)to(d0);
					\draw [blue, thick] (d1)to(d4)to(d3);
					\node[noeud,cross out] (e0) at (3,6) {};
					\node[noeud,fill=blue] (e1) at (4,5) {};
					\node[noeud,cross out] (e2) at (3,4) {};
					\draw [blue, thick](d2)to(e0)to(e1)to(e2)to(d2);
					
					\node (u_1) at (-0.5,5) {$u_1$};
					\node (u_2) at (1,6.5) {$u_2$};
					\node (u_3) at (2,5.5) {$u_3$};
					\node (u_4) at (1,3.5) {$u_4$};
					\node (u_5) at (0.6,4.8) {$u_5$};
					\node (u_6) at (4.5,5) {$u_6$};
					\node (u_7) at (3,3.5) {$u_7$};
					\node (u_8) at (3,6.5) {$u_8$};
					
					\draw[line width=0.25mm,dotted,red] (d2) circle (0.25);
					\draw[line width=0.25mm,dashed,blue] (e1) circle (0.25);
				\end{tikzpicture}
			};
			\node (opt2) at (0,-2.5) {
				\begin{tikzpicture}[scale=0.85, transform shape]
					\node[noeud] (d0) at (0,5) {};
					\node[noeud,cross out] (d1) at (1,6) {};
					\node[noeud,cross out] (d2) at (2,5) {};
					\node[noeud,cross out] (d3) at (1,4) {};
					\node[noeud,fill=blue] (d4) at (1,5) {};
					\draw [blue, thick](d0)to(d1)to(d2)to(d3)to(d0);
					\draw [blue, thick] (d1)to(d4)to(d3);
					\node[noeud,fill=red] (e0) at (3,6) {};
					\node[noeud] (e1) at (4,5) {};
					\node[noeud] (e2) at (3,4) {};
					\draw [blue, thick](d2)to(e0)to(e1)to(e2)to(d2);
					
					\node (u_1) at (-0.5,5) {$u_1$};
					\node (u_2) at (1,6.5) {$u_2$};
					\node (u_3) at (2,5.5) {$u_3$};
					\node (u_4) at (1,3.5) {$u_4$};
					\node (u_5) at (0.6,4.8) {$u_5$};
					\node (u_6) at (4.5,5) {$u_6$};
					\node (u_7) at (3,3.5) {$u_7$};
					\node (u_8) at (3,6.5) {$u_8$};
					
					\draw[line width=0.25mm,dashed,blue] (d4) circle (0.25);
					\draw[line width=0.25mm,dotted,red] (e0) circle (0.25);
				\end{tikzpicture}
			};
			
			\node (opt3) at (5.5,-2.5) {
				\begin{tikzpicture}[scale=0.85, transform shape]
					\node[noeud] (d0) at (0,5) {};
					\node[noeud,fill=red] (d1) at (1,6) {};
					\node[noeud,fill=blue] (d2) at (2,5) {};
					\node[noeud] (d3) at (1,4) {};
					\node[noeud] (d4) at (1,5) {};
					\draw [blue, thick](d0)to(d1)to(d2)to(d3)to(d0);
					\draw [blue, thick] (d1)to(d4)to(d3);
					\node[noeud] (e0) at (3,6) {};
					\node[noeud] (e1) at (4,5) {};
					\node[noeud] (e2) at (3,4) {};
					\draw [blue, thick](d2)to(e0)to(e1)to(e2)to(d2);
					
					\node (u_1) at (-0.5,5) {$u_1$};
					\node (u_2) at (1,6.5) {$u_2$};
					\node (u_3) at (2,5.5) {$u_3$};
					\node (u_4) at (1,3.5) {$u_4$};
					\node (u_5) at (0.6,4.8) {$u_5$};
					\node (u_6) at (4.5,5) {$u_6$};
					\node (u_7) at (3,3.5) {$u_7$};
					\node (u_8) at (3,6.5) {$u_8$};
					
					\draw[line width=0.25mm,dotted,red] (d1) circle (0.25);
					\draw[line width=0.25mm,dashed,blue] (d2) circle (0.25);
				\end{tikzpicture}
			};
			
			\foreach \I in {1,2,3}{\draw[->,double,line width=0.25mm] (initial)to(opt\I);}
			
			\node (opt4) at (0,-6) {
				\begin{tikzpicture}[scale=0.85, transform shape]
					\node[noeud,fill=red] (d0) at (0,5) {};
					\node[noeud,fill=red] (d1) at (1,6) {};
					\node[noeud,fill=blue] (d2) at (2,5) {};
					\node[noeud,cross out] (d3) at (1,4) {};
					\node[noeud] (d4) at (1,5) {};
					\draw [blue, thick](d0)to(d1)to(d2)to(d3)to(d0);
					\draw [blue, thick] (d1)to(d4)to(d3);
					\node[noeud,fill=blue] (e0) at (3,6) {};
					\node[noeud] (e1) at (4,5) {};
					\node[noeud] (e2) at (3,4) {};
					\draw [blue, thick](d2)to(e0)to(e1)to(e2)to(d2);
					
					\node (u_1) at (-0.5,5) {$u_1$};
					\node (u_2) at (1,6.5) {$u_2$};
					\node (u_3) at (2,5.5) {$u_3$};
					\node (u_4) at (1,3.5) {$u_4$};
					\node (u_5) at (0.6,4.8) {$u_5$};
					\node (u_6) at (4.5,5) {$u_6$};
					\node (u_7) at (3,3.5) {$u_7$};
					\node (u_8) at (3,6.5) {$u_8$};

					\draw[line width=0.25mm,dotted,red] (d0) circle (0.25);
					\draw[line width=0.25mm,dashed,blue] (e0) circle (0.25);
				\end{tikzpicture}
			};
			
			\node (opt5) at (5.5,-6) {
				\begin{tikzpicture}[scale=0.85, transform shape]
					\node[noeud,cross out] (d0) at (0,5) {};
					\node[noeud,fill=red] (d1) at (1,6) {};
					\node[noeud,fill=blue] (d2) at (2,5) {};
					\node[noeud,fill=blue] (d3) at (1,4) {};
					\node[noeud,cross out] (d4) at (1,5) {};
					\draw [blue, thick](d0)to(d1)to(d2)to(d3)to(d0);
					\draw [blue, thick] (d1)to(d4)to(d3);
					\node[noeud,cross out] (e0) at (3,6) {};
					\node[noeud,fill=red] (e1) at (4,5) {};
					\node[noeud,cross out] (e2) at (3,4) {};
					\draw [blue, thick](d2)to(e0)to(e1)to(e2)to(d2);
					
					\node (u_1) at (-0.5,5) {$u_1$};
					\node (u_2) at (1,6.5) {$u_2$};
					\node (u_3) at (2,5.5) {$u_3$};
					\node (u_4) at (1,3.5) {$u_4$};
					\node (u_5) at (0.6,4.8) {$u_5$};
					\node (u_6) at (4.5,5) {$u_6$};
					\node (u_7) at (3,3.5) {$u_7$};
					\node (u_8) at (3,6.5) {$u_8$};
					
					\draw[line width=0.25mm,dotted,red] (e1) circle (0.25);
					\draw[line width=0.25mm,dashed,blue] (d3) circle (0.25);
				\end{tikzpicture}
			};
			\foreach \I in {4,5}{\draw[->,double,line width=0.25mm] (opt3)to(opt\I);}
		\end{tikzpicture}
		\caption{Illustration of the second player’s winning strategy when the first player selects a vertex in component B.
		}
		\label{fig-6}
	\end{figure}

	By Proposition~\ref{P_ntimes P_m}, the tensor product $F \times G = P_{5} \times P_{3}$ is $\mathcal{N}$, so the first player has a winning strategy. It remains to show that $G \times H = P_{3} \times S_{1,1,3}$ is $\mathcal{P}$.
	
	Since both $P_{3}$ and $S_{1,1,3}$ are connected bipartite graphs, Lemma~\ref{tensorlemma2} implies that their tensor product decomposes into two connected components. It suffices to show that each component is an $\mathcal{P}$-position. Let $B$ and $R$ denote these two components.
	
	First, we prove that $B$ has a $\mathcal{P}$-position. The graph $B$ contains a unique articulation vertex $u_3$ such that $B \setminus \{u_3\}$ splits into two induced subgraphs:
	\begin{itemize}
		\item a cycle $C$ on vertices $\{u_1,u_2,u_3,u_4,u_5\}$, and
		\item a path $P$ on vertices $\{u_6,u_7,u_8\}$ attached at $u_3$.
	\end{itemize}
	where $C$ induces a cycle and $P$ induces a path attached to $u_3$. In particular, every path between a vertex of $C$ and a vertex of $P$ passes through $u_3$.
	
	We show that the second player has a winning strategy by responding to any first move so that, after the first two moves, the remaining position requires an even number of moves to terminate.
	
	\smallskip\noindent
	We describe a winning strategy:
	\smallskip\noindent
	We now verify that in each case the position reduces to a forced path.
	
	\begin{itemize}
		\item [(i)] If the first player selects $u_3$, the second player selects $u_6$. Then the shortest paths between $u_3$ and $u_6$ cover $u_7$ and $u_8$. The remaining uncovered vertices are $\{u_1,u_2,u_4,u_5\}$, which form a cycle. Since, any two selected vertices on this cycle cover all these vertices, the number of remaining moves is 2, in which the second player must make the last move. 
		
		\item [(ii)] If the first player selects a vertex in $P$, say $u_8$ (the case $u_7$ is symmetric), the second player selects a vertex in $C$ at distance $2$ from $u_3$, for instance $u_1$ or $u_5$. Suppose the first player selects $u_8$ (the case $u_7$ is similar), and the second player selects $u_1$.
		
		Every shortest path from $u_8$ to $u_1$ must pass through $u_3$. Moreover, there are two shortest paths from $u_8$ to $u_1$, namely
		\(
		u_8 - u_3 - u_2 - u_1
		\quad \text{and} \quad
		u_8-u_3 - u_4 - u_1,
		\) so the vertices $\{u_1,u_2,u_3,u_4,u_8\}$ are contained in the geodetic closure. Thus, the remaining uncovered vertices are $\{u_5,u_6,u_7\}$. Observe that: among these, the vertex $u_5$ does not lie on any shortest path between already selected vertices, hence it must be selected explicitly in a future move and the remaining vertices $\{u_6,u_7\}$ require exactly one additional move to be fully covered.
		Therefore, exactly two further moves are required to finish the game, and hence the second player makes the last move.
		
		\item [(iii)]
		Suppose the first player selects $u_2$ (the case $u_4$ is symmetric), and the second player selects $u_3$.
		Since $u_2$ and $u_3$ are adjacent, no additional vertices are covered at this stage. Moreover, because $u_3$ is an articulation vertex, any subsequent move in the path component must interact through $u_3$, and similarly any move in the cycle remains confined within the cycle.
		Hence the remaining play decomposes into independent play on:
		\begin{itemize}
			\item the cycle $C_4$ with one vertex already selected, and
			\item the path $P_3= \{u_6,u_7,u_8\}$ attached at $u_3$.
		\end{itemize}
		Since the players play the optimal moves, in the cycle component, selecting any remaining vertex forces coverage in at most two additional moves, while the path component requires also at most two move to be fully covered. Thus the total number of remaining moves is fixed and even. Therefore, the second player makes the last move.
		
	\end{itemize}
	
	\smallskip
	In all cases, the second player ensures that the remaining position has an even number of moves. Therefore, $(B,\emptyset)$ is an $\mathcal{P}$-position.
	
	\medskip
	
	Now we consider that the first player starts in component $R$.
	The graph $R$ contains two articulation vertices, namely $v_5$ and $v_7$, together with pendant vertices  $v_1,\ldots,v_4$ to $v_5$ and $v_9,v_{10}$ attached to $v_7$. 
	
	If the first player selects one of the articulation vertices, say $v_5$, then the second player selects $v_7$. The geodetic closure of $\{v_5,v_7\}$ contains all vertices on shortest paths between them, and the remaining vertices are pendant vertices. These vertices can be covered in pairs, so the total number of remaining moves is even.
	
	If the first player selects a pendant vertex, then the second player selects a pendant vertex attached to the other articulation vertex, thereby forcing both articulation vertices and the vertices on the connecting paths into the closure. The remaining vertices again form disjoint pendant structures that require an even number of moves.
	
	If the first player selects a non-pendant, non-articulation vertex, then the second player responds at a symmetric vertex so that the resulting position reduces to one of the previous cases. In each situation, the second player ensures that the number of remaining moves is even, and hence makes the final move. Therefore $(R,\emptyset)$ is an $\mathcal{P}$-position.
	
	\medskip\noindent Since both components $B$ and $R$ are $\mathcal{P}$-positions, it follows that $G \times H$ is an $\mathcal{P}$-position.
	
\end{proof}

\section{\lovesme situations}
\label{sec-lovesMeLovesMeNot}

\lovesme is one of the simplest combinatorial games: the two players alternate removing blossoms from a flower until it is empty. It has become a shorthand for easy, binary combinatorial games where every possible move from the initial position will be played. Note that, in those cases, the Sprague-Grundy values alternate between~0 (for an even number of moves) and~1 (for an odd number of moves), since no real strategy is involved. Some complex combinatorial games can be equivalent to \lovesme depending on the initial position. It is in particular the case for \AGAG in some cases:

\begin{proposition}
	\label{prop-activeGeodeticSheLovesMe}
	A game of \AGAG on the complete graph $K_n$ or a star $K_{1,n}$ will see all its vertices selected.
\end{proposition}

\begin{proof}
	Since every vertex in $K_n$ is simplicial, each vertex of $K_n$ will be selected in \AGAG (\Cref{lem-simplicial}). Hence, $K_n$ is a \lovesme, and thus $\gr(K_n)=n \bmod 2$.
	
	For a star $K_{1,n}$ with center $c$ and $n$ leaves, note that every leaf will need to be selected (each leaf is a simplicial vertex). If $n$ is even, then, the first player selects $c$ and win (being the second player with an even number of remaining moves). If $n$ is odd, then, there are two possible moves for the first player: either selecting $c$, in which case every vertex will be selected and they will lose, or selecting a leaf, in which case the second player selects $c$ and leave the first player in a losing position where every vertex will be selected. Hence, $K_{1,n}$ is a \lovesme, and thus $\gr(K_{1,n}) = 1-(n \bmod 2)$.
\end{proof}

Note that those two graphs have a unique minimum-size geodetic set, so one could think that such parity results could be obtained on any graph family with a unique, well-identified minimum-size geodetic set on an underlying structure. However, such a general structural result would be difficult to achieve, as illustrated by \Cref{fig-trianglesWithLeaves}: those two graphs are based on the same underlying structure (a triangle with leaves), have the same order and unique minimum-size geodetic sets of the same size, but different outcomes for \AGAG.
The detail is shown in \Cref{fig-trianglesWithLeaves-1-detail,fig-trianglesWithLeaves-2-detail}, and uses a \emph{pairing strategy}. A \emph{paired position} is a position on which the second player (from this position, so they can be the first player in the general game) can \emph{pair} vertices together: if the first player plays on one vertex of the pair, then the second player plays on the other vertex of the pair. Sometimes, several vertices can be paired with the same vertex, in which case they are identified as one for the purpose of the pair. This guarantees that the second player can always answer to the first player's move, and as such wins on a paired position.

\begin{figure}[h]
	\centering
	\begin{subfigure}[b]{0.4\textwidth}
		\centering
		\begin{tikzpicture}
			\node[noeud] (i1) at (0,0) {};
			\node[noeud] (i2) at (1,0) {};
			\node[noeud] (i3) at (0.5,0.5) {};
			\node[noeud] (l1) at (0,1) {};
			\node[noeud] (l2) at (0.5,1) {};
			\node[noeud] (l3) at (1,1) {};
			\node[noeud] (l4) at (-0.5,0) {};
			\node[noeud] (l5) at (1.5,0) {};
			\draw (l4)to(i1)to(i2)to(l5);
			\draw (l1)to(i3)to(i1);
			\draw (l3)to(i3)to(i2);
			\draw (l2)to(i3);
		\end{tikzpicture}
		\caption{On this graph, the second player wins, see \Cref{fig-trianglesWithLeaves-1-detail}.}
		\label{fig-trianglesWithLeaves-1}
	\end{subfigure}
	\hfil
	\begin{subfigure}[b]{0.5\textwidth}
		\centering
		\begin{tikzpicture}
			\node[noeud,fill=black] (i1) at (0,0) {};
			\node[noeud] (i2) at (1,0) {};
			\node[noeud] (i3) at (0.5,0.5) {};
			\node[noeud] (l1) at (0.5,1) {};
			\node[noeud] (l2) at (-0.5,-0.25) {};
			\node[noeud] (l3) at (-0.5,0.25) {};
			\node[noeud] (l4) at (1.5,-0.25) {};
			\node[noeud] (l5) at (1.5,0.25) {};
			\draw (l2)to(i1)to(i2)to(l4);
			\draw (l3)to(i1)to(i3)to(l1);
			\draw (i3)to(i2)to(l5);
		\end{tikzpicture}
		\caption{On this graph, the first player wins by picking the black vertex, see \Cref{fig-trianglesWithLeaves-2-detail}.}
		\label{fig-trianglesWithLeaves-2}
	\end{subfigure}
	\caption{Two graphs of the same order (8) with unique minimum-size geodetic sets of same size (5) but different outcomes for \AGAG.}
	\label{fig-trianglesWithLeaves}
\end{figure}

\begin{figure}[h]
	\centering
	\begin{tikzpicture}
		\node (initial) at (0,0) {
			\begin{tikzpicture}
				\node[noeud] (i1) at (0,0) {};
				\node[noeud] (i2) at (1,0) {};
				\node[noeud] (i3) at (0.5,0.5) {};
				\node[noeud] (l1) at (0,1) {};
				\node[noeud] (l2) at (0.5,1) {};
				\node[noeud] (l3) at (1,1) {};
				\node[noeud] (l4) at (-0.5,0) {};
				\node[noeud] (l5) at (1.5,0) {};
				\draw (l4)to(i1)to(i2)to(l5);
				\draw (l1)to(i3)to(i1);
				\draw (l3)to(i3)to(i2);
				\draw (l2)to(i3);
			\end{tikzpicture}
		};
		
		\node (opt1) at (-4.5,-2.5) {
			\begin{tikzpicture}
				\node[noeud] (i1) at (0,0) {};
				\node[noeud] (i2) at (1,0) {};
				\node[noeud] (i3) at (0.5,0.5) {};
				\node[noeud] (l1) at (0,1) {};
				\node[noeud,fill=black] (l2) at (0.5,1) {};
				\node[noeud] (l3) at (1,1) {};
				\node[noeud] (l4) at (-0.5,0) {};
				\node[noeud] (l5) at (1.5,0) {};
				\draw (l4)to(i1)to(i2)to(l5);
				\draw (l1)to(i3)to(i1);
				\draw (l3)to(i3)to(i2);
				\draw (l2)to(i3);
			\end{tikzpicture}
		};
		\node (opt2) at (-1.5,-2.5) {
			\begin{tikzpicture}
				\node[noeud] (i1) at (0,0) {};
				\node[noeud] (i2) at (1,0) {};
				\node[noeud,fill=black] (i3) at (0.5,0.5) {};
				\node[noeud] (l1) at (0,1) {};
				\node[noeud] (l2) at (0.5,1) {};
				\node[noeud] (l3) at (1,1) {};
				\node[noeud] (l4) at (-0.5,0) {};
				\node[noeud] (l5) at (1.5,0) {};
				\draw (l4)to(i1)to(i2)to(l5);
				\draw (l1)to(i3)to(i1);
				\draw (l3)to(i3)to(i2);
				\draw (l2)to(i3);
			\end{tikzpicture}
		};
		
		\node (ans1) at (-3,-5) {
			\begin{tikzpicture}
				\node[noeud] (i1) at (0,0) {};
				\node[noeud] (i2) at (1,0) {};
				\node[noeud,fill=black] (i3) at (0.5,0.5) {};
				\node[noeud] (l1) at (0,1) {};
				\node[noeud,fill=black] (l2) at (0.5,1) {};
				\node[noeud] (l3) at (1,1) {};
				\node[noeud] (l4) at (-0.5,0) {};
				\node[noeud] (l5) at (1.5,0) {};
				\draw (l4)to(i1)to(i2)to(l5);
				\draw (l1)to(i3)to(i1);
				\draw (l3)to(i3)to(i2);
				\draw (l2)to(i3);
				
				\draw[<->,bend right=60] (i1)to(i2);
				\draw[<->,bend right=60] (l4)to(l5);
				\draw[<->,bend right=60] (l3)to(l1);
			\end{tikzpicture}
		};
		
		\node (opt3) at (1.5,-2.5) {
			\begin{tikzpicture}
				\node[noeud,fill=black] (i1) at (0,0) {};
				\node[noeud] (i2) at (1,0) {};
				\node[noeud] (i3) at (0.5,0.5) {};
				\node[noeud] (l1) at (0,1) {};
				\node[noeud] (l2) at (0.5,1) {};
				\node[noeud] (l3) at (1,1) {};
				\node[noeud] (l4) at (-0.5,0) {};
				\node[noeud] (l5) at (1.5,0) {};
				\draw (l4)to(i1)to(i2)to(l5);
				\draw (l1)to(i3)to(i1);
				\draw (l3)to(i3)to(i2);
				\draw (l2)to(i3);
			\end{tikzpicture}
		};
		\node (opt4) at (4.5,-2.5) {
			\begin{tikzpicture}
				\node[noeud] (i1) at (0,0) {};
				\node[noeud] (i2) at (1,0) {};
				\node[noeud] (i3) at (0.5,0.5) {};
				\node[noeud] (l1) at (0,1) {};
				\node[noeud] (l2) at (0.5,1) {};
				\node[noeud] (l3) at (1,1) {};
				\node[noeud,fill=black] (l4) at (-0.5,0) {};
				\node[noeud] (l5) at (1.5,0) {};
				\draw (l4)to(i1)to(i2)to(l5);
				\draw (l1)to(i3)to(i1);
				\draw (l3)to(i3)to(i2);
				\draw (l2)to(i3);
			\end{tikzpicture}
		};
		
		\node (ans2) at (3,-5) {
			\begin{tikzpicture}
				\node[noeud,fill=black] (i1) at (0,0) {};
				\node[noeud] (i2) at (1,0) {};
				\node[noeud] (i3) at (0.5,0.5) {};
				\node[noeud] (l1) at (0,1) {};
				\node[noeud] (l2) at (0.5,1) {};
				\node[noeud] (l3) at (1,1) {};
				\node[noeud,fill=black] (l4) at (-0.5,0) {};
				\node[noeud] (l5) at (1.5,0) {};
				\draw (l4)to(i1)to(i2)to(l5);
				\draw (l1)to(i3)to(i1);
				\draw (l3)to(i3)to(i2);
				\draw (l2)to(i3);
				
				\draw[dashed] \convexpath{l1,l2,l3}{0.2cm};
				
				\draw[<->,bend right=60] (l5)to(1.2,1);
				\draw[<->,bend right=60] (i2)to(i3);
			\end{tikzpicture}
		};
		
		\foreach \I in {1,2,3,4}{\draw[->,double,line width=0.25mm] (initial)to(opt\I);}
		\foreach \I in {1,2}{\draw[->,double,line width=0.25mm] (opt\I)to(ans1);}
		\foreach \I in {3,4}{\draw[->,double,line width=0.25mm] (opt\I)to(ans2);}
	\end{tikzpicture}
	\caption{The second player wins on the graph of \Cref{fig-trianglesWithLeaves-1}, as seen with the possible (non-isomorphic) options of the first player, to which they can always answer to a \emph{paired position}. The paired vertices are linked with double arrows: if the first player picks one vertex of a given pair, the second player picks the other vertex of the same pair (grouped vertices are considered as one for the sake of a pair).}
	\label{fig-trianglesWithLeaves-1-detail}
\end{figure}

\begin{figure}
	\centering
	\begin{tikzpicture}
		\node (0) at (0,2.5) {
			\begin{tikzpicture}
				\node[noeud] (i1) at (0,0) {};
				\node[noeud] (i2) at (1,0) {};
				\node[noeud] (i3) at (0.5,0.5) {};
				\node[noeud] (l1) at (0.5,1) {};
				\node[noeud] (l2) at (-0.5,-0.25) {};
				\node[noeud] (l3) at (-0.5,0.25) {};
				\node[noeud] (l4) at (1.5,-0.25) {};
				\node[noeud] (l5) at (1.5,0.25) {};
				\draw (l2)to(i1)to(i2)to(l4);
				\draw (l3)to(i1)to(i3)to(l1);
				\draw (i3)to(i2)to(l5);
			\end{tikzpicture}
		};
		
		\node (initial) at (0,0) {
			\begin{tikzpicture}
				\node[noeud,fill=black] (i1) at (0,0) {};
				\node[noeud] (i2) at (1,0) {};
				\node[noeud] (i3) at (0.5,0.5) {};
				\node[noeud] (l1) at (0.5,1) {};
				\node[noeud] (l2) at (-0.5,-0.25) {};
				\node[noeud] (l3) at (-0.5,0.25) {};
				\node[noeud] (l4) at (1.5,-0.25) {};
				\node[noeud] (l5) at (1.5,0.25) {};
				\draw (l2)to(i1)to(i2)to(l4);
				\draw (l3)to(i1)to(i3)to(l1);
				\draw (i3)to(i2)to(l5);
			\end{tikzpicture}
		};
		
		\node (opt1) at (-6,-2.5) {
			\begin{tikzpicture}
				\node[noeud,fill=black] (i1) at (0,0) {};
				\node[noeud] (i2) at (1,0) {};
				\node[noeud,cross out] (i3) at (0.5,0.5) {};
				\node[noeud,fill=black] (l1) at (0.5,1) {};
				\node[noeud] (l2) at (-0.5,-0.25) {};
				\node[noeud] (l3) at (-0.5,0.25) {};
				\node[noeud] (l4) at (1.5,-0.25) {};
				\node[noeud] (l5) at (1.5,0.25) {};
				\draw (l2)to(i1)to(i2)to(l4);
				\draw (l3)to(i1)to(i3)to(l1);
				\draw (i3)to(i2)to(l5);
			\end{tikzpicture}
		};
		\node (opt2) at (-3,-2.5) {
			\begin{tikzpicture}
				\node[noeud,fill=black] (i1) at (0,0) {};
				\node[noeud,fill=black] (i2) at (1,0) {};
				\node[noeud] (i3) at (0.5,0.5) {};
				\node[noeud] (l1) at (0.5,1) {};
				\node[noeud] (l2) at (-0.5,-0.25) {};
				\node[noeud] (l3) at (-0.5,0.25) {};
				\node[noeud] (l4) at (1.5,-0.25) {};
				\node[noeud] (l5) at (1.5,0.25) {};
				\draw (l2)to(i1)to(i2)to(l4);
				\draw (l3)to(i1)to(i3)to(l1);
				\draw (i3)to(i2)to(l5);
			\end{tikzpicture}
		};
		
		\node (ans1) at (-4.5,-5) {
			\begin{tikzpicture}
				\node[noeud,fill=black] (i1) at (0,0) {};
				\node[noeud,fill=black] (i2) at (1,0) {};
				\node[noeud,cross out] (i3) at (0.5,0.5) {};
				\node[noeud,fill=black] (l1) at (0.5,1) {};
				\node[noeud] (l2) at (-0.5,-0.25) {};
				\node[noeud] (l3) at (-0.5,0.25) {};
				\node[noeud] (l4) at (1.5,-0.25) {};
				\node[noeud] (l5) at (1.5,0.25) {};
				\draw (l2)to(i1)to(i2)to(l4);
				\draw (l3)to(i1)to(i3)to(l1);
				\draw (i3)to(i2)to(l5);
			\end{tikzpicture}
		};
		
		\node (opt3) at (3,-2.5) {
			\begin{tikzpicture}
				\node[noeud,fill=black] (i1) at (0,0) {};
				\node[noeud] (i2) at (1,0) {};
				\node[noeud,fill=black] (i3) at (0.5,0.5) {};
				\node[noeud] (l1) at (0.5,1) {};
				\node[noeud] (l2) at (-0.5,-0.25) {};
				\node[noeud] (l3) at (-0.5,0.25) {};
				\node[noeud] (l4) at (1.5,-0.25) {};
				\node[noeud] (l5) at (1.5,0.25) {};
				\draw (l2)to(i1)to(i2)to(l4);
				\draw (l3)to(i1)to(i3)to(l1);
				\draw (i3)to(i2)to(l5);
			\end{tikzpicture}
		};
		\node (opt4) at (6,-2.5) {
			\begin{tikzpicture}
				\node[noeud,fill=black] (i1) at (0,0) {};
				\node[noeud,cross out] (i2) at (1,0) {};
				\node[noeud] (i3) at (0.5,0.5) {};
				\node[noeud] (l1) at (0.5,1) {};
				\node[noeud] (l2) at (-0.5,-0.25) {};
				\node[noeud] (l3) at (-0.5,0.25) {};
				\node[noeud] (l4) at (1.5,-0.25) {};
				\node[noeud,fill=black] (l5) at (1.5,0.25) {};
				\draw (l2)to(i1)to(i2)to(l4);
				\draw (l3)to(i1)to(i3)to(l1);
				\draw (i3)to(i2)to(l5);
			\end{tikzpicture}
		};
		
		\node (ans2) at (4.5,-5) {
			\begin{tikzpicture}
				\node[noeud,fill=black] (i1) at (0,0) {};
				\node[noeud,cross out] (i2) at (1,0) {};
				\node[noeud,fill=black] (i3) at (0.5,0.5) {};
				\node[noeud] (l1) at (0.5,1) {};
				\node[noeud] (l2) at (-0.5,-0.25) {};
				\node[noeud] (l3) at (-0.5,0.25) {};
				\node[noeud] (l4) at (1.5,-0.25) {};
				\node[noeud,fill=black] (l5) at (1.5,0.25) {};
				\draw (l2)to(i1)to(i2)to(l4);
				\draw (l3)to(i1)to(i3)to(l1);
				\draw (i3)to(i2)to(l5);
			\end{tikzpicture}
		};
		
		\node (opt5) at (0,-2.5) {
			\begin{tikzpicture}
				\node[noeud,fill=black] (i1) at (0,0) {};
				\node[noeud] (i2) at (1,0) {};
				\node[noeud] (i3) at (0.5,0.5) {};
				\node[noeud] (l1) at (0.5,1) {};
				\node[noeud] (l2) at (-0.5,-0.25) {};
				\node[noeud,fill=black] (l3) at (-0.5,0.25) {};
				\node[noeud] (l4) at (1.5,-0.25) {};
				\node[noeud] (l5) at (1.5,0.25) {};
				\draw (l2)to(i1)to(i2)to(l4);
				\draw (l3)to(i1)to(i3)to(l1);
				\draw (i3)to(i2)to(l5);
			\end{tikzpicture}
		};
		
		\node (ans3) at (0,-5) {
			\begin{tikzpicture}
				\node[noeud,fill=black] (i1) at (0,0) {};
				\node[noeud] (i2) at (1,0) {};
				\node[noeud] (i3) at (0.5,0.5) {};
				\node[noeud] (l1) at (0.5,1) {};
				\node[noeud,fill=black] (l2) at (-0.5,-0.25) {};
				\node[noeud,fill=black] (l3) at (-0.5,0.25) {};
				\node[noeud] (l4) at (1.5,-0.25) {};
				\node[noeud] (l5) at (1.5,0.25) {};
				\draw (l2)to(i1)to(i2)to(l4);
				\draw (l3)to(i1)to(i3)to(l1);
				\draw (i3)to(i2)to(l5);
				
				\draw[dashed] \convexpath{l4,l5}{0.2cm};
				
				\draw[<->,bend right=60] (i2)to(l1);
				\draw[<->,bend right=60] (1.5,0.45)to(i3);
			\end{tikzpicture}
		};
		
		\draw[->,double,line width=0.25mm] (0)to(initial);
		\foreach \I in {1,2,3,4,5}{\draw[->,double,line width=0.25mm] (initial)to(opt\I);}
		\foreach \I in {1,2}{\draw[->,double,line width=0.25mm] (opt\I)to(ans1);}
		\foreach \I in {3,4}{\draw[->,double,line width=0.25mm] (opt\I)to(ans2);}
		\draw[->,double,line width=0.25mm] (opt5)to(ans3);
	\end{tikzpicture}
	\caption{The first player wins on the graph of \Cref{fig-trianglesWithLeaves-2} by the move shown here. For each possible option of the second player starting from this position, they can always answer either to a \emph{paired position} or to a position where the only possible moves are leaves, of which an even number remains.}
	\label{fig-trianglesWithLeaves-2-detail}
\end{figure}

For other graphs, the game is not fully reducible to a \lovesme, but is still close since, under optimal play, we can determine exactly which vertices will be selected.

\begin{proposition}
	\label{prop-K2N}
	Let $n \geq 2$. A game of \AGAG on the complete bipartite graph $K_{2,n}$ will see all the vertices of its larger part selected.
\end{proposition}

\begin{proof}
	Selecting the two vertices in the part of size~2 will end the game (the whole other part will be in the geodetic closure). Hence, no player wants to select the first of those two vertices: by doing so, they would let their opponent win. Thus, the losing player will always select vertices in the larger part, and the winning player also will, until the larger part has been wholly selected (note that selecting two vertices in one part prevents from selecting vertices in the other part, ensuring that once the first player picks a part, the second player will be able to enforce that it remains in that part). Doing so will also end the game.
\end{proof}

We can also compute the Sprague-Grundy values of complete bipartite graphs, for which one part will always be fully selected:

\begin{proposition}
	\label{prop-Kmn}
	Let $m$ and $n$ be two integers such that $m,n \geq 2$. If $m$ and $n$ have the same parity, then, $\gr(K_{m,n})=0$. Otherwise, $\gr(K_{m,n})=2$.
\end{proposition}

\begin{proof}
	Note that, as for the proposition above, selecting two vertices in one part prevents from selecting vertices in the other part, and thus leaves a fixed number of moves that will all be played (selecting all the vertices in the part where two vertices have been selected).
	
	Assume first that $m$ and $n$ are both even. The first player selects a vertex in a part, and the second player selects a vertex in the same part, ensuring that an even number of moves will remain. Hence, $K_{m,n}$ is \outcomeP and thus $\gr(K_{m,n})=0$.
	
	Assume now that $m$ and $n$ are both odd. The first player selects a vertex in a part, and the second player selects a vertex in the other part. Now, whichever part the first player selects a vertex in, there will be an odd number of moves remaining, thus the second player will win. Hence, $K_{m,n}$ is \outcomeP and thus $\gr(K_{m,n})=0$.
	
	Assume finally (without loss of generality) that $m$ is even and $n$ is odd. Let us analyze the options of the first player.
	\begin{enumerate}
		\item Selecting a vertex $u_{odd}$ in the odd part. In this case, the second player has two answers:
		\begin{enumerate}
			\item selecting another vertex $v_{odd}$ in the odd part leaves an odd number of moves, and thus $\gr(K_{m,n},\{u_{odd},v_{odd}\})=1$;
			\item selecting a vertex $u_{even}$ in the even part, in which case there are again two options for the first player: selecting a vertex $v_{odd}$ in the odd part leaves an odd number of moves, so $(K_{m,n},\{u_{odd},u_{even},v_{odd}\})$ is \outcomeN and thus $\gr(K_{m,n},\{u_{odd},u_{even},v_{odd}\})=1$, while selecting a vertex $v_{even}$ in the even part leaves an even number of moves, so $(K_{m,n},\{u_{odd},u_{even},v_{even}\})$ is \outcomeP and thus $\gr(K_{m,n},\{u_{odd},u_{even},v_{even}\})=0$.
			
			Hence, $\gr(K_{m,n},\{u_{odd},u_{even}\})=\mex(\{0,1\})=2$.
		\end{enumerate}
		This implies that $\gr(K_{m,n},\{u_{odd}\})=\mex(\{1,2\})=0$.
		
		\item Selecting a vertex $u_{even}$ in the even part. In this case, the second player has two answers:
		\begin{enumerate}
			\item selecting another vertex $v_{even}$ in the even part leaves an even number of moves, and thus $\gr(K_{m,n},\{u_{even},v_{even}\})=0$;
			\item selecting a vertex $u_{odd}$ in the odd part, and as seen above $\gr(K_{m,n},\{u_{odd},u_{even}\})=2$.
		\end{enumerate}
		This implies that $\gr(K_{m,n},\{u_{even}\})=\mex(\{0,2\})=1$.
	\end{enumerate}
	Altogether, we have $\gr(K_{m,n})=\mex(\{0,1\})=2$.
\end{proof}

\section{Symmetry strategies}
\label{sec-symmetry}

The outcome of cycles can be fully derived from a symmetry strategy. Furthermore, we can also characterize the Sprague-Grundy value of a cycle:

\begin{theorem}
	\label{thm-cycles}
	For any positive integer $n$, $\gr(C_n)=n \bmod 2$.
\end{theorem}

\begin{proof}
	Denote the vertices of the cycle by $u_0,\ldots,u_{n-1}$, in cyclic order.
	
	If $n=2k$, then the second player has a winning strategy: when the first player selects a vertex $u_i$, they can select the opposite vertex (which is either $u_{i+k}$ or $u_{i-k}$), which ends the game. Hence, $\gr(C_n)=0$.
	
	Now, let $n=2k+1$. The first player selects a vertex, say $u_0$. Note that, for $i \in \{0,\ldots,k-1\}$, selecting $u_{k+1+i}$ is equivalent as selecting $u_{k-i}$ for the second move, as can be seen on \Cref{fig-equivalentMovesOnOddCycles}. Now, assume that the second player selects a vertex $u_{i}$ for $i \in \{1,\ldots,k\}$, then the first player can select $u_{k+i}$, which ends the game. Hence, $\gr(C_n)=1$ (since all possible first moves are equivalent, up to renaming the vertices, there is only one option for the first player).
\end{proof}

\begin{figure}[h]
	\centering
	\begin{tikzpicture}
		\node (c1) at (0,0) {
			\begin{tikzpicture}
				\draw (0,0) circle (1.5);
				\foreach \I in {0,4} {
					\node[noeud,fill=black] (\I) at (\I*360/11:1.5) {};
					\draw (\I*360/11:1.875) node {$u_{\I}$};
				}
				\foreach \I in {1,2,3} {
					\node[noeud,cross out] (\I) at (\I*360/11:1.5) {};
					\draw (\I*360/11:1.875) node {$u_{\I}$};
				}
				\foreach \I in {5,...,10} {
					\node[noeud] (\I) at (\I*360/11:1.5) {};
					\draw (\I*360/11:1.875) node {$u_{\I}$};
				}
			\end{tikzpicture}
		};
		\node (c2) at (6,0) {
			\begin{tikzpicture}
				\draw (0,0) circle (1.5);
				\foreach \I in {0,7} {
					\node[noeud,fill=black] (\I) at (\I*360/11:1.5) {};
					\draw (\I*360/11:1.875) node {$u_{\I}$};
				}
				\foreach \I in {8,9,10} {
					\node[noeud,cross out] (\I) at (\I*360/11:1.5) {};
					\draw (\I*360/11:1.875) node {$u_{\I}$};
				}
				\foreach \I in {1,...,6} {
					\node[noeud] (\I) at (\I*360/11:1.5) {};
					\draw (\I*360/11:1.875) node {$u_{\I}$};
				}
			\end{tikzpicture}
		};
		\draw (3,0) node {$\equiv$};
	\end{tikzpicture}
	\caption{Selecting $u_{k-i}$ and $u_{k+1+i}$ are equivalent moves on $C_{2k+1},\{u_0\}$ (here, with $k=5$ and $i=1$).}
	\label{fig-equivalentMovesOnOddCycles}
\end{figure}

We can also characterize the Sprague-Grundy value of a cycle with one or two selected vertices:

\begin{proposition}
	\label{prop-cyclesOneVertex}
	For any positive integer $k$:
	\begin{itemize}
		\item $\gr(C_{2k},\{u\})=k$ and $\gr(C_{2k+1},\{u\})=0$ for any vertex $u$;
		\item $\gr(C_{2k+1},\{u_0,u_i\}) = k+1-i$ for $i \in \{1,\ldots,k\}$ (denoting the vertices by \\$u_0,\ldots,u_{2k}$ in the cyclic order).
	\end{itemize}
\end{proposition}

\begin{proof}
	Let $k$ be a positive integer. Note that, for the first move, all the vertices are equivalent. Clearly, since $\gr(C_{2k+1})=1$, we have $\gr(C_{2k+1},\{u\})=0$. For $C_{2k}$ ($k \geq 2$), we have to analyze the moves of the second player. Denote the vertices by $u_0,\ldots,u_{2k-1}$ in cyclic order, and assume, without loss of generality, that the first player selected $u_0$. Note that selecting $u_{k+i}$ and $u_{k-i}$ are equivalent moves for $i \in \{1,\ldots,k-1\}$, and that selecting $u_{k}$ is a winning move (and thus $\gr(C_{2k},\{u_0,u_{k}\})=0$). We can prove by induction on $i$ that $\gr(C_{2k},\{u_0,u_{k-i}\})=i$. If $i=1$, then, any move of the first player will end the game, so $\gr(C_{2k},\{u_0,u_{k-1}\})=1$. Now, for $i \geq 2$, assume that the second player selected $u_{k-i}$. The first player can select $u_{k-j}$ for $j < i$, so there are options with Sprague-Grundy values $1,\ldots,i-1$ by induction hypothesis, and they can select $u_{k}$, which ends the game, so $\gr(C_{2k},\{u_0,u_{k-i}\}) \geq i$. Furthermore, the first player can select $u_{j}$ for $j \in \{k+1,\ldots,2k-i\}$; those options end the game. Finally, the first player can select $u_{2k-i+j}$ for $j \in \{1,\ldots,i-1\}$; and we can clearly by a symmetry argument see that $\gr(C_{2k},\{u_0,u_{k-i},u_{2k-i+j}\})=\gr(C_{2k},\{u_0,u_{k-j}\})=j$ by induction hypothesis.
	Hence, $\gr(C_{2k},\{u_0,u_{k-i}\})=\mex(\{0,1,\ldots,i-1\})=i$ for $i \in \{1,\ldots,k-1\}$, which implies that $\gr(C_{2k},\{u\})=\mex(\{0,\ldots,k-1\})=k$. We can apply a similar reasoning to obtain $\gr(C_{2k+1},\{u_0,u_i\}) = k+1-i$ for $i \in \{1,\ldots,k\}$.
\end{proof}

\Cref{thm-cartesianProductOfNandN} combined with \Cref{thm-paths} already give us the outcome for multidimensional grids, but we can also prove it with a symmetry strategy:

\begin{theorem}
	\label{thm-grids}
	\AGAG has outcome \outcomeN for a multidimensional grid if and only if all its dimensions are odd.
\end{theorem}

\begin{proof}
	Let $n_1,\ldots,n_k$ be the dimensions of the grid, we can associate each vertex to a $k$-vector of coordinates $(x_1,\ldots,x_k)$ with $x_i \in \{1,\ldots,n_i\}$. The central point of the grid is the point of coordinates $(\frac{n_1+1}{2},\ldots,\frac{n_k+1}{2})$. If this point corresponds to a vertex (that is, all dimensions are odd), then, the first player can select it. Now, if the second player selects the vertex of coordinates $(x_1,\ldots,x_k)$, the first player can select the opposite vertex of coordinates $(n_1+1-x_1,\ldots,n_k+1-x_k)$. This move is always possible, and thus the first player always has an answer to the second player's move, hence the game is \outcomeN. Conversely, if the central point is not a vertex (that is, any dimension is even), the first player can only select a vertex and the second player can select its opposite vertex, as described above, and hence the game is \outcomeP.
\end{proof}

Paths are another graph class where the outcome can be derived from a symmetry strategy.
Going further, we hereby give a full characterization of the Sprague-Grundy values of paths (which complements the possible computation using the algorithm for trees from Araujo \emph{et al.}~\cite{araujo2024graph}). Note that, while the value is expected, the proof itself is more involved.

\begin{theorem}
	\label{thm-paths}
	For any positive integer $n$, $\gr(P_n) = n \bmod 2$.
\end{theorem}

\begin{proof}
	Denote the vertices of $P_n$ by $u_1,\ldots,u_n$. Note that if the first two moves consisted in selecting two vertices $u_i$ and $u_j$ (with $1 \leq i < j \leq n$), the game on $P_n$ with $u_i$ and $u_j$ selected is equivalent to the game on $P_{n-(j-i)}$ with $u_i$ selected since none of the $j-i-1$ vertices between $u_i$ and $u_j$ can be selected anymore.
	
	If $n$ is even, then, the second player has a winning strategy by symmetry: whenever the first player selects $u_i$, they can always select $u_{n+1-i}$, which implies that $P_n$ is \outcomeP and thus that $\gr(P_n)=0$.
	
	Assume now that $n$ is odd. First, note that the first player has a winning strategy, by selecting the middle vertex and then applying the symmetry strategy described above. Hence, $\gr(P_n)>0$. We now prove that no option of $P_n$ has Sprague-Grundy value~1. In order to do this, we will prove that every option of $P_n$ has an option with Sprague-Grundy value~1. We first prove the following claim:
	
	\begin{claim}
		\label{clm-evenPathsPlayedInTheMiddle}
		For every nonnegative integer $m$, $\gr(P_{4m+2},\{u_{2m+1}\})=1$.
	\end{claim}
	
	\begin{proof}
		The result is proved by induction. It trivially holds if $m=0$, since only one move is possible. Assume now that the result holds for every $k < m$. Note that the first player has a winning move, by selecting $u_{2m+2}$, so $\gr(P_{4m+2},\{u_{2m+1}\})>0$. The other possible moves of the first player are to select $u_x$ for either $1 \leq x \leq 2m$ or $2m+3 \leq x \leq 4m+2$. If $x \leq 2m$, then, the resulting position is $P_{2m+1+x},u_x$; the second player can select either $u_{2m+1}$ (resulting in $P_{2x},u_x$) or $u_{2m+3}$ (resulting in $P_{2x-2},u_x \equiv P_{2x-2},u_{x-1}$), depending on whether or not $2x$ is a multiple of~4. Both options are always available, and thus at least one of them is of the form $P_{4m'+2},u_{2m'+1}$ with $m'<m$, thus having Sprague-Grundy value~1 by induction hypothesis. A similar reasoning can be applied if $x \geq 2m+3$, and thus every option of $P_{4m+2},\{u_{2m+1}\}$ has an option with Sprague-Grundy value~1, which implies that $\gr(P_{4m+2},\{u_{2m+1}\}) = 1$.
	\end{proof}
	
	Going back to $P_n$, the non-winning options consist in selecting a vertex $u_i$ with $i \neq \frac{n+1}{2}$. If $2i$ is a multiple of~4, then, the second player can select $u_{n-i+2}$, resulting in $P_{2(i-1)},\{u_i\} \equiv P_{2(i-1)},\{u_{i-1}\}$, which has Sprague-Grundy value~1 by \Cref{clm-evenPathsPlayedInTheMiddle}. Otherwise, $2i$ is not a multiple of~4, and the second player can select $u_{n-i}$, resulting in $P_{2i},u_i$, which has Sprague-Grundy value~1 by \Cref{clm-evenPathsPlayedInTheMiddle}. Hence, every option of $P_n$ has an option with Sprague-Grundy value~1, thus no option of $P_n$ can have Sprague-Grundy value~1, which proves that $\gr(P_n)=1$.
\end{proof}

\section{Sprague-Grundy values of block graphs and cacti}
\label{sec-blockAndCacti}

In their paper~\cite[Algorithm 2]{araujo2024graph}, Araujo \emph{et al.} presented a linear-time algorithm for computing the Sprague-Grundy value of a given tree for \AGAG. In order to do this, they use dynamic programming with subtrees that decrease in size at each step: for any given selected vertex $u$, each of the subtrees rooted at a neighbour of $u$ are completely independent in \AGAG, since by \Cref{cut_vertex_in_closure} selecting a vertex in a subtree will never impact any of the other subtrees. The base case happens when the subtree is either empty or a single vertex. Hence, it is possible to define a recursive algorithm computing, for each vertex of the tree, the Sprague-Grundy value of the subtrees rooted at each of its neighbours, and then apply the $\mex$ operation to the resulting set. Since each subtree is independent, and the intermediate values can be stored, the algorithm returns the result in linear-time.

In this section, we extend this idea to other tree-like graph classes: block graphs and cacti. Indeed, the building blocks of those graphs are simple graphs for which we can easily compute the Sprague-Grundy values: complete graphs and cycles. We will present polynomial-time algorithms for those classes, adapted from the trees' algorithm in~\cite{araujo2024graph}. Note that the building blocks have specific properties, and are close to the class of \emph{geodetic graphs}, that is, graphs for which there is a unique shortest path between any pair of vertices (this class includes trees, block graphs, and cacti with only odd cycles); the only non-geodetic class is cacti with even cycles, for which the shortest paths between two vertices will share common articulation points and thus are almost as constrained.

\subsection{Block graphs}

A graph is a \emph{block graph} if every biconnected component is a clique. Informally, we can see a block graph as a tree-like structure in which (inclusion-wise) maximal cliques are connected to each other by one of their vertices. Due to this property, we use the shorthand \emph{clique of $G$} to mean \emph{inclusion-wise maximal clique in $G$}.

\begin{theorem}
	\label{thm-blockGraphs}
	There is a quadratic-time algorithm computing the Sprague-Grundy value of a block graph for \AGAG.
\end{theorem}

\begin{proof}
	We will present the algorithm, then show that it computes the desired Sprague-Grundy value, before analyzing its time complexity.
	Throughout this proof, $G \setminus v$ (resp. $G \setminus S$) denotes the graph $G$ with the vertex $v$ (resp. the set of vertices $S$) removed, while $G + (w,vw)$ (where $v$ is a vertex of $G$) denotes the graph $G$ where a vertex $w$ is created and connected to only $v$ in $G$.
	
	The main algorithm \texttt{\grFBG($G$)} takes as an input a given block graph $G$ and outputs its Sprague-Grundy value $\gr(G)$. It uses a subalgorithm called \texttt{\grFBV($G$,$u$)}, which takes as an input a given block graph $G$ and a vertex $u$, and outputs $\gr(G,\{u\})$ (that is, the Sprague-Grundy value of $G$ with vertex $u$ selected). The principle of both algorithms is to decompose the input in order to apply recursion and compute values of higher-order graphs with the nim-sum and the $\mex$, while maintaining that at most one vertex of the input graph is selected in the geodetic set.
	Note that intermediate values computed by the recursive calls can be stored using an appropriate hash function, although this will not be detailed.
	For readability reasons, we write the main structure of \texttt{\grFBG($G$)} and \texttt{\grFBV($G$,$u$)} in text rather than in pseudo-code, allowing us to use slightly less formal language. Both descriptions are also backed by figures.
	
	\medskip
	\noindent\textbf{The algorithm \texttt{\grFBG($G$)}.} The algorithm works as follows:
	
	\begin{enumerate}
		\item If $G$ is a clique $K_n$, return $n \bmod 2$.
		
		\item Create a list $L$ of candidates for the first move:
		\begin{enumerate}
			\item $L$ initially contains every articulation point of~$G$;
			\item For every clique $A$ of $G$ containing at least one vertex that is not an articulation point, add one vertex of $A$ to $L$.
		\end{enumerate}
		
		\item Create a list $O$ of options and a list $T$ of integers, both initially empty.
		
		\item For every vertex $u \in L$ that is an articulation point:
		\begin{enumerate}
			\item Create the list $O_u$, initially empty. Let $H_1,\ldots,H_k$ be the components separated by $u$.
			
			\item Add $H_i,u$ for every $i \in \{1,\ldots,k\}$ to $O_u$.
			
			\item Add $O_u$ to $O$.
		\end{enumerate}
		
		\item For every vertex $u \in L$ that is not an articulation point:
		\begin{enumerate}
			\item Create the list $O_u$, initially empty. Let $A$ be the clique containing $u$, $v_1,\ldots,v_k$ be the vertices of $A$ that are articulation points, and $H_1,\ldots,H_k$ be the subgraphs separated from $A$ by those (note that $v_i$ may be an articulation point within $H_i$).
			
			\item Add $A \setminus \{v_1,\ldots,v_k\},\{u\}$, as well as $H_i+(w_i,v_iw_i),\{w_i\}$ for every $i \in \{1,\ldots,k\}$, to $O_u$.
			
			\item Add $O_u$ to $O$.
		\end{enumerate}
		
		\item For every option $o$ in $O$:
		\begin{enumerate}
			\item If $o$ is of the form $H$ where $H$ is a graph, add \texttt{\grFBG($H$)} to $T$.
			
			\item If $o$ is of the form $H,u$ where $H$ is a graph and $u$ a vertex, add \texttt{\grFBV($H$,$u$)} to $T$.
			
			\item If $o$ is of the form of $O_u$, where $O_u$ is a sum of options, then $O_u$ is a list of options $H_1,\ldots,H_k$ and $(H_{k+1},u_{k+1}),\ldots,(H_\ell,u_\ell)$ where each $H_i$ is a graph and each $u_i$ a vertex. Add \begin{center}\texttt{\grFBG($H_1$)} $\oplus$ $\ldots$ $\oplus$ \texttt{\grFBG($H_k$)} $\oplus$\\ \texttt{\grFBV($H_{k+1}$,$u_{k+1}$)} $\oplus$ $\ldots$ $\oplus$ \texttt{\grFBV($H_\ell$,$u_\ell$)}\end{center} to $T$.
		\end{enumerate}
		
		\item Return $\mex(T)$
	\end{enumerate}
	
	The algorithm \texttt{\grFBG($G$)} iterates over all possible moves for the first player, and decomposes the graph afterwards. The decomposition is depicted on \Cref{fig-blockGraphs-1}. This gives all possible options, on which a $\mex$ can then be computed. Note that it will mostly be called in the induction to manage cliques.
	The first decomposition is given by \Cref{cut_vertex_in_closure}, while the second one is given by the following claim:
	
	\begin{claim}
		\label{clm-blockGraphs-SelectingANonArtPoint}
		Let $u$ be a vertex that is not an articulation point in a clique $A$ of the block graph $G$. Let $v$ be an articulation point of $A$, separating subgraphs $H_1$ and $H_2$ (we consider that $v$ is in both $H_1$ and $H_2$) such that $u \in H_1$ and that $H_2$ is maximal (that is, $H_2$ contains $v$ and every vertex $x$ such that $v$ is in every shortest path between $u$ and $x$).
		
		We have $\gr(G,\{u\})=\gr(H_1 \setminus v,\{u\}) \oplus \gr(H_2+(w,vw),\{w\})$.
	\end{claim}
	
	\begin{proof}
		Denote by $H'_1$ the graph $H_1 \setminus v$, that is, the graph obtained by removing $v$ from $H_1$; and by $H'_2$ the graph $H_2+(w,vw)$, that is, the graph obtained by adding to $H_2$ the vertex $w$ and the edge $vw$.
		In order to prove the claim, it suffices to show that the two subgraphs are independent from each other in $G$, and that $w$ plays the role of $u$ in $H'_2$.
		
		Selecting a vertex $x$ from $H'_1$ in $G$ will never affect the component $H'_2$: let $y$ be a vertex in $H'_2$, the distance between $u$ and $y$ is smaller or equal than the distance between $x$ and $y$, and so selecting $x$ affects no vertex in $H'_2$.
		Similarly, selecting a vertex from $H'_2$ in $G$ will never affect the component $H'_1$.
		Finally,  in $H'_2,\{w\}$, $w$ plays the role of $u$: for every vertex $x \in H'_2$, every vertex on a shortest path in $G$ between $u$ and $x$ is also on a shortest path in $H'_2$ between $w$ and $x$.
		
		Hence, the two components are independent and $w$ acts exactly as $u$ in $H'_2$, so we can decompose $G,\{u\}$ into the disjoint sum of the two games $H'_1,\{u\}$ and $H'_2,\{w\}$.
	\end{proof}
	
	\begin{claim}
		\label{clm-blockGraphs-FirstMove}
		Provided that the Sprague-Grundy values obtained by the recursive calls to \texttt{\grFBG} and \texttt{\grFBV} are correct, the algorithm above outputs the Sprague-Grundy value of a given block graph $G$.
	\end{claim}
	
	\begin{proof}
		Assume that the recursive calls return the desired Sprague-Grundy values. First, if $G$ is a complete graph, we are done by \Cref{prop-activeGeodeticSheLovesMe}. Assume that $G$ contains at least two cliques.
		
		The algorithm begins at Step~2 by creating the list of all distinct moves on $G$. Note that all articulation points are distinct moves, but for all vertices in the same clique that are not articulation points, selecting one or the other does not have any effect on other vertices. Consider $u_1$ and $u_2$ two such vertices, for any other vertex $v$, the shortest paths between $u_1$ and $v$ and between $u_2$ and $v$ only differ when they reach $u_1$ and $u_2$. Hence, only one such vertex is selected for each clique.
		
		For each possible move, Steps~4 and~5 then compute the option obtained by selecting each given vertex. The option is a decomposition, depending on the nature of the vertex.
		In Step~4, we consider the option given by selecting an articulation point. The Sprague-Grundy value of this option is obtained with \Cref{cut_vertex_in_closure}.
		In Step~5, we consider the option given by selecting a vertex that is not an articulation point. The Sprague-Grundy value of this option is obtained by applying \Cref{clm-blockGraphs-SelectingANonArtPoint} to each articulation point of the clique containing the selected vertex.
		
		The algorithm will then, for each of these options, compute the Sprague-Grundy value of the option obtained by selecting the vertex (Step~6) with recursive calls to the two algorithms (assumed to be correct). Finally, a $\mex$ computation over those options will, by definition, give the correct Sprague-Grundy value of $G$ (Step~7).
	\end{proof}
	
	\medskip
	\noindent\textbf{The algorithm \texttt{\grFBV($G$,$u$)}.} The algorithm works as follows:
	
	\begin{enumerate}
		\item If $G$ is a clique $K_n$, return $(n-1) \bmod 2$.
		
		\item Create a list $L$ of candidates for the second move:
		\begin{enumerate}
			\item $L$  initially contains every articulation point of~$G$;
			\item For every clique $A$ of $G$ containing at least one vertex that is neither $u$ nor an articulation point, add one vertex of $A$ to $L$.
		\end{enumerate}
		
		\item Create a list $O$ of options and a list $T$ of integers, both initially empty.
		
		\item If there is a vertex $v \in L$ that is in the same clique as $u$ and is not an articulation point, add $G \setminus v, \{u\}$ to $O$.
		
		\item For every vertex $v \in L$ that is in the same clique as $u$ and is an articulation point:
		\begin{enumerate}
			\item Create the list $O_v$, initially empty. Let $H$ be the component containing $u$, and $H_1,\ldots,H_k$ be the other components separated from $H$ by $v$.
			
			\item Add $H \setminus v, \{u\}$, as well as $H_i, \{v\}$ for every $i \in \{1,\ldots,k\}$, to $O_v$.
			
			\item Add $O_v$ to $O$.
		\end{enumerate}
		
		\item For every vertex $v \in L$ that is not in the same clique as $u$:
		\begin{enumerate}
			\item Create the list $O_v$, initially empty. Let $p$ be the shortest path between $u$ and $v$, and let $w_1,\ldots,w_k$ be the articulation points of $G$ on $p$ when starting from $u$ and ending at $v$ (we have $w_1 \neq u$ and $w_k \neq v$).
			
			\item Let $H$ and $H'$ be the maximal (inclusion-wise) connected subgraphs containing, respectively, $u$ and not $w_1$ for $H$, and $v$ and not $w_k$ for $H'$. Add $H,u$ and $H',v$ to $O_v$ (the latter only if $v$ is not an articulation point).
			
			\item For every $i \in \{1,\ldots,k\}$, if $w_i$ is an articulation point separating components $H$, $H'$, $H_1,\ldots,H_\ell$ such that $u \in H$ and $v \in H'$, then, add $H_j,\{w_i\}$ to $O_v$ for every $j \in \{1,\ldots,\ell\}$.
			
			\item For every $i \in \{1,\ldots,k-1\}$, let $A$ be the clique containing $w_i$ and $w_{i+1}$, $x_1,\ldots,x_\ell$ be the other articulation points of $A$, and $H_1,\ldots,H_\ell$ be the subgraphs separated from $A$ by those (note that $x_j$ may be an articulation point within $H_j$). Add $A \setminus \{w_i,w_{i+1},x_1,\ldots,x_\ell\}$, as well as $H_j + (y_j,x_jy_j), \{y_j\}$ for every $j \in \{1,\ldots,\ell\}$, to $O_v$.
			
			\item If $v$ is an articulation point separating components $H$, $H_1,\ldots,H_\ell$ such that $u \in H$, then, add $H_j,\{v\}$ to $O_v$ for every $j \in \{1,\ldots,\ell\}$.
			
			\item Add $O_v$ to $O$.
		\end{enumerate}
		
		\item For every option $o$ in $O$:
		\begin{enumerate}
			\item If $o$ is of the form $H$ where $H$ is a graph, add \texttt{\grFBG($H$)} to $T$.
			
			\item If $o$ is of the form $H,v$ where $H$ is a graph and $v$ a vertex, add \texttt{\grFBV($H$,$v$)} to $T$.
			
			\item If $o$ is of the form of $O_v$, where $O_v$ is a sum of options, then $O_v$ is a list of options $H_1,\ldots,H_k$ and $(H_{k+1},v_{k+1}),\ldots,(H_\ell,v_\ell)$ where each $H_i$ is a graph and each $v_i$ a vertex. Add \begin{center}\texttt{\grFBG($H_1$)} $\oplus$ $\ldots$ $\oplus$ \texttt{\grFBG($H_k$)} $\oplus$\\ \texttt{\grFBV($H_{k+1}$,$v_{k+1}$)} $\oplus$ $\ldots$ $\oplus$ \texttt{\grFBV($H_\ell$,$v_\ell$)}\end{center} to $T$.
		\end{enumerate}
		
		\item Return $\mex(T)$
	\end{enumerate}
	
	Assume that $u$ is selected on the block graph $G$.
	The algorithm \texttt{\grFBV($G$,$u$)} iterates over all possible moves for the second player, and decomposes the graph afterwards. The decomposition is depicted on \Cref{fig-blockGraphs-2}. This gives all possible options, on which a $\mex$ can then be computed.
	We use the following claim to manage moves in the same clique that are not articulation points:
	
	\begin{claim}
		\label{clm-blockGraphs-InTheSameClique}
		Let $v$ be a vertex in the same clique as $u$, that is not an articulation point. We have $\gr(G,\{u,v\})=\gr(G \setminus v,\{u\})$.
	\end{claim}
	
	\begin{proof}
		The vertices $u$ and $v$ are both simplicial vertices, and thus will always both be in the final geodetic set. So it suffices to show that selecting $v$ and removing it from $G$ are equivalent. Consider any vertex $w$ in a different clique: the shortest path between $w$ and $u$ and the shortest path between $w$ and $v$ are identical except for $u$ and $v$. Hence, selecting $v$ when $u$ has been selected is equivalent to removing $v$ from $G$.
	\end{proof}
	
	The following claim can be applied repeatedly to decompose subgraphs attached to articulation points on the shortest path between $u$ a second move in another clique (see $C$ in \Cref{fig-blockGraphs-1-2}):
	
	\begin{claim}
		\label{clm-blockGraphs-InAnotherClique-CliquesSeparatedByTheShortestPath}
		Let $v$ be a vertex in a different clique from $u$, and let $w$ be an articulation point (distinct from $u$ and $v$) on the shortest path between $u$ and $v$. 
		Let $H_2$ be another clique joined at $w$, and let $x$ be a vertex from $H_2$ distinct from $w$.
		Denote by $H'_2$ the subgraph of $G$ containing $w$ and all vertices $y$ such that $w$ is not on the shortest path between $x$ and $y$. Denote by $H'_1$ the subgraph of $G$ containing $w$ and all vertices that are not in $H'_2$.
		
		We have $\gr(G,\{u,v\}) = \gr(H'_1,\{u,v\}) \oplus \gr(H'_2,\{w\})$.
	\end{claim}
	
	\begin{proof}
		It suffices to show that the two subgraphs are independent from each other in $G$, that all forbidden moves in one are also forbidden in the other, and that $w$ plays the role of $u$ and $v$ in $H'_2$.
		Selecting a vertex $x$ from $H'_1$ in $G$ will never affect the component $H'_2$: since $w$ is on the shortest path between $x$ and any vertex $y \in H'_2$, from the point of view of $y$, the vertices $x$ and $u$ (or $v$) are undistinguishable and the shortest paths between $x$ and $y$ and between $u$ (or $v$) and $y$ will cover the exact same vertices in $H'_2$.
		Similarly, selecting a vertex from $H'_2$ in $G$ will never affect the component $H'_1$.
		Furthermore, $w$ is on the shortest path between $u$ and $v$, and thus cannot be selected in $G$, hence it is correct to consider it as already-selected in $H'_2$.
		Finally,  in $H'_2,\{w\}$, $w$ plays the role of both $u$ and $v$: for every vertex $x \in H'_2$, every vertex on a shortest path in $G$ between $u$ and $x$ (or between $v$ and $x$) is also on a shortest path in $H'_2$ between $w$ and $x$.
		
		Hence, the two components are independent and $w$ acts exactly as $u$ and $v$ in $H'_2$, so we can decompose $G,\{u,v\}$ into the disjoint sum of the two games $H'_1,\{u,v\}$ and $H'_2,\{w\}$.
	\end{proof}
	
	The following two claims can also be applied (repeatedly for the second one) to simplify consecutive articulation points on the shortest path between $u$ and a second move in another clique (see, in \Cref{fig-blockGraphs-1-2}, $B_2,B_3$ for \Cref{clm-blockGraphs-InAnotherClique-ConsecutiveArticulationPointsOnTheShortestPathOutsideClique} and $B_1$ for \Cref{clm-blockGraphs-InAnotherClique-ConsecutiveArticulationPointsOnTheShortestPathInsideClique}):
	
	\begin{claim}
		\label{clm-blockGraphs-InAnotherClique-ConsecutiveArticulationPointsOnTheShortestPathOutsideClique}
		Let $v$ be a vertex in a different clique from $u$, and let $w_1$ and $w_2$ be two articulation points (distinct from $u$ and $v$) on the shortest path between $u$ and $v$. Denote by $A$ the clique containing $w_1$ and $w_2$, and let $x_1,\ldots,x_k$ be the other articulation points of $A$, separating it from cliques $B_1,\ldots,B_k$. Let $H_2^1,\ldots,H_2^k$ be the maximal connected subgraphs of $G$ each containing $x_i$ and $B_i$ but no vertex from $A$ other than $x_i$ for $i \in \{1,\ldots,k\}$. Denote by $H_1$ the subgraph of $G$ containing all vertices that are not in $H_2^1,\ldots,H_2^k$.
		
		We have $\gr(G,\{u,v\})=\gr(H_1,\{u,v\}) \oplus \gr(H_2^1 + (y_1,x_1y_1), \{y_1\}) \oplus \ldots \oplus \gr(H_2^k + (y_k,x_ky_k), \{y_k\})$.
	\end{claim}
	
	\begin{proof}
		Denote, for $i \in \{1,\ldots,k\}$, by $H'^i_2$ the graph obtained by adding vertex $y_i$ and edge $x_iy_i$ to $H_2^i$.
		It suffices to show that the $k+1$ subgraphs are independent from each other in $G$, and that $y_i$ plays the role of $u$ and $v$ in $H'^i_2$.
		Clearly, selecting a vertex from $H_1$ in $G$ will never affect any component $H_2^i$, since they are completely separated by $x^i$. Any selected vertex $y \in H_1$ is thus undistinguishable from $u$ or $v$ from this point of view: every shortest path from $y$ to a vertex in $H_2^i$ will have to pass through $x^i$, exactly as all shortest paths from $u$ or $v$ to the same vertex of $H_2^i$.
		
		Conversely, selecting a vertex $z \in H_2^i$ does not affect the component $H_1$: for every vertex $t \in A$, $z$ itself will be the only vertex of $H_1$ on the shortest path between $z$ and $t$, and for every vertex $t \in H_1 \setminus A$, all shortest paths go through $x_i$ and then (without loss of generality) $w_1$, and thus it is undistinguishable from $v$ (the shortest path from $v$ to $t$ will go through $w_1$).
		Finally, for $i \neq j$, if we select a vertex in $H_2^i$, $H_2^j$ is not affected as they are separated in $G$ by two articulation points $x_i$ and $x_j$, neither of which is on the shortest path between $u$ and $v$.
		Hence, all components are independent from each other.
		
		Let $i \in \{1,\ldots,k\}$ and consider the position $H'^i_2,\{y_i\}$. Clearly, for any vertex $z \in H_2^i$, every vertex of $H_2^i$ on the shortest path between $u$ (or $v$) and $z$ is also on the shortest path between $y_i$ and $z$.
		Note in particular that two moves in components $H'^i_2$ and $H'^j_2$ for $i \neq j$ will, in $G$, cover the two articulation points $x_i$ and $x_j$ and no vertex outside of $H_2^i \cup H_2^j$, and so these moves are perfectly replicated on $H'^i_2,\{y_i\}$ and $H'^j_2,\{y_j\}$.
		
		Hence, the $k+1$ components are independent and each $y_i$ acts exactly as $u$ and $v$ in $H'^i_2$, so we can decompose $G,\{u,v\}$ into the disjoint sum of the $k+1$ games $H_1,\{u,v\}$ and $H'^i_2,\{y_i\}$ for $i \in \{1,\ldots,k\}$.
	\end{proof}
	
	\begin{claim}
		\label{clm-blockGraphs-InAnotherClique-ConsecutiveArticulationPointsOnTheShortestPathInsideClique}
		Let $v$ be a vertex in a different clique from $u$, and let $w_1$ and $w_2$ be two articulation points (distinct from $u$ and $v$) on the shortest path between $u$ and $v$. Denote by $A$ the clique containing $w_1$ and $w_2$, and assume that $A$ contains no other articulation point and that the only cliques joined at $w_1$ and $w_2$ contain two vertices on the shortest path between $u$ and $v$ ($w_1$ or $w_2$ and another one).
		Denote by $H$ the graph obtained by removing from $G$ all vertices of $A$ except for $w_1$ and $w_2$ and then merging $w_1$ and $w_2$ (thus reattaching the two components separated by $A$).
		
		We have $\gr(G,\{u,v\}) = \gr(H,\{u,v\}) \oplus \gr(A \setminus \{w_1,w_2\})$.
	\end{claim}
	
	\begin{proof}
		All the vertices of $A$ will be selected at the end of the game (since they are simplicial vertices), except for $w_1$ and $w_2$ which are on the shortest path between $u$ and $v$. Furthermore, moves on $A$ have no effect on the rest of the graph, since any shortest path between $x \in A$ and $y \not\in A$ will go through (without loss of generality) $w_1$ and thus from the perspective of $y$ the vertex $x$ is undistinguishable from $v$.
		Thus, $A \setminus \{w_1,w_2\}$ can be removed from $G$ and considered as an independent game. Doing this, the clique becomes a simple edge, and the two vertices can be fused together without changing anything to the shortest paths in $G$ (since both vertices are already in the geodetic closure).
		
		Hence, we can decompose and simplify $G,\{u,v\}$ into the disjoint sum of the two games $H,\{u,v\}$ and $A \setminus \{w_1,w_2\}$.
	\end{proof}
	
	When \Cref{clm-blockGraphs-InAnotherClique-CliquesSeparatedByTheShortestPath}, \Cref{clm-blockGraphs-InAnotherClique-ConsecutiveArticulationPointsOnTheShortestPathOutsideClique} and \Cref{clm-blockGraphs-InAnotherClique-ConsecutiveArticulationPointsOnTheShortestPathInsideClique} have been repeatedly applied, the following claim gives the final decomposition (see $A$ and $D$ in \Cref{fig-blockGraphs-1-2}):
	
	\begin{claim}
		\label{clm-blockGraphs-InAnotherClique-NeighbouringCliques}
		Let $v$ be a vertex in a different clique from $u$, and assume that those two cliques are joined at one articulation point $w$. Assume furthermore than no other clique is joined at $w$. Denote by $H_u$ and $H_v$ be the subgraphs separated by $w$ such that $H_u$ (resp. $H_v$) contains $u$ (resp. $v$).  We have $\gr(G,\{u,v\})=\gr(H_u \setminus w,\{u\}) \oplus \gr(H_v \setminus w,\{v\})$.
	\end{claim}
	
	\begin{proof}
		Denote $H_u \setminus w$ and $H_v \setminus w$ by $H'_u$ and $H'_v$, respectively. Note that we cannot play on $w$ since it is on the shortest path between $u$ and $v$.
		As in previous claims, we prove that the two parts are independent from each other. Let $x \in H'_u$ and $y \in H'_v$. The shortest path between $x$ and $y$ goes through $w$, and thus is in $H'_u$ identical to the shortest path between $x$ and $v$, and in $H'_v$ identical to the shortest path between $u$ and $y$. Thus, selecting a vertex from $H'_u$ has no effect on $H'_v$, and conversely.
		Furthermore, when selecting $x \in H'_u$, the shortest path between $x$ and $u$ and the shortest path between $x$ and $w$ are identical except for $u$ and $w$. The same reasoning applies for $H'_v$. Thus, no other vertex needs to be selected to simulate the other vertex for this part.
		
		Hence, we can decompose $G,\{u,v\}$ into the disjoint sum of the two games $H'_u,\{u\}$ and $H'_v,\{v\}$.
	\end{proof}
	
	We are now ready to prove that the algorithm \texttt{\grFBV($G$,$u$)} outputs the correct value:
	
	\begin{claim}
		\label{clm-blockGraphs-SecondMove}
		Provided that the Sprague-Grundy values obtained by the recursive calls to \texttt{\grFBG} and \texttt{\grFBV} are correct, the algorithm above outputs the Sprague-Grundy value of the position $G,\{u\}$.
	\end{claim}
	
	\begin{proof}
		Assume that the recursive calls return the desired Sprague-Grundy values. First, if $G$ is a complete graph, we are done by \Cref{prop-activeGeodeticSheLovesMe}. Assume that $G$ contains at least two cliques.
		
		Again, the algorithm begins at Step~2 by creating the list of all distinct moves on $G,\{u\}$, which are articulation points and one vertex per clique. The algorithm again is a decomposition in Steps~4 to~6. Step~7 is the computation of the Sprague-Grundy value of each option with the nim-sum and the recursive calls to the two algorithms (assumed to be correct), and Step~8 is the $\mex$ computation, giving us the Sprague-Grundy value of $G,\{u\}$.
		
		We now detail the three possible types of options for the second player. First (Step~4), the second player can select a vertex that is not an articulation point in the same clique as $u$, in which case \Cref{clm-blockGraphs-InTheSameClique} gives us the option (simply remove the vertex from the graph).
		In Step~5, we consider the options when selecting an articulation point in the same clique as $u$. In this case, \Cref{cut_vertex_in_closure} gives us the decomposition (see \Cref{fig-blockGraphs-1-1}).
		
		Step~6 is the most complex decomposition, which occurs when selecting a vertex $v$ in a different clique (see \Cref{fig-blockGraphs-1-2}). Several steps occur for decomposing the graph: first, consider the shortest path between $u$ and $v$, and denote by $w_1,\ldots,w_k$ the articulation points along it (all distinct from $u$ and $v$). We will now create a large disjoint sum of positions.
		We begin by decomposing subgraphs attached to each $w_i$ (\Cref{clm-blockGraphs-InAnotherClique-CliquesSeparatedByTheShortestPath}). Now, all $w_i$'s separate only cliques on the shortest path between $u$ and $v$.
		For each pair of consecutive articulation points $w_i,w_{i+1}$, consider the clique $A$ to which they both belong. We decompose all subgraphs attached to $A$ by other articulation points (\Cref{clm-blockGraphs-InAnotherClique-ConsecutiveArticulationPointsOnTheShortestPathOutsideClique}). We can now separate $A$ from the graph and fuse together $w_i$ and $w_{i+1}$ (\Cref{clm-blockGraphs-InAnotherClique-ConsecutiveArticulationPointsOnTheShortestPathInsideClique}), and iterate until $u$ and $v$ are separated by only one articulation point. The final decomposition is given by \Cref{clm-blockGraphs-InAnotherClique-NeighbouringCliques}.
	\end{proof}
	
	Consider now the decomposition process that happens when using the algorithms \texttt{\grFBG($G$)} and \texttt{\grFBV($G$,$u$)}, starting from a given block graph $G$. Note that every position is the decomposition has at most one selected vertex, allowing us to repeatedly apply the algorithms. This gives a game directed graph (that is, a directed graph such that every vertex is a position, and every arc represents an option), for which every terminal position (that is, a position that cannot be decomposed anymore) is a complete graph. Finally, note that the same vertex cannot be selected twice in the resulting decomposition.
	\Cref{clm-blockGraphs-FirstMove} and \Cref{clm-blockGraphs-SecondMove} prove that the computation of the Sprague-Grundy of a given position is correct.
	
	The algorithm works in quadratic time: there is a quadratic number of possible decomposed subgraphs (each clique can have a linear number of articulation points, and can be decomposed a linear number of times over each of them), which can be generated beforehand. The game directed graph can then be constructed in linear time at each recursive step, and resolving the graph (and thus the recursive calls) can finally be done in quadratic time. Note in particular that some blocks which appear in decompositions of larger graphs also appear in decompositions of smaller graphs; as an example, see on \Cref{fig-blockGraphs-1}, the block with cliques $B_3^a$, $B_3^b$ and $B_3^c$ appears whenever the selected vertex is $v_1$, $v_2$ or a vertex in the cliques $A$, $B_1^a$, $B_1^b$ or $B_2$, so its Sprague-Grundy value will only be computed once.
\end{proof}

\begin{figure}[h]
	\centering
	\makebox[\textwidth][c]{
		\scalebox{1}{
			\begin{tikzpicture}
				\node (graph) at (0,0) {
					\begin{tikzpicture}
						\draw (0,0) ellipse (1 and 0.5);
						\draw (0,0) node {$A$};
						\node[noeud,fill=black] (u) at (0,0.5) {};
						\draw (0,0.75) node {$u$};
						
						\draw[rotate around={-35:(-1.75,0)}] (-1.9,0.375) ellipse (0.75 and 0.375);
						\draw[rotate around={35:(-1.5,0)}] (-1.6,-0.375) ellipse (0.5 and 0.25);
						\draw (-1.75,0.375) node {$B_1^a$};
						\draw (-1.5,-0.375) node {$B_1^b$};
						\node[noeud] (v1) at (-1,0) {};
						\draw (-0.7,0) node {$v_1$};
						
						\draw (0,-1.25) ellipse (0.5 and 0.75);
						\draw (0,-1.25) node {$B_2$};
						\node[noeud] (v2) at (0,-0.5) {};
						\draw (0,-0.75) node {$v_2$};
						
						\draw (2,0) ellipse (1 and 0.5);
						\draw[rotate around={25:(3.75,0.25)}] (3.75,0.5) ellipse (0.75 and 0.375);
						\draw[rotate around={-25:(3.75,-0.25)}] (3.5,-0.5) ellipse (0.5 and 0.25);
						\draw (2,0) node {$B_3^a$};
						\draw (3.75,0.5) node {$B_3^b$};
						\draw (3.5,-0.5) node {$B_3^c$};
						\node[noeud] (v3) at (1,0) {};
						\draw (0.7,0) node {$v_3$};
					\end{tikzpicture}
				};
				
				\draw (3.5,0) node {$\equiv$};
				
				\node (graphs) at (9,0) {
					\begin{tikzpicture}
						\draw (0,0) ellipse (1 and 0.5);
						\draw (0,0) node[scale=0.85] {$A \setminus \{v_1,v_2,v_3\}$};
						\node[noeud,fill=black] (u) at (0,0.5) {};
						\draw (0,0.75) node {$u$};
						
						\draw (-2.25,0.75) ellipse (0.375 and 0.75);
						\draw (-2.25,-0.5) ellipse (0.25 and 0.5);
						\draw (-2.25,0.75) node {$B_1^a$};
						\draw (-2.25,-0.5) node {$B_1^b$};
						\node[noeud] (v1) at (-2.25,0) {};
						\draw (-2.625,0) node {$v_1$};
						\node[noeud,fill=black] (v1p) at (-1.75,0) {};
						\draw (v1)to(v1p);
						\draw (-1.3,0) node {+};
						
						\draw (6.5,0) ellipse (0.5 and 0.75);
						\draw (6.5,0) node {$B_2$};
						\node[noeud] (v2) at (6.5,0.75) {};
						\draw (6.5,0.5) node {$v_2$};
						\node[noeud,fill=black] (v2p) at (6.5,1.25) {};
						\draw (v2)to(v2p);
						\draw (5.75,0) node {+};
						
						\draw (3.25,0) ellipse (1 and 0.5);
						\draw[rotate around={25:(5,0.25)}] (5,0.5) ellipse (0.75 and 0.375);
						\draw[rotate around={-25:(5,-0.25)}] (4.75,-0.5) ellipse (0.5 and 0.25);
						\draw (3.25,0) node {$B_3^a$};
						\draw (5,0.5) node {$B_3^b$};
						\draw (4.75,-0.5) node {$B_3^c$};
						\node[noeud] (v3) at (2.25,0) {};
						\draw (2.55,0) node {$v_3$};
						\node[noeud,fill=black] (v3p) at (1.75,0) {};
						\draw (v3)to(v3p);
						\draw (1.3,0) node {+};
					\end{tikzpicture}
				};
			\end{tikzpicture}
		}
	}
	\caption{Selecting a vertex that is not an articulation point in a block graph allows us to decompose the graph as a sum of subgraphs (each ellipsis is a clique).}
	\label{fig-blockGraphs-1}
\end{figure}

\begin{figure}[h]
	\centering
	\begin{subfigure}{0.4\linewidth}
		\centering
		\scalebox{1}{
			\begin{tikzpicture}
				\node (graph) at (0,0) {
					\begin{tikzpicture}
						\draw (0,0) ellipse (1 and 0.5);
						\draw (0,0) node {$A$};
						\node[noeud,fill=black] (u) at (0,0.5) {};
						\draw (u) node[above,yshift=0.5mm] {$u$};
						\node[noeud,fill=black] (v) at (1,0) {}; 
						\draw (v) node[left,xshift=-0.5mm] {$v$};
						\draw[rotate around={35:(1.75,0)}] (1.9,0.375) ellipse (0.75 and 0.375);
						\draw[rotate around={-35:(1.5,0)}] (1.6,-0.375) ellipse (0.5 and 0.25);
						\draw (1.625,0.375) node {$B_1$};
						\draw (1.375,-0.375) node {$B_2$};
						\draw (-1.5,0) ellipse (0.5 and 0.25);
						\draw (-1.5,0) node {$C$};
					\end{tikzpicture}
				};
				\node (decomposition1) at (0,-2.5) {
					\begin{tikzpicture}
						\draw (0,0) ellipse (1 and 0.5);
						\draw (0,0) node {$A$};
						\node[noeud,fill=black] (u) at (0,0.5) {};
						\draw (u) node[above,yshift=0.5mm] {$u$};
						\node[noeud,cross out] (exv) at (1,0) {}; 
						\draw (exv) node[right,xshift=0.5mm] {$v$};
						\draw (-1.5,0) ellipse (0.5 and 0.25);
						\draw (-1.5,0) node {$C$};
						
					\end{tikzpicture}
				};
				\node (decomposition2) at (-1,-4) {
					\begin{tikzpicture}
						\draw (0,0) ellipse (0.75 and 0.375);
						\draw (0,0) node {$B_1$};
						\node[noeud,fill=black] (v) at (-0.75,0) {}; 
						\draw (v) node[left,xshift=-0.5mm] {$v$};
					\end{tikzpicture}
				};
				\node (decomposition3) at (1.375,-4) {
					\begin{tikzpicture}
						\draw (0,0) ellipse (0.5 and 0.25);
						\draw (0,0) node {$B_2$};
						\node[noeud,fill=black] (v) at (-0.5,0) {}; 
						\draw (v) node[left,xshift=-0.5mm] {$v$};
					\end{tikzpicture}
				};
				\draw (0,-1.25) node[scale=1.] {$\equiv$};
				\draw (-2.25,-4) node[scale=1.] {+};
				\draw (0.375,-4) node[scale=1.] {+};
			\end{tikzpicture}
		}
		\caption{Selecting, as a second vertex, an articulation point in the same clique as the first selected vertex, gives a simple decomposition.}
		\label{fig-blockGraphs-1-1}
	\end{subfigure}
	\hfil
	\begin{subfigure}{0.55\linewidth}
		\centering
		\scalebox{1}{
			\begin{tikzpicture}
				\node (graph) at (0,0) {
					\begin{tikzpicture}
						\draw (0,0) ellipse (1 and 0.5);
						\draw (0,0) node {$A$};
						\node[noeud,fill=black] (u) at (0,0.5) {};
						\draw (u) node[above,yshift=0.5mm] {$u$};
						\draw (2,0) ellipse (1 and 0.5);
						\draw (2,-0.25) node {$B_1$};
						\node[noeud] (w1) at (1,0) {};
						\draw (w1) node[below,yshift=-2.5mm] {$w_1$};
						\draw[rotate around={35:(1.7,1)}] (1.7,1) ellipse (0.375 and 0.5);
						\draw (1.7,1) node {$B_2$};
						\draw[rotate around={-35:(2.3,1)}] (2.3,1) ellipse (0.375 and 0.5);
						\draw (2.3,1) node {$B_3$};
						\node[noeud] (x1) at (2,0.5) {};
						\draw (x1) node[below,yshift=-0.5mm] {$x_1$};
						\draw (3,0.75) ellipse (0.5 and 0.75);
						\draw (3,0.75) node {$C$};
						\node[noeud] (w2) at (3,0) {};
						\draw (w2) node[below,yshift=-2.5mm] {$w_2$};
						\draw (4,0) ellipse (1 and 0.5);
						\draw (4,0) node {$D$};
						\node[noeud,fill=black] (v) at (4,0.5) {};
						\draw (v) node[above,yshift=0.5mm] {$v$};
					\end{tikzpicture}
				};
				\node (decomposition1) at (-2.5,-2.5) {
					\begin{tikzpicture}
						\draw (0,0) ellipse (1 and 0.5);
						\draw (0,0) node {$A$};
						\node[noeud,fill=black] (u) at (0,0.5) {};
						\draw (u) node[above,yshift=0.5mm] {$u$};
						\node[noeud,cross out] (w1) at (1,0) {};
						\draw (w1) node[below,yshift=-2.5mm] {$w_1$};
					\end{tikzpicture}
				};
				\node (decomposition2) at (1.75,-4.5) {
					\begin{tikzpicture}
						\draw (0,0) ellipse (1 and 0.5);
						\draw (0,0) node {$D$};
						\node[noeud,fill=black] (v) at (0,0.5) {};
						\draw (v) node[above,yshift=0.5mm] {$v$};
						\node[noeud,cross out] (w2) at (-1,0) {};
						\draw (w2) node[left,xshift=-0.5mm] {$w_2$};
					\end{tikzpicture}
				};
				\node (decomposition3) at (0.75,-2.5) {
					\begin{tikzpicture}
						\draw (0,0) ellipse (1 and 0.5);
						\draw (0,0) node {$B_1$};
						\node[noeud,cross out] (w1) at (-1,0) {};
						\node[noeud,cross out] (w2) at (1,0) {};
						\draw (w1) node[below,yshift=-2.5mm] {$w_1$};
						\draw (w2) node[below,yshift=-2.5mm] {$w_2$};
						\node[noeud,cross out] (x1) at (0,0.5) {};
						\draw (x1) node[above,yshift=0.5mm] {$x_1$};
					\end{tikzpicture}
				};
				\node (decomposition4) at (-1.75,-4.5) {
					\begin{tikzpicture}
						\draw (0,0) ellipse (1 and 0.5);
						\draw (0,0) node {$C$};
						\node[noeud,fill=black] (w2) at (-1,0) {};
						\draw (w2) node[left,xshift=-0.5mm] {$w_2$};
					\end{tikzpicture}
				};
				\node (decomposition5) at (3.5,-2.5) {
					\begin{tikzpicture}
						\draw (0,0) ellipse (0.375 and 0.5);
						\draw (0,0) node {$B_2$};
						\draw (0,-1) ellipse (0.375 and 0.5);
						\draw (0,-1) node {$B_3$};
						\node[noeud] (x1) at (0,-0.5) {};
						\draw (x1) node[right,xshift=1.5mm] {$x_1$};
						\node[noeud,fill=black] (y1) at (-0.5,-0.5) {};
						%						\draw (y1) node[left,xshift=-0.5mm] {$y_1$};
						\draw (x1)to(y1);
					\end{tikzpicture}
				};
				\draw (0,-1.5) node[scale=1.] {$\equiv$};
				\draw (-1,-2.5) node[scale=1.] {+};
				\draw (2.375,-2.5) node[scale=1.] {+};
				\draw (-3.25,-4.5) node[scale=1.] {+};
				\draw (0,-4.5) node[scale=1.] {+};
			\end{tikzpicture}
		}
		\caption{Selecting a second vertex in a different clique allows for a more complex decomposition (here, $B_1$, $C$ and $D$ are joined at the same articulation point $w_2$; $B_1$, $B_2$ and $B_3$ are joined at articulation point $x_1$).}
		\label{fig-blockGraphs-1-2}
	\end{subfigure}
	\caption{Selecting a second vertex in a block graph allows us to decompose it further (each ellipsis is a clique, and each crossed out vertex has been removed from the subgraph).}
	\label{fig-blockGraphs-2}
\end{figure}

\subsection{Cacti}

A \emph{cactus} (plural: cacti) is a graph where every edge is in at most one simple cycle. Alternately, cacti have a tree-like structure in which cycles and trees are connected to each other by one of their vertices. All cycles that we consider in this section are simple.

\begin{theorem}
	\label{thm-cactus}
	There is a polynomial-time algorithm computing the Sprague-Grundy value of a cactus for \AGAG.
\end{theorem}

\begin{proof}
	As with block graphs, we will present the algorithm, which uses two subroutines this time. Again, $G \setminus u$ (resp. $G \setminus S$) denotes the graph $G$ with the vertex $u$ (resp. the set of vertices $S$) removed, while $G + (v,uv)$ (where $u$ is a vertex of $G$) denotes the graph $G$ where a vertex $v$ is created and connected to only $u$ in $G$.
	Furthermore, given two vertices $u$ and $v$ in the same even cycle on $2k$ vertices, we say that they are \emph{diametrically opposed} if they are at distance $k$ from each other. Note that selecting two diametrically opposed vertices in a cycle prevents from playing on this cycle in later moves.
	
	The main algorithm \texttt{\grFCG($G$)} takes as input a cactus $G$ and outputs its Sprague-Grundy value $\gr(G)$. The two subalgorithms are \texttt{\grFCVA($G$,$u$)} and \texttt{\grFCVB($G$,$u$,$v$)}, which take as input a block graph $G$ and either one vertex $u$ (which is either a leaf or a vertex in a cycle that is not an articulation point) or two vertices $u$ and $v$ (which are two vertices in the same cycle, are not articulation points, and are not diametrically opposed), and output $\gr(G,\{u\})$ and $\gr(G,\{u,v\})$, respectively.
	Again, the principle will be to decompose over every possible move, compute the Sprague-Grundy value of each option (possibly using a nim-sum), and then apply the $\mex$ over all options.
	Note that, unlike in block graphs where every non-articulation point in a clique was equivalent, this is not the case for vertices in the same cycle for cacti. This increases the number of options.
	
	\medskip
	\noindent\textbf{The algorithm \texttt{\grFCG($G$)}.} It iterates over all possible options (that is, vertices), and if a vertex is an articulation point, it decomposes $G$ into a sum of games (see \Cref{fig-cactus-1}).
	
	\begin{enumerate}
		\item If $G$ is a tree, use~\cite[Algorithm 2]{araujo2024graph} and return the value. If $G$ is a cycle on $n$ vertices, return $n \bmod 2$.
		
		\item Create a list $O$ of options and a list $T$ of integers, both initially empty.
		
		\item For every vertex $u$ of $G$:
		\begin{enumerate}
			\item If $u$ is an articulation point, denote by $H_1,\ldots,H_k$ the components separated by $u$. Create the list $O_u$, initially empty. Add $H_i,u$ for every $i \in \{1,\ldots,k\}$ to $O_u$. Add $O_u$ to $O$.
			
			\item Else, add $G,u$ to $O$.
		\end{enumerate}
		
		\item For every option $o$ in $O$:
		\begin{enumerate}		
			\item If $o$ is of the form $H,u$ where $H$ is a graph and $u$ a vertex, add \texttt{\grFCVA($H$,$u$)} to $T$.
			
			\item If $o$ is of the form of $O_u$, where $O_u$ is a list of options of the form $H_i,u_i$ for $i \in \{1,\ldots,k\}$ where each $H_i$ is a graph and each $u_i$ a vertex. Add \begin{center}\texttt{\grFCVA($H_1,u_1$)} $\oplus$ $\ldots$ $\oplus$ \texttt{\grFCVA($H_k,u_k$)}\end{center} to $T$.
		\end{enumerate}
		
		\item Return $\mex(T)$
	\end{enumerate}
	
	\medskip
	\noindent\textbf{The algorithm \texttt{Decompose($G$,$u$,$v$,$w$)}.} This subroutine takes as input a cactus $G$, two vertices $u$ and $v$, and a vertex $w$ that is an articulation point on a shortest path between $u$ and $v$ (with $w \not\in \{u,v\}$). It outputs the list of decomposed subgraphs separated by $w$, except for those containing $u$ and $v$. Those are disjoint games resulting from selecting $u$ and $v$. See \Cref{fig-cactus-2} for a depiction.
	
	\begin{enumerate}
		\item Create a list $O$, initially empty.
		
		\item Denote by $H_u$ and $H_v$ the components separated by $w$ containing $u$ and $v$, respectively (note that we may have $H_u=H_v$). Denote by $H_1,\ldots,H_k$ the other components separated by $w$.
		
		\item Add $H_i,\{w\}$ for every $i \in \{1,\ldots,k\}$ to $O$.
		
		\item Return $O$.
	\end{enumerate}
	
	\medskip
	\noindent\textbf{The algorithm \texttt{\grFCVA($G$,$u$)}.} Recall that $u$ is either a leaf or a vertex in a cycle that is not an articulation point (see \Cref{fig-cactus-3}).
	
	\begin{enumerate}
		\item If $G$ is a tree, use~\cite[Algorithm 2]{araujo2024graph} and return the value. If $G$ is an odd cycle, return 0. If $G$ is an even cycle on $2k$ vertices, return $k$.
		
		\item Create a list $O$ of options and a list $T$ of integers, both initially empty.
		
		\item For every vertex $v$ of $G$ distinct from $u$:
		\begin{enumerate}
			\item If $u$ is a leaf and $v$ is its neighbour, then, let $H$ be a copy of $G \setminus u$.
			\begin{enumerate}
				\item If $v$ is an articulation point in $H$, denote by $H_1,\ldots,H_k$ the components separated by $v$. Add $H_i,v$ for every $i \in \{1,\ldots,k\}$ to $O$.
				\item Otherwise, add $H,\{v\}$ to $O$.
			\end{enumerate}
			
			\item If $u$ and $v$ are in the same cycle $C$:
			\begin{enumerate}
				\item Create a list $O_v$, initially empty.
				
				\item If $v$ is an articulation point separating $H_u$ (the subgraph containing $u$) and $H_1,\ldots,H_k$, add $H_i,\{v\}$ for every $i \in \{1,\ldots,k\}$ to $O_v$.
				
				\item If $C$ is an even cycle, and $u$ and $v$ are diametrically opposed, then for every vertex $w \in C$ that is an articulation point (with $w \neq v$), add \texttt{Decompose($G$,$u$,$v$,$w$)} to $O_v$.
				
				\item Otherwise, for every vertex $w \in C$ on the shortest path between $u$ and $v$ that is an articulation point (with $w \neq v$), denote the output of \texttt{Decompose($G$,$u$,$v$,$w$)} by $H^w$, and add every element of $H^w$ to $O_v$. Finally, let $H$ be the graph obtained by removing every element of every $H^w$ from $G$. Add $H,\{u,v\}$ to $O_v$.
				
				\item Add $O_v$ to $O$.
			\end{enumerate}
			
			\item Otherwise, denote by $P$ a shortest path between $u$ and $v$. Denote its vertices by $w_0,w_1,\ldots,w_k$ with $w_0=u$ and $w_k=v$. Let $C_u$ be the cycle containing $u$ (if it exists), and $C_v$ be the cycle containing $v$ (if it exists). Then:
			\begin{enumerate}
				\item Create a list $O_v$, initially empty.
				
				\item If $v$ is an articulation point separating $H_u$ (the subgraph containing $u$) and $H_1,\ldots,H_k$, add $H_i,\{v\}$ for every $i \in \{1,\ldots,k\}$ to $O_v$.
				
				\item Assume that at least one of $C_u,C_v$ exists. Let $w_{i_u}$ be the last vertex of $P$ in $C_u$, and $w_{i_v}$ be the first vertex of $P$ in $C_v$. Denote by $H_u$ the subgraph containing $u$ separated by $w_{i_u}$, and by $H_v$ the subgraph containing $v$ separated by $w_{i_v}$ (if $v$ is an articulation point, remove from $H_v$ all the vertices separated from $w_{i_v}$ by $v$).
				
				For every $j \in \{1,\ldots,i_u\}$, if $w_j$ is an articulation point, remove from $H_u$ every vertex separated from $u$ by $w_j$. Denote by $H'_u$ the graph thus obtained.
				
				For every $j$ in $\{i_v,\ldots,k-1\}$, if $w_j$ is an articulation point, remove from $H_v$ every vertex separated from $v$ by $w_j$. Denote by $H'_v$ the graph thus obtained.
				
				If $C_u$ is an odd cycle, or if it is an even cycle but $u$ and $w_{i_u}$ are not diametrically opposed, then, add $H'_u,\{u,w_{i_u}\}$.
				
				If $C_v$ is an odd cycle, or if it is an even cycle but $v$ and $w_{i_v}$ are not diametrically opposed, then, add $H'_v,\{v,w_{i_v}\}$ to $O_v$.
				
				(Note: If $u$ is a leaf, skip the $H_u$ part. If $v$ is not in a cycle, skip the $H_v$ part.)
				
				\item For every $i \in \{1,\ldots,k-1\}$, if $w_i$ is an articulation point, add \texttt{Decompose($G$,$u$,$v$,$w_i$)} to $O_v$.
				
				\item For every maximal subpath of vertices $w_{i_b},\ldots,w_{i_e}$ in the same cycle $C$ (with $i_b < i_e$, $i_b,i_e \in \{0,\ldots,k\}$):
				\begin{enumerate}
					\item If $C$ is even, and $w_{i_b}$ and $w_{i_e}$ are diametrically opposed, denote the vertices of $C$ that are not in $P$ by $x_1,\ldots,x_\ell$. For every $j \in \{1,\ldots,\ell\}$, add \texttt{Decompose($G$,$u$,$v$,$x_j$)} to $O_v$.
					\item Otherwise, denote by $H$ the maximal connected subgraph of $G$ containing all vertices of $C \setminus P$. Denote by $H'$ the subgraph of $G$ obtained by adding $w_{i_b}$ and $w_{i_e}$ to $H$. Add $H',\{w_{i_b},w_{i_e}\}$ to $O_v$.
				\end{enumerate} 
				
				\item Add $O_v$ to $O$.
			\end{enumerate}
		\end{enumerate}
		
		\item For every option $o$ in $O$:
		\begin{enumerate}
			\item Create a list $T_o$, initially empty.
			
			\item For every element $o_i$ of $o$:
			\begin{enumerate}
				\item If $o_i$ is of the form $H,u$ where $H$ is a graph and $u$ a vertex, add \texttt{\grFCVA($H$,$u$)} to $T_o$.
				\item If $o_i$ is of the form $H,\{u,v\}$ where $H$ is a graph and $u$ and $v$ two vertices, add \texttt{\grFCVB($H$,$u$,$v$)} to $T_o$.
			\end{enumerate}
			
			\item Let $x_1,\ldots,x_k$ be the elements of $T_o$. Add $x_1 \oplus \ldots \oplus x_k$ to $T$.
		\end{enumerate}
		
		\item Return $\mex(T)$
	\end{enumerate}
	
	\medskip
	\noindent\textbf{The algorithm \texttt{\grFCVB($G$,$u$,$v$)}.} Recall that $u$ and $v$ are in the same cycle, and are neither articulation points nor diametrically opposed. Furthermore, no vertex on the shortest path between $u$ and $v$ is an articulation point.
	Denote the cycle containing $u$ and $v$ by $C$, and its vertices by $c_0,\ldots,c_{\ell-1}$. Let $u=c_0$ and $v=c_k$ with $k < \lceil \frac{\ell}{2} \rceil$.
	See \Cref{fig-cactus-4} for a depiction of the decomposition when it differs from the one in \texttt{\grFCVA($G$,$u$)}.
	
	\begin{enumerate}
		\item If $G$ is a tree, use~\cite[Algorithm 2]{araujo2024graph} and return the value. If $G$ is an odd cycle, return $\lfloor \frac{\ell}{2} \rfloor + 1 - k$. If $G$ is an even cycle, return $\frac{\ell}{2} - k$.
		
		\item Create a list $O$ of options and a list $T$ of integers, both initially empty.
		
		\item For every vertex $w$ of $G$ that is not on the shortest path between $u$ and $v$:
		\begin{enumerate}
			\item Create a list $O_w$, initially empty. Let $H$ be a copy of $G$.
			
			\item Denote by $x$ the vertex of $C$ that is the closest from $w$ (we may have $x=w$ if $w \in C$). For index purposes, denote $x = c_m$.
			
			\item If $x$ is an articulation point, let $H_1,\ldots,H_n$ be the components separated by $x$ that contain neither $u$ nor $w$. Add $H_i,\{x\}$ for every $i \in \{1,\ldots,n\}$ to $O_w$, and remove $H_i \setminus x$ from $H$.
			
			\item If $C$ is even and $x$ is diametrically opposed from either $u$ or $v$, then, for every $i \in \{k+1,\ldots,m-1,m+1,\ldots,\ell-1\}$ such that $c_i$ is an articulation point, add \texttt{Decompose($H$,$u$,$x$,$c_i$)} (or \texttt{Decompose($H$,$v$,$x$,$c_i$)} as needed) to $O_w$.
			
			\item Otherwise, denote by $P_{ux}$ and $P_{vx}$ the two unique shortest paths from $u$ to $x$ and from $v$ to $x$, respectively.
			\begin{enumerate}
				
				\item If $v \in P_{ux}$, then, let $H'$ be the maximal subgraph of $H$ containing no vertex $y$ separated from $u$ by $x$. For every $i \in \{k+1,\ldots,m-1\}$ such that $c_i$ is an articulation point, add \texttt{Decompose($H$,$v$,$x$,$c_i$)} to $O_w$, and remove from $H'$ every vertex that is separated from $u$ by $c_i$. Finally, add $H',\{u,x\}$ to $O_w$.
				
				\item Otherwise, if $u \in P_{vx}$, then, let $H'$ be the maximal subgraph of $H$ containing no vertex $y$ separated from $u$ by $x$. For every $i \in \{m+1,\ldots,\ell-1\}$ such that $c_i$ is an articulation point, add \texttt{Decompose($H$,$u$,$x$,$c_i$)} to $O_w$, and remove from $H'$ every vertex that is separated from $u$ by $c_i$. Finally, add $H',\{v,x\}$ to $O_w$.
				
				\item Otherwise: for every $i \in \{k+1,\ldots,m-1\}$ such that $c_i$ is an articulation point, add \texttt{Decompose($H$,$v$,$x$,$c_i$)} to $O_w$; for every $i \in \{m+1,\ldots,\ell-1\}$, such that $c_i$ is an articulation point, add \texttt{Decompose($H$,$u$,$x$,$c_i$)} to $O_w$.
			\end{enumerate}
			
			\item If $w \not\in C$, denote by $P$ a shortest path between $x$ and $w$. Denote its vertices by $y_0,y_1,\ldots,y_n$ with $y_0=x$ and $y_n=w$. Let Let $C_w$ be the cycle containing $w$ (if it exists). Then:
			\begin{enumerate}
				\item If $w$ is an articulation point separating $H_x$ (the subgraph containing $x$) and $H_1,\ldots,H_n$, add $H_i,\{w\}$ for every $i \in \{1,\ldots,n\}$ to $O_w$.
				
				\item If $C_w$ exists, let $y_{i_w}$ be the first vertex of $P$ in $C_w$. Denote by $H_w$ the subgraph containing $w$ separated by $y_{i_v}$ (if $w$ is an articulation point, remove from $H_w$ all the vertices separated from $y_{i_v}$ by $w$).
				
				For every $j$ in $\{i_w,\ldots,n-1\}$, if $y_j$ is an articulation point, remove from $H_w$ every vertex separated from $w$ by $y_j$. Denote by $H'_w$ the graph thus obtained.
				
				Add $H'_w,\{w,y_{i_v}\}$ to $O_w$.
				
				\item For every $i \in \{1,\ldots,n-1\}$, if $y_i$ is an articulation point, add \texttt{Decompose($H$,$x$,$w$,$y_i$)} to $O_w$.
				
				\item For every maximal subpath of vertices $y_{i_b},\ldots,y_{i_e}$ in the same cycle $D$ (with $i_b < i_e$, $i_b,i_e \in \{0,\ldots,n\}$):
				\begin{enumerate}
					\item If $D$ is even, and $y_{i_b}$ and $y_{i_e}$ are diametrically opposed, denote the vertices of $D$ that are not in $P$ by $z_1,\ldots,z_p$. For every $j \in \{1,\ldots,p\}$, add \texttt{Decompose($G$,$x$,$w$,$z_j$)} to $O_w$.
					\item Otherwise, denote by $H'$ the maximal connected subgraph of $H$ containing all vertices of $D \setminus P$. Denote by $H''$ the subgraph of $H$ obtained by adding $y_{i_b}$ and $y_{i_e}$ to $H'$. Add $H'',\{y_{i_b},y_{i_e}\}$ to $O_w$.
				\end{enumerate} 
			\end{enumerate}
			
			\item Add $O_w$ to $O$.
		\end{enumerate}
		
		\item For every option $o$ in $O$:
		\begin{enumerate}
			\item Create a list $T_o$, initially empty.
			
			\item For every element $o_i$ of $o$:
			\begin{enumerate}
				\item If $o_i$ is of the form $H,u$ where $H$ is a graph and $u$ a vertex, add \texttt{\grFCVA($H$,$u$)} to $T_o$.
				\item If $o_i$ is of the form $H,\{u,v\}$ where $H$ is a graph and $u$ and $v$ two vertices, add \texttt{\grFCVB($H$,$u$,$v$)} to $T_o$.
			\end{enumerate}
			
			\item Let $x_1,\ldots,x_k$ be the elements of $T_o$. Add $x_1 \oplus \ldots \oplus x_k$ to $T$.
		\end{enumerate}
		
		\item Return $\mex(T)$
	\end{enumerate}
	
	We now prove that the algorithms are correct. Since many arguments are close to those used for proving \Cref{thm-blockGraphs}, we are going to be less detailed.
	First, note that for all three of \texttt{\grFCG($G$)}, \texttt{\grFCVA($G$,$u$)} and \texttt{\grFCVB($G$,$u$,$v$)}, if $G$ is a tree, then~\cite[Algorithm~2]{araujo2024graph} outputs a correct solution, and if $G$ is a cycle, then \Cref{thm-cycles} and \Cref{prop-cyclesOneVertex} ensure that the output is correct.
	
	We now prove that the decompositions of each algorithm are correct, and that the correct properties hold for calling the subalgorithms.
	
	The selected articulation points are all managed with \Cref{cut_vertex_in_closure}. Furthermore, let $u$ and $v$ be two selected vertices, \Cref{cut_vertex_in_closure} can be applied on any vertex $w$ on a shortest path between $u$ and $v$ to get a disjoint sum of the components disconnected from $u$ and $v$ by $w$, each of those with $w$ selected. Indeed, let $x$ be a vertex inside of one such component $C$, then having $w$ selected is the same for $x$ as having any vertex outside of $C$ selected; conversely, let $y$ be a vertex in $G \setminus C$, all the shortest paths between $x$ and $y$ cover the same vertices of $G \setminus C$ than the shortest paths between $w$ and $y$ (since $w$ is an articulation point separating $x$ from $y$) and thus than shortest paths between $y$ and at least one of $u$ and $v$.
	This proves that \texttt{Decompose($G$,$u$,$v$,$w$)} is correct.
	
	Consider now the main algorithm and its two subroutines.
	\begin{itemize}
		\item \texttt{\grFCG($G$)} iterates over all possible options, which are all the different vertices. Each subroutine call has only one vertex selected, thus the calls are correct. Note furthermore that, in every call, if $u$ is an articulation point, then the graph is decomposed through \Cref{cut_vertex_in_closure}. Hence, $u$ can never be an articulation point, and thus it is either in a cycle or a leaf, allowing us to do the recursive calls to \texttt{\grFCVA}.
		
		Hence, for all $u$, $O_u$ contains the corresponding option (which can be a disjoint sum). Using the nim-sum gives the Sprague-Grundy value of each option, and computing the $\mex$ over all those values gives $\gr(G)$.
		
		\item \texttt{\grFCVA($G$,$u$)} iterates over all possible options, which are all vertices different from $u$. If $u$ is a leaf and the option $v$ is its neighbour, then, we can remove $u$ from the graph as $v$ is on all shortest paths between $u$ and any other vertex of $G$, and thus we iterate on the selection of $v$ (Step 3(a), with substep i if $v$ is an articulation point and ii otherwise).
		
		If $u$ and $v$ are in the same cycle (Step 3(b)), we first manage the case of $v$ being an articulation point (substep ii), and then there are two cases. First (substep iii), if the cycle is even and $u$ and $v$ are diametrically opposed, then no other vertex from the cycle can be selected anymore since the union of the two shortest paths between $u$ and $v$ is the cycle, so we only have to decompose along all articulation points of the cycle and add the sum to the options set. Second (substep iv) is all other cases: there is only one shortest path between $u$ and $v$, and thus we decompose along it, and the remaining non-decomposed graph $H$ (containing both $u$ and $v$) contains exactly two selected vertices, so we add the sum of $H$ and all decomposed graphs to the options set.
		
		Finally, if $u$ and $v$ are in different cycles (Step 3(c)), we need to follow a shortest path $P$ between them, which can go throughout subtrees and other cycles, some of which can be even and where the first and last vertex of $P$ in these even cycles can be diametrically opposed (so there may be many different shortest paths between $u$ and $v$).
		We first decompose along $v$ if it is an articulation point (substep ii), adding those subgraphs to the options set.
		We then assume that one of $u$ and $v$ is in a cycle (substep iii). Without loss of generality, assume that $u$ is in a cycle $C_u$, and let $w_{i_u}$ be the last vertex of $P$ in $C_u$. If $C_u$ is even, and $u$ and $w_{i_u}$ are diametrically opposed, then again it is impossible to play on $C_u$ anymore in future moves, so we decompose along all vertices of $P \cap C_u$ (the other vertices of $C_u$ will be managed in substep v), and we only add those decomposed graphs to the options set ($C_u$ will not appear in the options set). Otherwise, we again decompose along vertices of $P \cap C_u$ and add all those to the options set, but we also add to this set the maximal subgraph of $G$ containing $u$ and $w_{i_u}$ but no vertex separated from $u$ by a vertex of $P \cap C_u$. Denote this last subgraph by $H'_u$ and denote by $H$ the graph separated from $H'_u$ by $w_{i_u}$, clearly selecting a vertex $x$ in $H'_u$ will be equivalent to selecting $u$ for any vertex in $H$; and conversely selecting a vertex $x$ outside of $H'_u$ is equivalent to having $w_{i_u}$ selected, as $w_{i_u}$ is on every shortest path between $u$ and $x$. So $H'_u,\{u,w_{i_u}\}$ is an option.
		The same reasoning applies if $v$ is in a cycle.
		
		Substep iv decomposes along articulation points of $P$ and adds correct options to the set by using \texttt{Decompose}.
		As for substep v, it consists in managing intermediary cycles along $P$. Consider a maximal subpath $w_{i_b} \ldots w_{i_e}$ of $P$ that intersects a given cycle $C$. We apply the same reasoning as the paragraph above, substituting $u$ with $w_{i_b}$ and $w_{i_u}$ with $w_{i_e}$.
		
		Thus, for each vertex $v$, at the end of Step 3, the set $O_v$ contains a disjoint sum of options that have been decomposed correctly.
		Again, all recursive calls to \texttt{\grFCVA} are on vertices that are not articulation points (by prealably decomposing if needed). Furthermore, in all calls to \texttt{\grFCVB}, the two vertices are in the same cycle and are not articulation points, and no vertex in the shortest path between them is an articulation point (also by decomposing if needed).
		Hence, the calls verify all properties of the subalgorithms, so the recursive calls will output correct Sprague-Grundy values and allow us to compute the nim-sum for each set $O_v$, over which the $\mex$ will give $\gr(G,\{u\})$.
		
		\item Consider now \texttt{\grFCVB($G$,$u$,$v$)}. Since it is very similar to the previous subalgorithm, we are going to give fewer details. Recall that $u$ and $v$ are in the same cycle $C$, and that no vertex on their shortest path (including them) is an articulation point. We again iterate on every option $w$, which can be any vertex in $G$ that is not on the shortest path between $u$ and $v$. Let $x$ be the vertex of $C$ that is closest to $w$.
		The principle is to manage independently $C$ and the path between $x$ and $w$, the latter (Step 3(f)) being the same as Step 3(c) from \texttt{\grFCVA} (substituting $u$ with $x$ and $v$ with $w$). For the former, there are three possibilities.
		
		The first possibility (Step 3(d)) is that $C$ is even and $x$ is diametrically opposed to one of $u$ or $v$, in which case no vertex from $C$ can be played on anymore in future options. In this case, we decompose along the articulation points of $C$ to obtain a disjoint sum.
		
		The second possibility (Step (3e), substep iii) is that $x$ is diametrically opposed to a vertex on the shortest path between $u$ and $v$, in which case the shortest path between $u$ and $x$ and the shortest path between $v$ and $x$ cover together all vertices in $C$. Again, we cannot select a vertex of $C$ in future moves, and we decompose along the articulation points of $C$ to obtain a disjoint sum.
		
		The third possibility (Step (3e), substeps i and ii) is that (without loss of generality) $v$ is on the shortest path between $u$ and $x$. In this case, we decompose along the shortest path between $v$ and $x$ (those vertices cannot be selected anymore) to obtain a disjoint sum, and we then manage the cycle using \texttt{\grFCVB} (substituting $v$ with $x$), since all other vertices are still options for future moves.
		
		Thus, for each vertex $w$, at the end of Step 3, the set $O_w$ contains a disjoint sum of options that have been decomposed correctly.
		Again, all recursive calls to \texttt{\grFCVA} are on vertices that are not articulation points (by prealably decomposing if needed). Furthermore, in all calls to \texttt{\grFCVB}, the two vertices are in the same cycle and are not articulation points, and no vertex in the shortest path between them is an articulation point (also by decomposing if needed).
		Hence, the calls verify all properties of the subalgorithms, so the recursive calls will output correct Sprague-Grundy values and allow us to compute the nim-sum for each set $O_w$, over which the $\mex$ will give $\gr(G,\{u,v\})$.
	\end{itemize}
	
	Since the base cases (cycles, trees, empty graphs) are correct and the recursive calls are allowed due to verifying the required properties and the decompositions being correct, the algorithm \texttt{\grFCG($G$)} outputs $\gr(G)$. All that remains is to show that the running time is polynomial.
	
	As with block graphs, we now have an algorithm that will create the game graph of $G$: each vertex is a position, which allows to compute only once each Sprague-Grundy value. The running time hence depends on the cost and number of decompositions. By encoding cycles and vertices in trees in an arborescent manner, it becomes easier to manage the decomposition at an articulation point. Hence, the running time of the decomposition along an articulation point (such as step 3(a) of \texttt{GrundyCactus($G$)} or \texttt{Decompose($G$,$u$,$v$,$w$)}) is linear; managing steps along a potentially linear number of vertices (such as following the shortest path between $u$ and $v$, and thus steps 3(b) or 3(c) of \texttt{\grFCVA($G$,$u$)}) is thus quadratic. As the number of options as well as the number of cycles is at most linear, this means that we maintain a polynomial running time.
\end{proof}

\begin{figure}[h]
	\centering
	\begin{tikzpicture}
		\node (base) at (0,0) {
			\begin{tikzpicture}
				\draw (0,0) ellipse (1 and 0.5);
				\draw (0,0) node {$C_1$};
				\draw (1.75,0) ellipse (0.75 and 0.375);
				\draw (1.75,0) node {$C_2$};
				\node[noeud,fill=black] (u) at (1,0) {};
				\draw (u) node[below,yshift=-1.5mm] {$u$};
				\draw (u)to(1,1);
				\draw (0.25,1.75)to(1,1)to(1.5,1.5);
				\draw (2,1.5) ellipse (0.5 and 0.25);
				\draw (2,1.5) node {$C_3$};
			\end{tikzpicture}
		};
		
		\node (sum1) at (5,1) {
			\begin{tikzpicture}
				\draw (0,0) ellipse (1 and 0.5);
				\draw (0,0) node {$C_1$};
				\node[noeud,fill=black] (u) at (1,0) {};
				\draw (u) node[below,yshift=-1.5mm] {$u$};
			\end{tikzpicture}
		};
		
		\node (sum2) at (5,-1) {
			\begin{tikzpicture}
				\draw (1.75,0) ellipse (0.75 and 0.375);
				\draw (1.75,0) node {$C_2$};
				\node[noeud,fill=black] (u) at (1,0) {};
				\draw (u) node[below,yshift=-1.5mm] {$u$};
			\end{tikzpicture}
		};
		
		\node (sum3) at (8,0) {
			\begin{tikzpicture}
				\node[noeud,fill=black] (u) at (1,0) {};
				\draw (u) node[below,yshift=-1.5mm] {$u$};
				\draw (u)to(1,1);
				\draw (0.25,1.75)to(1,1)to(1.5,1.5);
				\draw (2,1.5) ellipse (0.5 and 0.25);
				\draw (2,1.5) node {$C_3$};
			\end{tikzpicture}
		};
		
		\draw (2.75,0) node {$\equiv$};
		\draw (5,0) node {$+$};
		\draw (6.5,0) node {$+$};
	\end{tikzpicture}
	\caption{The decomposition of a cactus (each ellipsis is a cycle): selecting an articulation point.}
	\label{fig-cactus-1}
\end{figure}

\begin{figure}[h]
	\centering
	\begin{tikzpicture}
		\node (input) at (0,0) {
			\begin{tikzpicture}
				\node[noeud,fill=black] (u) at (0,0) {};
				\node[noeud,fill=black] (v) at (3,0) {};
				\node[noeud,cross out] (w) at (1.5,0) {};
				\draw (u) node[below,yshift=-1.5mm] {$u$};
				\draw (v) node[below,yshift=-1.5mm] {$v$};
				\draw (w) node[below,xshift=-1.5mm,yshift=-1.5mm] {$w$};
				\draw (u)to(w)to(v);
				
				\draw (1.5,0.75) ellipse (0.5 and 0.75);
				\draw (2.5,-0.75) ellipse (0.5 and 0.25);
				\draw (w)to(2,-0.75);
			\end{tikzpicture}
		};
		
		\node (output) at (6,0) {
			\begin{tikzpicture}
				\draw (0,0) ellipse (0.5 and 0.75);
				\node[noeud,fill=black] (w1) at (0,-0.75) {};
				\draw (w1) node[below,yshift=-1.5mm] {$w$};
				
				\node[noeud,fill=black] (w2) at (1.5,0) {};
				\draw (2.5,-0.75) ellipse (0.5 and 0.25);
				\draw (w2)to(2,-0.75);
				\draw (w2) node[above,yshift=1.5mm] {$w$};
				
				\draw (1,-0.25) node {+};
			\end{tikzpicture}
		};
		
		\draw[->,double,line width=0.25mm] (input)to(output);
	\end{tikzpicture}
	\caption{The subroutine \texttt{Decompose($G$,$u$,$v$,$w$)} where $u$ and $v$ are selected, and $w$ is on a shortest path between $u$ and $v$. It outputs a list of subgames.}
	\label{fig-cactus-2}
\end{figure}

\begin{figure}[h]
	\centering
	\begin{subfigure}{0.49\linewidth}
		\centering
		\begin{tikzpicture}
			\node (first) at (0,0) {
				\begin{tikzpicture}
					\draw (0,0) ellipse (1.5 and 0.75);
					\node[noeud,fill=black] (u) at (-1.5,0) {};
					\draw (u) node[left,xshift=-1.5mm] {$u$};
					\node[noeud,fill=black] (v) at (1.5,0) {};
					\draw (v) node[right,xshift=1.5mm] {$v$};
					
					\draw (-0.5,1.2) ellipse (0.25 and 0.5);
					\node[noeud] (w1) at (-0.5,0.7) {};
					\draw (w1) node[below,yshift=-1.5mm] {$w_1$};
					
					\node[noeud] (w2) at (0.5,0.7) {};
					\draw (w2) node[below,yshift=-1.5mm] {$w_2$};
					\draw (w2)to(0.5,1.2);
					\draw (0.2,1.4)to(0.5,1.2)to(1,1.7);
					
					\draw (0,-1.25) ellipse (0.25 and 0.5);
					\draw (0.25,-1.25)to(0.75,-1.25);
					\draw (1.05,-1.25) ellipse (0.3 and 0.15);
					\node[noeud] (w3) at (0,-0.75) {};
					\draw (w3) node[above,yshift=1.5mm] {$w_3$};
				\end{tikzpicture}
			};
			
			\node (second1) at (-2.5,-4) {
				\begin{tikzpicture}
					\draw (-0.5,1.2) ellipse (0.25 and 0.5);
					\node[noeud,fill=black] (w1) at (-0.5,0.7) {};
					\draw (w1) node[below,yshift=-1.5mm] {$w_1$};
				\end{tikzpicture}
			};
			
			\node (second2) at (-0.5,-4) {
				\begin{tikzpicture}
					\node[noeud,fill=black] (w2) at (0.5,0.7) {};
					\draw (w2) node[below,yshift=-1.5mm] {$w_2$};
					\draw (w2)to(0.5,1.2);
					\draw (0.2,1.4)to(0.5,1.2)to(1,1.7);
				\end{tikzpicture}
			};
			
			\node (second3) at (2.5,-4) {
				\begin{tikzpicture}
					\draw (0,-1.25) ellipse (0.25 and 0.5);
					\draw (0.25,-1.25)to(0.75,-1.25);
					\draw (1.05,-1.25) ellipse (0.3 and 0.15);
					\node[noeud,fill=black] (w3) at (0,-0.75) {};
					\draw (w3) node[above,yshift=1.5mm] {$w_3$};
				\end{tikzpicture}
			};
			
			\draw (0,-2.5) node {$\equiv$};
			\draw (-1.5,-4) node {$+$};
			\draw (0.75,-4) node {$+$};
		\end{tikzpicture}
		\caption{Depiction of substep (b)-iii: $u$ and $v$ are in the same cycle, which is even, and are diametrically opposed. In this case, playing on the cycle will be impossible in future moves.}
		\label{fig-cactus-3-1}
	\end{subfigure}
	\hfil
	\begin{subfigure}{0.49\linewidth}
		\centering
		\begin{tikzpicture}
			\node (first) at (0,0) {
				\begin{tikzpicture}
					\draw (0,0) ellipse (1.5 and 0.75);
					\node[noeud,fill=black] (u) at (-1.25,0.4) {};
					\draw (u) node[left,xshift=-1.5mm] {$u$};
					\node[noeud,fill=black] (v) at (1.5,0) {};
					\draw (v) node[right,xshift=1.5mm] {$v$};
					
					\draw (-0.5,1.2) ellipse (0.25 and 0.5);
					\node[noeud] (w1) at (-0.5,0.7) {};
					\draw (w1) node[below,yshift=-1.5mm] {$w_1$};
					
					\node[noeud] (w2) at (0.5,0.7) {};
					\draw (w2) node[below,yshift=-1.5mm] {$w_2$};
					\draw (w2)to(0.5,1.2);
					\draw (0.2,1.4)to(0.5,1.2)to(1,1.7);
					
					\draw (0,-1.25) ellipse (0.25 and 0.5);
					\draw (0.25,-1.25)to(0.75,-1.25);
					\draw (1.05,-1.25) ellipse (0.3 and 0.15);
					\node[noeud] (w3) at (0,-0.75) {};
					\draw (w3) node[above,yshift=1.5mm] {$w_3$};
				\end{tikzpicture}
			};
			
			\node (second1) at (-1.75,-4.5) {
				\begin{tikzpicture}
					\draw (0,0) ellipse (1.5 and 0.75);
					\node[noeud,fill=black] (u) at (-1.25,0.4) {};
					\draw (u) node[left,xshift=-1.5mm] {$u$};
					\node[noeud,fill=black] (v) at (1.5,0) {};
					\draw (v) node[right,xshift=1.5mm] {$v$};
					
					\node (w1) at (-0.5,0.7) {};
					\draw (w1) node[below,yshift=-1mm] {$w_1$};
					\node (w2) at (0.5,0.7) {};
					\draw (w2) node[below,yshift=-1mm] {$w_2$};
					
					\foreach \I in {(-1,0.55),(-0.75,0.65),(-0.5,0.7),(-0.25,0.725),(0,0.75),(0.25,0.725),(0.5,0.7),(0.75,0.65),(1,0.55),(1.25,0.4),(1.425,0.2)} {
						\node[noeud,cross out] at \I {};
					}
					
					\draw (0,-1.25) ellipse (0.25 and 0.5);
					\draw (0.25,-1.25)to(0.75,-1.25);
					\draw (1.05,-1.25) ellipse (0.3 and 0.15);
					\node[noeud] (w3) at (0,-0.75) {};
					\draw (w3) node[above,yshift=1.5mm] {$w_3$};
				\end{tikzpicture}
			};
			
			\node (second2) at (1.25,-4.5) {
				\begin{tikzpicture}
					\draw (-0.5,1.2) ellipse (0.25 and 0.5);
					\node[noeud,fill=black] (w1) at (-0.5,0.7) {};
					\draw (w1) node[below,yshift=-1.5mm] {$w_1$};
				\end{tikzpicture}
			};
			
			\node (second3) at (3,-4.5) {
				\begin{tikzpicture}
					\node[noeud,fill=black] (w2) at (0.5,0.7) {};
					\draw (w2) node[below,yshift=-1.5mm] {$w_2$};
					\draw (w2)to(0.5,1.2);
					\draw (0.2,1.4)to(0.5,1.2)to(1,1.7);
				\end{tikzpicture}
			};
			
			\draw (0,-2.5) node {$\equiv$};
			\draw (0.5,-4.5) node {$+$};
			\draw (2.125,-4.5) node {$+$};
		\end{tikzpicture}
		\caption{Depiction of substep (b)-iv: $u$ and $v$ are in the same cycle and are not diametrically opposed. In this case, playing on the crossed out vertices of the cycle will be impossible in future moves.}
		\label{fig-cactus-3-2}
	\end{subfigure}
	
	\begin{subfigure}{0.49\linewidth}
		\centering
		\begin{tikzpicture}
			\node (first) at (0,0) {
				\begin{tikzpicture}
					\draw[line width=0.5mm] (1,-1)to(1,1)to(3,1)to(3,-1)to(1,-1);
					\draw (0,0) ellipse (1 and 0.5);
					\draw (-0.5,0) node {$C_u$};
					\node[noeud,fill=black] (u) at (0,0.5) {};
					\draw (u) node[above,yshift=1mm] {$u$};
					\node[noeud] (wiu) at (1,0) {};
					\draw (wiu) node[left,xshift=0.5mm] {$w_{i_u}$};
					\node[noeud] (wi) at (0.5,0.425) {};
					\draw (wi) node[above right,xshift=-0.5mm,yshift=-0.5mm] {$w_i$};
					\draw (wi)to(0.5,0.9);
					\draw  (0.2,1.2)to(0.5,0.9)to(0.9,1.5);
					\draw (4,0) ellipse (1 and 0.5);
					\draw (4.5,0) node {$C_v$};
					\node[noeud,fill=black] (v) at (4,0.5) {};
					\draw (v) node[above,yshift=1mm] {$v$};
					\node[noeud] (wiv) at (3,0) {};
					\draw (wiv) node[right,xshift=0.5mm] {$w_{i_v}$};
					\draw[line width=0.5mm] (wiu)to(wiv);
					\draw (2,0.25) node {$P$};
				\end{tikzpicture}
			};
			
			\node (second) at (0,-4) {
				\begin{tikzpicture}
					\draw (0,0) ellipse (1 and 0.5);
					\draw (-0.5,0) node {$C_u$};
					\node[noeud,fill=black] (u) at (0,0.5) {};
					\draw (u) node[above,yshift=1mm] {$u$};
					\node[noeud,fill=black] (wiu) at (1,0) {};
					\draw (wiu) node[left,xshift=0.5mm] {$w_{i_u}$};
					\node[noeud,cross out] (wi) at (0.5,0.425) {};
					\draw (wi) node[above,yshift=-0.5mm] {$w_i$};
					\draw (4,0) ellipse (1 and 0.5);
					\draw (4.5,0) node {$C_v$};
					\node[noeud,fill=black] (v) at (4,0.5) {};
					\draw (v) node[above,yshift=1mm] {$v$};
					\node[noeud,fill=black] (wiv) at (3,0) {};
					\draw (wiv) node[right,xshift=0.5mm] {$w_{i_v}$};
					\foreach \I in {(0.25,0.475),(0.8,0.3),(3.2,0.3),(3.5,0.425),(3.75,0.475)} {
						\node[noeud,cross out] at \I {};
					}
					
					\node[noeud,fill=black] (wip) at (2,0) {};
					\draw (wip) node[below,yshift=-1.5mm] {$w_i$};
					\draw (wip)to(2,0.5);
					\draw  (1.7,0.8)to(2,0.5)to(2.4,1);
					
					\draw (1.5,0) node {+};
					\draw (2.5,0) node {+};
				\end{tikzpicture}
			};
			
			\draw[->,double,line width=0.25mm] (first)to(second);
		\end{tikzpicture}
		\caption{Depiction of substep (c)-iii: $u$ and $v$ are in different cycles, and we manage their own cycles $C_u$ and $C_v$ (which both exist in this example). The box is managed in substep (c)-v.}
		\label{fig-cactus-3-3}
	\end{subfigure}
	\hfil
	\begin{subfigure}{0.49\linewidth}
		\centering
		\begin{tikzpicture}
			\node (first) at (0,0) {
				\begin{tikzpicture}
					\draw[line width=0.5mm] (-3.5,-0.5)to(-3.5,1)to(-0.8,1)to(-0.8,0.3)to(-1.25,0)to(-1.25,-0.5)to(-3.5,-0.5);
					\draw[line width=0.5mm] (3.5,-0.5)to(3.5,1)to(1,1)to(1,-0.5)to(3.5,-0.5);
					\draw (0,0) ellipse (1 and 0.5);
					\node[noeud] (wib) at (-0.8,0.3) {};
					\draw (wib) node[below right] {$w_{i_b}$};
					\node[noeud] (wie) at (1,0) {};
					\draw (wie) node[left] {$w_{i_e}$};
					\draw (0,1) ellipse (0.25 and 0.5);
					\node[noeud] (wi) at (0,0.5) {};
					\draw (0,-0.5)to(0,-1);
					\draw (-0.75,-1) ellipse (0.5 and 0.25);
					\draw (1,-1) ellipse (0.75 and 0.375);
					\draw (-0.25,-1)to(0.25,-1);
					\node[noeud,fill=black] (u) at (-3,0) {};
					\draw (u) node[above,yshift=1mm] {$u$};
					\node[noeud,fill=black] (v) at (3,0) {};
					\draw (v) node[above,yshift=1mm] {$v$};
					\draw[line width=0.5mm] (u)to(-2,0)to(-2,0.3)to(wib);
					\draw[line width=0.5mm] (v)to(wie);
					\draw (-2,0.55) node {$P$};
					\draw (2,0.25) node {$P$};
				\end{tikzpicture}
			};
			
			\node (second) at (0,-4) {
				\begin{tikzpicture}
					\draw (0,0) ellipse (1 and 0.5);
					\node[noeud,fill=black] (wib) at (-0.8,0.3) {};
					\draw (wib) node[below right] {$w_{i_b}$};
					\node[noeud,fill=black] (wie) at (1,0) {};
					\draw (wie) node[left] {$w_{i_e}$};
					\draw (0,-0.5)to(0,-1);
					\draw (-0.75,-1) ellipse (0.5 and 0.25);
					\draw (1,-1) ellipse (0.75 and 0.375);
					\draw (-0.25,-1)to(0.25,-1);
					\foreach \I in {(-0.5,0.425),(-0.25,0.475),(0,0.5),(0.25,0.475),(0.5,0.425),(0.8,0.3)} {
						\node[noeud,cross out] at \I {};
					}
					
					\draw (3,0.25) ellipse (0.25 and 0.5);
					\node[noeud,fill=black] (wi) at (3,-0.25) {};
					
					\draw (2.25,0) node {+};
				\end{tikzpicture}
			};
			
			\draw[->,double,line width=0.25mm] (first)to(second);
		\end{tikzpicture}
		\caption{Depiction of substep (c)-v: $u$ and $v$ are in different cycles, and we manage, one by one, every cycle along a shortest path $P$ between them. Here, we assume that $w_{i_b}$ and $w_{i_e}$ are not diametrically opposed (otherwise, the whole cycle is removed from the graph and both sides are decomposed). The box on the left was already treated before, and the box on the right will be treated after, until we reach $C_v$ (or $v$, if $C_v$ does not exist).}
		\label{fig-cactus-3-4}
	\end{subfigure}
	
	\caption{Gradual decomposition given by step \textbf{3} of \texttt{\grFCVA($G$,$u$)} when $u$ is in a simple cycle. Here, we consider different vertices $v$.}
	\label{fig-cactus-3}
\end{figure}

\begin{figure}[h]
	\centering
	
	\begin{subfigure}{\linewidth}
		\centering
		\begin{tikzpicture}
			\node (first) at (0,0) {
				\begin{tikzpicture}
					\draw (0,0) ellipse (1.5 and 1);
					\draw (0,0) node {$C$};
					\node[noeud,fill=black] (u) at (-1.5,0) {};
					\draw (u) node[left,xshift=-1mm] {$u$};
					\node[noeud,fill=black] (v) at (0,1) {};
					\draw (v) node[above,yshift=0.5mm] {$v$};
					\foreach \I in {(-0.25,0.975),(-0.5,0.95),(-0.75,0.875),(-1,0.75),(-1.2,0.6),(-1.4,0.35)} {
						\node[noeud,cross out] at \I {};
					}
					\node[noeud] (x) at (1.5,0) {};
					\draw (x) node[left,xshift=-1mm] {$x$};
					\node[noeud,fill=black] (w) at (3,0) {};
					\draw (w) node[above,yshift=0.5mm] {$w$};
					\draw[line width=0.5mm] (w)to(x);
					\node[noeud] (y1) at (1,0.725) {};
					\draw (y1) node[below,yshift=-0.5mm] {$y_1$};
					\draw (y1)to(1.5,1.225);
					\draw (2,1.225) ellipse (0.5 and 0.25);
					\node[noeud] (y2) at (1,-0.725) {};
					\draw (y2) node[above,yshift=0.5mm] {$y_2$};
					\draw (y2)to(1,-1.025);
					\draw (0.8,-1.225)to(1,-1.025)to(1.3,-1.325);
					\draw (0,-1.5) ellipse (0.25 and 0.5);
					\node[noeud] (y3) at (0,-1) {};
					\draw (y3) node[above,yshift=0.5mm] {$y_3$};
				\end{tikzpicture}
			};
			
			\node (second) at (6,0) {
				\begin{tikzpicture}
					\node[noeud,fill=black] (y1) at (0,0) {};
					\draw (y1) node[below,yshift=-0.5mm] {$y_1$};
					\draw (y1)to(0.5,0.5);
					\draw (1,0.5) ellipse (0.5 and 0.25);
					
					\node[noeud,fill=black] (y2) at (2.5,0) {};
					\draw (y2) node[above,yshift=0.5mm] {$y_2$};
					\draw (y2)to(2.5,-0.3);
					\draw (2.3,-0.5)to(2.5,-0.3)to(2.8,-0.6);
					
					\draw (4,-0.5) ellipse (0.25 and 0.5);
					\node[noeud,fill=black] (y3) at (4,0) {};
					\draw (y3) node[above,yshift=0.5mm] {$y_3$};
					
					\draw (1.75,0) node {+};
					\draw (3.25,0) node {+};
				\end{tikzpicture}
			};
			
			\draw[->,double,line width=0.25mm] (first)to(second);
		\end{tikzpicture}
		\caption{Depiction of substep (d): $x$ is diametrically opposed to $u$. After this move, playing on $C$ is impossible.}
		\label{fig-cactus-4-1}
	\end{subfigure}
	
	\begin{subfigure}{\linewidth}
		\centering
		\begin{tikzpicture}
			\node (first) at (0,0) {
				\begin{tikzpicture}
					\draw (0,0) ellipse (1.5 and 1);
					\draw (0,0) node {$C$};
					\node[noeud,fill=black] (u) at (-1,0.75) {};
					\draw (u) node[left,xshift=-1mm] {$u$};
					\node[noeud,fill=black] (v) at (0,1) {};
					\draw (v) node[above,yshift=0.5mm] {$v$};
					\foreach \I in {(-0.25,0.975),(-0.5,0.95),(-0.75,0.875)} {
						\node[noeud,cross out] at \I {};
					}
					\node[noeud] (x) at (1.5,0) {};
					\draw (x) node[left,xshift=-1mm] {$x$};
					\node[noeud,fill=black] (w) at (3,0) {};
					\draw (w) node[above,yshift=0.5mm] {$w$};
					\draw[line width=0.5mm] (w)to(x);
					\node[noeud] (y1) at (1,0.725) {};
					\draw (y1) node[below,yshift=-0.5mm] {$y_1$};
					\draw (y1)to(1.5,1.225);
					\draw (2,1.225) ellipse (0.5 and 0.25);
					\node[noeud] (y2) at (1,-0.725) {};
					\draw (y2) node[above,yshift=0.5mm] {$y_2$};
					\draw (y2)to(1,-1.025);
					\draw (0.8,-1.225)to(1,-1.025)to(1.3,-1.325);
					\draw (0,-1.5) ellipse (0.25 and 0.5);
					\node[noeud] (y3) at (0,-1) {};
					\draw (y3) node[above,yshift=0.5mm] {$y_3$};
				\end{tikzpicture}
			};
			
			\node (second) at (6,0) {
				\begin{tikzpicture}
					\draw (0,0) ellipse (1.5 and 1);
					\draw (0,0) node {$C$};
					\node[noeud,fill=black] (u) at (-1,0.75) {};
					\draw (u) node[left,xshift=-1mm] {$u$};
					\node[noeud,cross out] (v) at (0,1) {};
					\draw (v) node[above,yshift=0.5mm] {$v$};
					\foreach \I in {(1.2,0.6),(1.4,0.35),(0.75,0.875),(0.5,0.95),(0.25,0.975),(-0.25,0.975),(-0.5,0.95),(-0.75,0.875)} {
						\node[noeud,cross out] at \I {};
					}
					\node[noeud,fill=black] (x) at (1.5,0) {};
					\draw (x) node[left,xshift=-1mm] {$x$};
					\node[noeud,cross out] (y1) at (1,0.725) {};
					\draw (y1) node[below,yshift=-0.5mm] {$y_1$};
					\node[noeud] (y2) at (1,-0.725) {};
					\draw (y2) node[above,yshift=0.5mm] {$y_2$};
					\draw (y2)to(1,-1.025);
					\draw (0.8,-1.225)to(1,-1.025)to(1.3,-1.325);
					\draw (0,-1.5) ellipse (0.25 and 0.5);
					\node[noeud] (y3) at (0,-1) {};
					\draw (y3) node[above,yshift=0.5mm] {$y_3$};
					
					\node[noeud,fill=black] (y1) at (3,0) {};
					\draw (y1) node[below,yshift=-0.5mm] {$y_1$};
					\draw (y1)to(3.5,0.5);
					\draw (4,0.5) ellipse (0.5 and 0.25);
					
					\draw (2.25,0) node {+};
				\end{tikzpicture}
			};
			
			\draw[->,double,line width=0.25mm] (first)to(second);
		\end{tikzpicture}
		\caption{Depiction of substep (e)-i: $v \in P_{ux}$.}
		\label{fig-cactus-4-2}
	\end{subfigure}
	
	\begin{subfigure}{\linewidth}
		\centering
		\begin{tikzpicture}
			\node (first) at (0,0) {
				\begin{tikzpicture}
					\draw (0,0) ellipse (1.5 and 1);
					\draw (0,0) node {$C$};
					\node[noeud,fill=black] (u) at (-1,-0.75) {};
					\draw (u) node[below,yshift=-0.5mm] {$u$};
					\node[noeud,fill=black] (v) at (0,1) {};
					\draw (v) node[above,yshift=0.5mm] {$v$};
					\foreach \I in {(-0.25,0.975),(-0.5,0.95),(-0.75,0.875),(-1,0.75),(-1.2,0.6),(-1.4,0.35),(-1.5,0),(-1.2,-0.6),(-1.4,-0.35)} {
						\node[noeud,cross out] at \I {};
					}
					\node[noeud] (x) at (1.5,0) {};
					\draw (x) node[left,xshift=-1mm] {$x$};
					\node[noeud,fill=black] (w) at (3,0) {};
					\draw (w) node[above,yshift=0.5mm] {$w$};
					\draw[line width=0.5mm] (w)to(x);
					\node[noeud] (y1) at (1,0.725) {};
					\draw (y1) node[below,yshift=-0.5mm] {$y_1$};
					\draw (y1)to(1.5,1.225);
					\draw (2,1.225) ellipse (0.5 and 0.25);
					\node[noeud] (y2) at (1,-0.725) {};
					\draw (y2) node[above,yshift=0.5mm] {$y_2$};
					\draw (y2)to(1,-1.025);
					\draw (0.8,-1.225)to(1,-1.025)to(1.3,-1.325);
					\draw (0,-1.5) ellipse (0.25 and 0.5);
					\node[noeud] (y3) at (0,-1) {};
					\draw (y3) node[above,yshift=0.5mm] {$y_3$};
				\end{tikzpicture}
			};
			
			\node (second) at (6,0) {
				\begin{tikzpicture}
					\node[noeud,fill=black] (y1) at (0,0) {};
					\draw (y1) node[below,yshift=-0.5mm] {$y_1$};
					\draw (y1)to(0.5,0.5);
					\draw (1,0.5) ellipse (0.5 and 0.25);
					
					\node[noeud,fill=black] (y2) at (2.5,0) {};
					\draw (y2) node[above,yshift=0.5mm] {$y_2$};
					\draw (y2)to(2.5,-0.3);
					\draw (2.3,-0.5)to(2.5,-0.3)to(2.8,-0.6);
					
					\draw (4,-0.5) ellipse (0.25 and 0.5);
					\node[noeud,fill=black] (y3) at (4,0) {};
					\draw (y3) node[above,yshift=0.5mm] {$y_3$};
					
					\draw (1.75,0) node {+};
					\draw (3.25,0) node {+};
				\end{tikzpicture}
			};
			
			\draw[->,double,line width=0.25mm] (first)to(second);
		\end{tikzpicture}
		\caption{Depiction of substep (e)-iii: after this move, playing on $C$ is impossible.}
		\label{fig-cactus-4-3}
	\end{subfigure}
	
	\caption{Gradual decomposition given by step \textbf{3} of \texttt{\grFCVB($G$,$u$,$v$)}. Recall that both vertices are in the same cycle $C$, and are not diametrically opposed. For every vertex $w$, we denote by $x$ the vertex in $C$ that is the closest to $w$. Substep (f) managing the path between $x$ and $w$ is similar to substep (c) managing the path between $u$ and $v$ from \texttt{\grFCVA($G$,$u$)}.}
	\label{fig-cactus-4}
\end{figure}

\section{Future research}

In this paper, we explored \AGAG in three directions: general results such as graph operations (how to get the outcome of the Cartesian product of two graphs in function of the outcomes of the individual graphs) and decomposition through articulation points, characterization of outcomes and Sprague-Grundy values, and polynomial-time algorithms based on dynamic programming for constrained graph classes. However, many more general graph classes remain open, algorithmically as well as structurally.

A first interesting research direction would be to fully characterize the graphs for which \AGAG is a \lovesme game, or those that are close to such a situation. We already exhibited complete graphs and stars as \lovesme games in \Cref{sec-lovesMeLovesMeNot}, but there might be other very structured classes with many simplicial vertices that are in the same situation. Note that the unicity of minimum-size geodetic sets and the structural backbone of the graph are not sufficient conditions for this property.

Another idea would be to further the work initiated in \Cref{sec-symmetry} and to define a type of graph symmetry adapted to \AGAG, that would allow us to obtain the outcomes of other structured graph classes. Expansions of cycles could be good candidates.

As explained in \Cref{sec-blockAndCacti}, the dynamic programming framework we used for block graphs and cacti requires strong constraints for shortest paths between vertices, which do not hold for more general graph classes such as outerplanar graphs. A research direction would be to try to express weaker constraints allowing a dynamic programming framework based on decomposition into independent subcomponents to apply for other classes. Outerplanar graphs are a natural candidate, since the geodetic set problem is polynomial for them but NP-hard for planar graphs.

Finally, other graph products than Cartesian (studied in \Cref{sec-generalAndProducts}) could be explored, such as the tensor, strong, lexicographical, modular or rooted products. However, those products may be more involved: for example, as seen in the same section, the outcome of the tensor product of two graphs is not necessarily a function of the outcomes of the two graphs. Some further conditions on the structure of the graphs or their product may be needed to obtain similar results.

\subsubsection*{Acknowledgements.} The authors would like to thank Sagnik Sen, who initiated this research project and for his advice and tremendous help; as well as the universities' libraries, without which the seminal paper~\cite{buckley1985closed} wouldn't be obtainable.

\bibliographystyle{splncs04}
\bibliography{newReferences}

\end{document}